 \documentclass[final,3p,times]{elsarticle}
\usepackage{amsmath,amssymb}
\usepackage{mathrsfs}
\usepackage{turnstile}
\usepackage{hyperref}
\usepackage{float}
\usepackage[normalem]{ulem}
\newtheorem{fact}{Fact}[section]
\newtheorem{prop}[fact]{Proposition}
\newtheorem{defn}[fact]{Definition}
\newtheorem{thm}[fact]{Theorem}

\newtheorem{lem}[fact]{Lemma}
\newtheorem{corr}[fact]{Corollary}
\newdefinition{rmk}[fact]{Remark}
\newdefinition{concl}[fact]{Conclusion}
\newproof{pf}{Proof}

\newtheorem{claim2}[fact]{\textbf{Claim}}

\usepackage{proof}
\usepackage{ulem}
\usepackage{accents}

\DeclareMathSymbol{\mhyph}{\mathalpha}{operators}{`-}

\journal{ }

\begin{document}

\begin{frontmatter}



\title{Theories of Frege structure equivalent to Feferman's system $\mathsf{T}_0$}


\author{Daichi Hayashi}

%\affiliation{organization={Hokkaido University},%Department and Organization
%           addressline={}, 

\begin{abstract}
Feferman \cite{feferman1975language} defines an impredicative system $\mathsf{T}_0$ of explicit mathematics, which is proof-theoretically equivalent to the subsystem $\Delta^1_2 \mhyph \mathsf{CA} + \mathsf{BI}$ of second-order arithmetic.
In this paper, we propose several systems of Frege structure with the same proof-theoretic strength as $\mathsf{T}_0$.
To be precise, we first consider the Kripke--Feferman theory, which is one of the most famous truth theories, and we extend it by two kinds of induction principles inspired by \cite{jager2001universes}.
In addition, we give similar results for the system based on Aczel's original Frege structure \cite{aczel1980frege}.
Finally, we equip Cantini's supervaluation-style theory with the notion of universes, the strength of which was an open problem in \cite{kahle2001truth}.

%In particular, those are based on Kripke--Feferman truth, Aczel's original definition of Frege structure, and Cantini's supervaluation-style truth, respectively. 
\end{abstract}




\end{frontmatter}


\tableofcontents
%\clearpage
\section{Introduction}\label{sec:introduction}
Aczel~\cite{aczel1980frege} introduced the framework of \emph{Frege structure}, which is essentially a model of lambda calculus augmented with the notions of \emph{truth} and \emph{proposition}, to analyse Russell's paradox in Frege's Grundgesetze.
While Aczel's study was model-theoretic, Frege structure can be seen as an axiomatic theory of truth formulated over an applicative theory (see Section~\ref{sec:preliminaries}). In addition, depending on what kinds of truth and proposition are assumed, various theories of Frege structure have been proposed (for an overview, see, e.g.,~\cite{cantini1996logical}).
As another applicative framework, Feferman formulated systems of \emph{explicit mathematics}~\cite{feferman1975language}, in which \emph{types}, second-order objects for sets, are generated by individual terms, called \emph{names}.

Although these two frameworks are different at a glance, there are, as Aczel anticipated \cite[p.~34]{aczel1980frege}, various technical correspondences between them.
In particular, for various systems of explicit mathematics, we can find proof-theoretically equivalent theories of Frege structure, as displayed in Table~\ref{table:applicative_systems},
where $\mathsf{EM}$ is a system of explicit mathematics (see Definition~\ref{def:EM}); $\mathsf{EMU}$ is an extension of $\mathsf{EM}$ by \emph{universes} (cf. \cite{strahm1999first}); $\mathsf{NEM}$ is an extension of $\mathsf{EM}$ by \emph{name induction} (cf. \cite{kahle2000theory}); $\mathsf{T}_0$ is an extension of $\mathsf{EM}$ by \emph{inductive generations} (see Definition~\ref{def:T_0}); 
%$\mathsf{NAI}$ is an extension of $\mathsf{EM}$ by $name \ induction$; $\mathsf{LUN}$ is an extension of $\mathsf{EM}$ by $least \ universes$;
$\mathsf{KF}$ and $\mathsf{PT}$ are theories of Frege structure based on Kripke-Feferman logic and Aczel-Feferman logic, respectively (see Definition~\ref{def:KF} and \ref{def:PT}); $\mathsf{KFU}$ and $\mathsf{PTU}$ are respectively extensions of $\mathsf{KF}$ and $\mathsf{PT}$ by \emph{universes} (see Definition~\ref{def:KFU} and \ref{def:PTU}); and $\mathsf{VF}$ is a theory of Frege structure based on supervaluation logic (see Definition~\ref{def:VF}). The table also contains the corresponding subsystems of second-order arithmetic and their proof-theoretic ordinals. 
The system $\Sigma^1_1 \mhyph \mathsf{AC}$ has the schema $\Sigma^1_1$ axiom of choice;
$\mathsf{ATR} + (\Sigma^{1}_{1} \mhyph \mathsf{DC})$
consists of the arithmetical transfinite recursion $\mathsf{ATR}$ with the $\Sigma^1_1$ dependent choice;
$\Pi^1_1 \mhyph \mathsf{CA}^{\mhyph\mhyph}_0$ is the parameter-free $\Pi^1_1$ comprehension schema;
and $\Delta^{1}_{2} \mathsf{\mhyph CA}+\mathsf{BI}$ is the $\Delta^1_2$ comprehension schema with the bar induction.
For details of their proof-theoretic ordinals, see, e.g., \cite{jager1999proof, pohlers1998subsystems}). 

\begin{table}[H]\label{table:applicative_systems}
\centering
\caption{Applicative theories}
\begin{tabular}{|c|c|c|c|}
\hline
Ordinal strength & Explicit mathematics            & Frege structure             & Second-order arithmetic                                    \\ \hline
$\psi_{\Omega}(\varepsilon_{I+1})$ & $\mathsf{T}_{0}$          & $?$                & $\Delta^{1}_{2} \mathsf{\mhyph CA}+\mathsf{BI}$                  \\ \hline
$\psi_{\Omega}(\varepsilon_{\Omega+1})$ & $\mathsf{NEM}$ & $\mathsf{VF}$             & $\Pi^{1}_{1} \mhyph \mathsf{CA}^{\mhyph \mhyph}_0$ \\ \hline
          
$\varphi 1\varepsilon_00$ & $\mathsf{EMU}$ & $\mathsf{KFU}, \mathsf{PTU}$  & $\mathsf{ATR} + (\Sigma^{1}_{1} \mhyph \mathsf{DC})$ \\ \hline

$\varphi\varepsilon_00$ & $\mathsf{EM}$              & $\mathsf{KF},\mathsf{PT}$  &     $\Sigma^{1}_{1} \mhyph \mathsf{AC}$             \\ \hline
\end{tabular}
\end{table}

As the table shows, the correspondence between explicit mathematics and Frege structure has so far been obtained only up to the strength of $(\Pi^{1}_{1} \mhyph \mathsf{CA})_0^{\mhyph \mhyph}$.
Therefore, this paper aims to provide well-motivated theories of Frege structure proof-theoretically equivalent to $\mathsf{T}_0$.

%When seeing Frege structure as theory of truth, this task also has importance from the foundational viewpoint.
This task also has importance from the foundational viewpoint if Frege structure is seen as a theory of truth.
In axiomatic theory of truth, it has been one of the central tasks to obtain stronger truth theories. For one thing, Halbach \cite{halbach2000truth} argues that an expressively strong truth theory can, to some extent, reduce ontological assumptions on sets to semantic assumptions. 
Thus, perhaps such a theory, if it is well motivated, can take the  place of set theory as a foundation for a large part of mathematics. 
As far as the author knows, Fujimoto's system $\mathsf{Aut(VF)}$ \cite{fujimoto_2011}, formulated over Peano arithmetic ($\mathsf{PA}$), is so far the strongest among well-motivated truth theories.\footnote{However, the author also remarks that Cantini gave a recursion-theoretically motivated system of Frege structure which is at least as strong as $\mathsf{T}_0$ \cite[p.~253]{cantini1996logical}. Over Peano arithmetic, Schindler \cite{schindler2015disquotational} formulated a disquotational truth theory of the same consistency strength as the full second-order arithmetic.} 
Since the strength of $\mathsf{Aut(VF)}$ lies strictly between $\Delta^1_2 \mhyph \mathsf{CA}$ and $\Delta^1_2 \mhyph \mathsf{CA} + \mathsf{BI}$, the author believes that the theories proposed in this paper break the record, though our base theory is not $\mathsf{PA}$.
Moreover, since $\mathsf{T}_0$ is expressively rich enough to interpret various set theories (cf.~\cite{kentaro2015new,tupailo2001realization,tupailo2003realization}), 
we can expect our theories of Frege structure to contribute to Halbach's programme.

The structure of this paper is as follows: In Section~\ref{sec:preliminaries}, we define the total applicative theory $\mathsf{TON}$ as the common base theory of explicit mathematics and Frege structure. We also define the theories of explicit mathematics, $\mathsf{EM}$ and its extension $\mathsf{T}_0$.
In Section~\ref{sec:FS_by_SK}, following Kahle's formulation \cite{kahle2001truth,kahle2003universes}, we introduce Kripke-Feferman-style theories of Frege structure, $\mathsf{KF}$ and its extension $\mathsf{KFU}$.
In Section~\ref{sec:extension_by_PI}, we expand $\mathsf{KFU}$ by the \emph{proposition induction} schema, to obtain a theory $\mathsf{KFUPI}$ that is as strong as $\mathsf{T}_0$.
In Section~\ref{sec:FS_by_AF}, we consider Aczel-Feferman-style theories of Frege structure, $\mathsf{PT}$ and its extension $\mathsf{PTU}$. In conclusion, we can obtain a theory $\mathsf{PTUPI}$ that has the same strength as $\mathsf{T}_0$.
In Section~\ref{sec:extension_by_LU}, we formulate an alternative principle \emph{least universes} schema, and then we show that this also gives the strength of $\mathsf{T}_0$ to both $\mathsf{KFU}$ and $\mathsf{PTU}$.
In Section~\ref{sec:FS_by_supervaluation}, we introduce a supervaluation-style Frege structure essentially based on Kahle's formulation \cite{kahle2001truth}.
Then, following Kahle's suggestion \cite[p.~124]{kahle2001truth},
we extend $\mathsf{VF}$ by universes and prove that the resulting theory $\mathsf{VFU}$ is proof-theoretically equivalent to $\mathsf{T}_0$.
Therefore, one of Kahle's open questions is solved.

\section{Technical preliminaries}\label{sec:preliminaries} \ 
This section defines Feferman's system $\mathsf{T}_0$ of explicit mathematics.
We first define the total applicative theory $\mathsf{TON}$,
which, in this paper, is used as the common base theory of explicit mathematics and Frege structure.
\subsection{Total applicative theory}\label{subsec:TON}
The first-order language $\mathcal{L}$ for the total operational theory $(\mathsf{TON})$ (see, e.g. \cite{kahle2003universes}) consists of
the standard logical symbols, individual variables 
($x,y,z,x_{0}, x_{1}, \dots, a,b,c, \dots,  f,g,h$),
individual constants $\mathsf{k}, \mathsf{s}$ (combinators), $\mathsf{p}, \mathsf{p_{0},p_{1}}$ (pairing and projections), $\mathsf{0}$ (zero),  $\mathsf{s}_\mathsf{N}$ (successor), $\mathsf{p}_\mathsf{N}$ (predecessor), $\mathsf{d}_\mathsf{N}$ (definition by numerical cases),
 a binary function symbol $\mathsf{App}(x,y)$ (application), and a unary predicate $\mathrm{N}(x)$ (natural numbers).
Formulae of $\mathcal{L}$ are constructed from atomic formulae by the logical symbols $\neg A$, $A \land B$, $A \to B$, and $\forall x. \ A$.
We assume that $\lor$ and $\exists$ are defined in a standard manner, whereas $\to$ is given as a primitive symbol.
 The meaning of each symbol is clear from its defining axioms below.
% For example, $\mathsf{d}_\mathrm{N}$ returns different objects according to whether a given equation holds or not. 
 
We shall use the following abbreviations in this paper:
\begin{enumerate}
\item $ab := (ab):= \mathsf{App}(a,b)$.
\item $a_{1}a_{2}a_{3} \dots a_{n} : \equiv ((\cdots((a_{1}a_{2})a_{3}\cdots )a_{n})$.
\item $(x,y) : = \mathsf{p}xy.$ \item $(x)_{i} : = \mathsf{p}_{i}x$ for $i \in \{ 0,1 \}.$
\item $s \neq t := \neg(s=t) $.
\item $\forall x_0, x_1, \dots x_n. \ A :\equiv \forall x_0. \forall x_1. \dots \forall x_n. \ A$.
\end{enumerate}

\begin{defn}
The $\mathcal{L}$-theory $\mathsf{TON}$ consists of the following axioms:
\begin{itemize}
\item $\mathsf{k}ab = a$,
\item $\mathsf{s}abc = ac(bc)$,
\item $(a,b)_0=a \land (a,b)_1 = b$,
\item $\mathrm{N}(\mathsf{0}) \land \forall x. \ \mathrm{N}(x) \to \mathrm{N}(\mathsf{s_{N}}x)$,
\item $\forall x. \ \mathrm{N}(x) \to \mathsf{s}_\mathsf{N}x \neq \mathrm{0} \land \mathsf{p}_\mathsf{N}(\mathsf{s}_\mathsf{N}x) = x$,
\item $\forall x. \ \mathrm{N}(x) \land x \neq \mathsf{0} \to \mathrm{N}(\mathsf{p}_\mathsf{N}x) \land \mathsf{s}_\mathsf{N}(\mathsf{p}_\mathsf{N}x) = x$,
\item $\mathrm{N}(a) \land \mathrm{N}(b) \land a = b \to \mathsf{d}_\mathsf{N}uvab = u$,
\item $\mathrm{N}(a) \land \mathrm{N}(b) \land a \neq b \to \mathsf{d}_\mathsf{N}uvab = v$,
\item $A(\mathsf{0}) \land [\forall x. \ \mathrm{N}(x) \land A(x) \to A(\mathsf{s}_\mathsf{N}x)] \to \forall x. \ \mathrm{N}(x) \to A(x)$,
for every $\mathcal{L}$-formula $A$.
\end{itemize}
\end{defn}

In this paper, we will repeatedly use the following well-known fact
:
\begin{fact}[cf. \cite{cantini1996logical}]\label{daft:TON} \ 
\begin{description}
\item[$\lambda$-abstraction.] For each variable $x$ and an $\mathcal{L}$-term $t$, we can find a term $\lambda x.t$ such that
$\mathsf{TON} \vdash (\lambda x.t)x = t.$

\item[Recursion.] There is a term $\mathsf{rec}$ such that
$\mathsf{TON} \vdash \forall f. \ \mathsf{rec}f = f (\mathsf{rec}f).$

\end{description}
\end{fact}

\subsection{Explicit mathematics}\label{subsec:EM}
We define systems of explicit mathematics over $\mathsf{TON}$.
Usually, theories of explicit mathematics are defined in second-order language. For simplicity, however, we shall formulate them as first-order ones, similar to \cite{hayashi1995new}.

The first-order language $\mathcal{L}_{EM}$ of explicit mathematics is an extension of $\mathcal{L}$ with 
two predicates: a unary predicate symbol $\mathrm{R}(x),$ meaning that $x$ represents a set,
and a binary predicate symbol $x \in y,$ meaning that $x$ is contained in $y$.
We define $\forall x \in y. \ A(x)$ to be $\forall x. \ x \in y \to A(x)$.
In addition, $\mathcal{L}_{EM}$ has individual constant symbols, called \emph{generators}: $\mathsf{int}$ (intersection), $\mathsf{j}$ (join), $\mathsf{nat}$ (natural numbers), $\mathsf{id}$ (identity), $\mathsf{dom}$ (domain), $\mathsf{inv}$ (inversion), and  $\mathsf{i}$ (inductive generations).
%Each generator is intended to construct a set from existing sets.
%For instance, given a set and a relation, $\mathsf{i}$ generates  the accessible part of the set with respect to the relation.

\begin{defn}\label{def:EM}
The $\mathcal{L}_{EM}$-theory $\mathsf{EM}$ consists of the following axioms:
\begin{description}
\item[Natural numbers.] $\mathrm{R}(\mathsf{nat}) \land \forall x. \ x \in \mathsf{nat} \leftrightarrow \mathrm{N}(x).$

\item[Identity.] $\mathrm{R}(\mathsf{id}) \land \forall x. \ x \in \mathsf{id} \leftrightarrow \exists y. \ x = (y,y).$

\item[Complements.] $\mathrm{R}(x) \to \mathrm{R}(\mathsf{co}(x)) \land \forall y. \ y \in \mathsf{co}(x) \leftrightarrow y \notin x$.

\item[Intersections.] $\mathrm{R}(x) \land \mathrm{R}(y) \to$
$ \mathrm{R}(\mathsf{int}(x,y)) \land \forall z. \ z \in \mathsf{int} (x,y) \leftrightarrow z \in x \land z \in y$.

\item[Domains.] $\mathrm{R}(a) \to \mathrm{R}(\mathsf{dom}(a)) \land \forall x. \ x \in \mathsf{dom}(a) \leftrightarrow \exists y. \ (x,y) \in a .$

\item[Inverse images.] $\mathrm{R}(a) \to \mathrm{R} (\mathsf{inv}(a,f)) \land \forall x. \ x \in \mathsf{inv}(a,f) \leftrightarrow fx \in a.$

\item[Joins.] $\mathrm{R}(x) \land [\forall y\in x. \ \mathrm{R}(fy)] \to \mathrm{R}(\mathsf{j} (x,f)) \land \Sigma(x,f,\mathsf{j}(x,f)),$
where
\[
\Sigma(x,f,y) : \equiv \forall u. \ u \in y \leftrightarrow \exists v, w. \ u = (v,w) \land v \in x \land w \in fv.
\]

\end{description} 
\end{defn}
The above axioms explain how a new set is generated from existing ones.
For instance, the join axiom says that given a set $x$ and a function $f$ whose domain is $x$, there exists the disjoint union $\mathsf{j}(x,f)$ of the range of $f$.

\begin{fact}[cf. \cite{beeson1985foundations}]\label{fact:strength_of_EM}
$\mathsf{EM}$ is proof-theoretically equivalent to $\Sigma^1_1 \mhyph \mathsf{AC}$.%$\widehat{\mathsf{ID}}_1.$
\end{fact}

\subsection{Universes and inductive generations}\label{subsec:U_and_IG}

A \emph{universe} in explicit mathematics is a set that is closed under the name-generating operations in Definition~\ref{def:EM}.
More formally, the fact that $a$ is a universe is expressed by the formula: 
\[
\mathrm{U}(a) : \equiv [\forall b. \ \mathcal{C}(a,b) \to b \in a] \to \forall b \in a. \ \mathrm{R}(b),
\] 
where $\mathcal{C}(a,b)$ is the disjunction of the following:
\begin{enumerate}
\item $a = \mathsf{nat} \lor a = \mathsf{id}$,
\item $\exists x. \ b = \mathsf{co}(x) \land x \in a$,
\item $\exists x,y. \ b = \mathsf{int}(x,y) \land x \in a \land y \in a$,
\item $\exists x. \ b = \mathsf{dom}(x) \land x \in a$,
\item $\exists f,x. \ b= \mathsf{inv}(x,f) \land x \in a$,
\item $\exists f,x. \ b= \mathsf{join}(x,f) \land x \in a \land \forall y \in x. fy \in a$.
\end{enumerate}

Of course, we require an additional axiom to assure the existence of universes.
Here, we introduce the \emph{limit axiom} \cite{strahm1999first}.
For a new constant symbol $\mathsf{l}$, the limit axiom is the following:
\[
\mathrm{R}(x) \to \mathrm{R}(\mathsf{l}x) \land x \in \mathsf{l}x. 
\]

Then, the system $\mathsf{EMU}$ is defined as $\mathsf{EM}$ equipped with the limit axiom.
From Strahm's proof \cite{strahm1999first}, it follows that $\mathsf{EMU}$ is proof-theoretically equivalent to %$\mathsf{ATR} + (\Sigma^{1}_{1} \mhyph \mathsf{DC})$ 
$\widehat{\mathsf{ID}}_{<\varepsilon_0}$, the arithmetical fixed-point theory iterated up to $\varepsilon_0$.
While the proof-theoretic strength of $\mathsf{EMU}$ is beyond the range of predicativity, that is, its proof-theoretic ordinal is larger than the Feferman--Sch\"{u}tte ordinal $\Gamma_0$,
it is still weaker than impredicative theories, such as the theory $\mathsf{ID}_1$ of arithmetical inductive definitions (for the definition, see \cite{pohlers2008proof}).
Therefore, the strength of $\mathsf{EMU}$ is called \emph{metapredicative} \cite{strahm1999first}.

Feferman \cite{feferman1975language} formulated a highly strong principle for explicit mathematics called \emph{inductive generation}:
\begin{align}
& \mathrm{R}(a) \land \mathrm{R}(b) \to \mathrm{R}(\mathsf{i}(a,b)) \land \mathsf{Closed}(a,b,\mathsf{i}(a,b)). \tag{$\mathsf{IG}.1$} \\
& \mathrm{R}(a) \land \mathrm{R}(b) \land \mathsf{Closed}(a,b,A) \to \forall x \in \mathsf{i}(a,b). \ A(x), \tag{$\mathsf{IG}.2$}
\end{align}

where $A$ is any formula, and
\begin{itemize}
\item $y <_{b} x : \equiv (y,x) \in b$;
%\item $\forall x \in y. \ A(x,y) :\equiv \forall x. \ x \in y \to A(x,y)$;
\item $\mathsf{Closed}(a,b,A(\bullet)) :\equiv \forall x \in a. \ [\forall y \in a. \ y <_{b} x \to A(y)] \to A(x)$.
\end{itemize}

\begin{defn}\label{def:T_0}
The $\mathcal{L}_{EM}$-theory $\mathsf{T}_0$ consists of $\mathsf{EM}$ with $(\mathsf{IG}.1)$ and $(\mathsf{IG}.2)$.
\end{defn}

The proof-theoretic strength of $\mathsf{T}_0$ is far beyond $\mathsf{EMU}.$
\begin{fact}[\cite{jager1983well}]\label{strength_of_T0}
$\mathsf{T}_0$ is proof-theoretically equivalent to $\Delta^1_2 \mhyph \mathsf{CA} + \mathsf{BI}$.% and to $\mathsf{KPi}$ (see Definition~\ref{def:KPi}).
\end{fact}

%%%%%%%%%%%%%%%%%%%%%%%%%%%%%%%%%%%%%%%%%%%%%%%%%%%%%%%%%%%%%%%%%%%%%%%%%%

\section{Frege structure by Strong Kleene schema}\label{sec:FS_by_SK}
In this section, we introduce the Kripke--Feferman theory $\mathsf{KF}$ of Frege structure and its extension by universes. 
We first define the base language $\mathcal{L}_{FS}$ over which our theories of Frege structure are formulated.

The language $\mathcal{L}_{FS}$ is $\mathcal{L}$ augmented with the following symbols:
\begin{itemize}
\item individual constants $\dot{=}, \dot{\neg}, \dot{\land}, \dot{\to}, \dot{\forall}, \dot{\mathsf{N}},  \mathsf{l}$;
\item unary predicates $\mathrm{T}(x), \mathrm{U}(x).$
\end{itemize}
Here, $\mathrm{T}(x)$ is intended to mean the truth predicate; $\mathrm{U}(x)$ means that $x$ is a \emph{universe} (see below for more details); and the constant $\mathsf{l}$ is used to generate universes, similar to the one in explicit mathematics.
% of the form $\mathrm{T}(lf(x))$.
 The other individual constants are used as sentence-constructing operators. For example, the term $x \dot{\land} y : = \dot{\land}(x,y)$ informally denotes (the code of) a conjunctive sentence consisting of $x$ and $y$. 
We can similarly understand the terms $\dot{\neg}x$, $\dot{\mathsf{N}}x$, and $\dot{\forall}f$.
We also use the notations
$(x \dot{=} y) := (\dot{=} (x,y))$ and
$(x \dot{\to} y) := (\dot{\to}(x,y))$, whose informal meaning should be clear.

\subsection{System $\mathsf{KF}$}\label{subsec:system_KF}
The theory $\mathsf{KF}$, as a theory of truth, was semantically introduced by Kripke \cite{kripke1976outline}; Feferman \cite{feferman1991reflecting} then gave its first-order axiomatisation. Cantini \cite{cantini1993extending} later formulated $\mathsf{KF}$ as a theory of Frege structure.
In $\mathsf{KF}$, each sentence is monotonically evaluated based on the Strong Kleene evaluation, as displayed in the following truth table.
Note that the conditional $A \to B$ is definable as $\neg A \lor B$.
\begin{center}
\begin{tabular}{|c|c|c|c|}
\hline
$\neg$        &     \\ \hline
$\mathrm{T}$  & $\mathrm{F}$    \\ \hline
$\mathrm{U}$ & $\mathrm{U}$  \\ \hline
$\mathrm{F}$  & $\mathrm{T}$    \\ \hline
\end{tabular} \ \ \ \ 
\begin{tabular}{|c|c|c|c|}
\hline
$\lor$        & $\mathrm{T}$  & $\mathrm{U}$ & $\mathrm{F}$  \\ \hline
$\mathrm{T}$  & $\mathrm{T}$  & $\mathrm{T}$ & $\mathrm{T}$  \\ \hline
$\mathrm{U}$ & $\mathrm{T}$ & $\mathrm{U}$ & $\mathrm{U}$ \\ \hline
$\mathrm{F}$  & $\mathrm{T}$  & $\mathrm{U}$ & $\mathrm{F}$  \\ \hline
\end{tabular} \ \ \ \ 
\begin{tabular}{|c|c|c|c|}
\hline
$\land$        & $\mathrm{T}$  & $\mathrm{U}$ & $\mathrm{F}$  \\ \hline
$\mathrm{T}$  & $\mathrm{T}$  & $\mathrm{U}$ & $\mathrm{F}$  \\ \hline
$\mathrm{U}$ & $\mathrm{U}$ & $\mathrm{U}$ & $\mathrm{F}$ \\ \hline
$\mathrm{F}$  & $\mathrm{F}$  & $\mathrm{F}$ & $\mathrm{F}$  \\ \hline
\end{tabular} \ \ \ \ 
%\end{table}
%\begin{table}[]
\begin{tabular}{|c|c|c|c|}
\hline
$\to$        & $\mathrm{T}$  & $\mathrm{U}$ & $\mathrm{F}$  \\ \hline
$\mathrm{T}$  & $\mathrm{T}$  & $\mathrm{U}$ & $\mathrm{F}$  \\ \hline
$\mathrm{U}$ & $\mathrm{T}$ & $\mathrm{U}$ & $\mathrm{U}$ \\ \hline
$\mathrm{F}$  & $\mathrm{T}$  & $\mathrm{T}$  & $\mathrm{T}$  \\ \hline
\end{tabular}
\end{center}

The formal system $\mathsf{KF}$, as an $\mathcal{L}_{FS}$-theory, is defined straightforwardly from the above table.
\begin{defn}[cf. \cite{cantini1996logical, cantini2016truth, kahle2001truth, kahle2003universes}]\label{def:KF}
The $\mathcal{L}_{FS}$-theory $\mathsf{KF}$ has $\mathsf{TON}$ with the full induction schema for $\mathcal{L}_{FS}$ and (the universal closure of) the following axioms:
\begin{description}
\item[Compositional axioms.] \

\begin{description}
\item[$(\mathsf{K}_=)$] $\mathrm{T}(x \dot{=} y) \leftrightarrow x = y$

\item[$(\mathsf{K}_{\neq})$] $\mathrm{T}(\dot{\neg}(x \dot{=} y)) \leftrightarrow x \neq y$

\item[$(\mathsf{K}_{\mathrm{N}})$] $\mathrm{T}(\dot{\mathsf{N}}x) \leftrightarrow \mathrm{N}(x)$

\item[$(\mathsf{K}_{\neg \mathrm{N}})$] $\mathrm{T}(\dot{\neg}(\dot{\mathsf{N}}x)) \leftrightarrow \neg \mathrm{N}(x)$

\item[$(\mathsf{K}_{\neg\neg})$] $\mathrm{T} (\dot{\neg} (\dot{\neg} x)) \leftrightarrow \mathrm{T}(x)$

\item[$(\mathsf{K}_{\land})$] $\mathrm{T} (x \dot{\land} y) \leftrightarrow \mathrm{T}(x) \land \mathrm{T} (y)$

\item[$(\mathsf{K}_{\neg\land})$] $\mathrm{T} (\dot{\neg} (x \dot{\land} y)) \leftrightarrow \mathrm{T}(\dot{\neg}x) \lor \mathrm{T} (\dot{\neg}y)$

\item[$(\mathsf{K}_{\to})$] $\mathrm{T} (x \dot{\to} y) \leftrightarrow \mathrm{T}(\dot{\neg}x) \lor \mathrm{T}(y)$

\item[$(\mathsf{K}_{\neg\to})$] $\mathrm{T}(\dot{\neg} (x \dot{\to} y)) \leftrightarrow \mathrm{T}(x) \land \mathrm{T}(\dot{\neg}y)$

\item[$(\mathsf{K}_{\forall})$] $\mathrm{T}(\dot{\forall} f) \leftrightarrow \forall x. \ \mathrm{T}(fx)$

\item[$(\mathsf{K}_{\neg\forall})$] $\mathrm{T}(\dot{\neg}(\dot{\forall} f)) \leftrightarrow \exists x. \ \mathrm{T}(\dot{\neg}(fx))$
\end{description}
\item[Consistency.] \

\begin{description}
\item[$(\mathsf{T} \mhyph \mathsf{Cons})$] $\neg [\mathrm{T}(x) \land \mathrm{T}(\dot{\neg}x)]$
\end{description}
\end{description}
%The axioms for $\land$ and $\exists$ are similarly defined.
\end{defn}

The system $\mathsf{KF}$ is proof-theoretically stronger than $\mathsf{TON}$.
In particular, $\mathsf{KF}$ can relatively interpret $\mathsf{EM}$.
To see this, we define a translation $': \mathcal{L}_{EM} \to \mathcal{L}_{FS}$.
First, the vocabularies of $\mathcal{L}$ are unchanged.
Second, we define $(x \in y)'$ to be $\mathrm{T}(yx)$,
which can be read as `$y$ is true at $x$.' 
As for the translation of $\mathrm{R}(x)$, 
we define predicates $\mathrm{C}(f)$ and $\mathrm{P}(x)$, where $\mathrm{P}(x)$ means that $x$ is a \emph{proposition}, that is, $x$ is determined to be true or false; and $\mathrm{C}(f)$ means that $f$ is a \emph{class} (propositional function), that is, $fx$ is a proposition for every object $x$.
\[ 
\mathrm{C}(f) := \forall x. \ \mathrm{P}(fx) := \forall x. \ \mathrm{T}(fx) \lor \mathrm{T} (\dot{\neg}(fx)). 
\]
Then, we let $\mathrm{R}'(x) :\equiv \mathrm{C}(x)$.
Finally, as is remarked in \cite[p.~59]{cantini1996logical}, for each generator of $\mathcal{L}_{EM}$, we can find a term such that the translation of the defining axiom of $c$ in $\mathsf{EM}$ is derivable in $\mathsf{KF}$ (see also the proof of Lemma~\ref{lem:interpretation_of_join_VFU}).
To summarise, the following is obtained:
\begin{fact}[\cite{cantini1996logical}]
\label{fact:strength_of_KF}
%$\mathsf{KF}$ and $\mathsf{EM}$ have the same $\mathcal{L}$-theorems.
For each $\mathcal{L}_{EM}$-formula $A$,
if $\mathsf{EM} \vdash A$, then $\mathsf{KF} \vdash A'.$
\end{fact}

We also remark that $\mathsf{EM}$ and $\mathsf{KF}$ are proof-theoretically equivalent (for the proof, see, e.g., \cite[Section~57]{cantini1996logical}).

\subsection{Universes for $\mathsf{KF}$}\label{subsec:universes_for_KF}
In order to strengthen a given truth theory, iterating truth predicates is effective in many cases (cf. \cite{feferman1991reflecting,jager1999proof, fujimoto_2011}).
In the framework of Frege structure, an analogous idea is realisable by using the notion of \emph{universes} (cf. \cite{cantini1996logical,kahle2003universes}). 
Following \cite{kahle2003universes}, %we introduce several notations to express universes. Firstly, we define predicates $\mathrm{C}(f)$ and $\mathrm{P}(x)$, where $\mathrm{P}(x)$ means that $x$ is a \emph{proposition}, that is, $x$ is determined to be true or false; $\mathrm{C}(f)$ means that $f$ is a \emph{class} (propositional function), that is, $fx$ is a proposition for every object $x$.
%\[ 
%\mathrm{C}(f) := \forall x. \mathrm{P}(fx) := \forall x. \mathrm{T}(fx) \lor \mathrm{T} (\dot{\neg}(fx)). 
%\]
%Next, 
we let $f \sqsubset g $ ($g$ $\emph{reflects}$ on $f$) be the following formula:
\[ 
[\forall x. \ \mathrm{T}(fx) \to \mathrm{T}(g(fx))] \ \land \ [\forall x. \ \mathrm{T}(\dot{\neg}fx) \to \mathrm{T}(g(\dot{\neg}fx))]. \]
The predicate $f \sqsubset g$ informally means that both positive and negative facts about $f$ are expressed as positive statements within $g$.
%The axiom $(\mathsf{Lim})$ thus generates a truth predicate $\mathsf{l}f$ for the language $\mathcal{L}$ with the class $f$ as a predicate.
 
%Similarly, for each formula $B(x)$, Let $f \sqsubset B $ be the following:
%\[
%[\forall x. \mathrm{T}(fx) \to B(fx)] \ \land \ [\forall x. \mathrm{T}(\dot{\neg}fx) \to B(\dot{\neg}fx)].
%\]

\begin{defn}[cf. \cite{kahle2001truth, kahle2003universes}]\label{def:KFU}
The $\mathcal{L}_{FS}$-theory $\mathsf{KFU}$ has $\mathsf{TON}$ and (the universal closure of) the following axioms:
\begin{description}
\item[Compositional axioms in $\mathrm{U}$.] \
\begin{description}
\item[$(\mathsf{U}_=)$] $\mathrm{U}(u) \to \forall x, y. \ \mathrm{T}(u(x \dot{=} y)) \leftrightarrow x = y$

\item[$(\mathsf{U}_{\neq})$] $\mathrm{U}(u) \to \forall x, y. \ \mathrm{T}(u(\dot{\neg}(x \dot{=} y))) \leftrightarrow x \neq y$

\item[$(\mathsf{U}_{\mathrm{N}})$] $\mathrm{U}(u) \to \forall x. \ \mathrm{T}(u(\dot{\mathsf{N}}x)) \leftrightarrow \mathrm{N}(x)$

\item[$(\mathsf{U}_{\neg\mathrm{N}})$] $\mathrm{U}(u) \to \forall x. \ \mathrm{T}(u(\dot{\neg}(\dot{\mathsf{N}}x))) \leftrightarrow \neg \mathrm{N}(x)$

\item[$(\mathsf{U}_{\neg\neg})$] $\mathrm{U}(u) \to \forall x. \ \mathrm{T} (u(\dot{\neg} (\dot{\neg} x))) \leftrightarrow \mathrm{T}(ux)$

\item[$(\mathsf{U}_{\land})$] $\mathrm{U}(u) \to \forall x, y. \ \mathrm{T} (u(x \dot{\land} y)) \leftrightarrow \mathrm{T}(ux) \land \mathrm{T} (uy)$

\item[$(\mathsf{U}_{\neg\land})$] $\mathrm{U}(u) \to \forall x, y. \ \mathrm{T} (u(\dot{\neg} (x \dot{\land} y))) \leftrightarrow \mathrm{T}(u(\dot{\neg}x)) \lor \mathrm{T} (u(\dot{\neg}y))$

\item[$(\mathsf{U}_{\to})$] $\mathrm{U}(u) \to \forall x, y. \ \mathrm{T} (u(x \dot{\to} y)) \leftrightarrow \mathrm{T}(u(\dot{\neg}x)) \lor \mathrm{T}(uy)$

\item[$(\mathsf{U}_{\neg\to})$] $\mathrm{U}(u) \to \forall x, y. \ \mathrm{T}(u(\dot{\neg} (x \dot{\to} y))) \leftrightarrow \mathrm{T}(ux) \land \mathrm{T}(u(\dot{\neg}y))$

\item[$(\mathsf{U}_{\forall})$] $\mathrm{U}(u) \to \forall f. \ \mathrm{T}(u(\dot{\forall} f)) \leftrightarrow \forall x. \ \mathrm{T}(u(fx))$

\item[$(\mathsf{U}_{\neg\forall})$] $\mathrm{U}(u) \to \forall f. \ \mathrm{T}(u(\dot{\neg}(\dot{\forall} f))) \leftrightarrow \exists x. \ \mathrm{T}(u(\dot{\neg}(fx)))$
\end{description}

\item[Consistency in $\mathrm{U}$.] \
\begin{description}
\item[$(\mathsf{U} \mhyph \mathsf{Cons})$] $\mathrm{U}(u) \to \forall x. \ \neg [\mathrm{T}(ux) \land \mathrm{T}(u(\dot{\neg}x))]$
\end{description}

\item[Structural properties of $\mathrm{U}$.] \
\begin{description}
\item[$(\mathsf{U} \mhyph \mathsf{Class})$] $\mathrm{U}(u) \to \mathrm{C}(u)$

\item[$(\mathsf{U} \mhyph \mathsf{True})$] $\mathrm{U}(u) \to \forall x. \ \mathrm{T}(ux) \to \mathrm{T}(x)$

\item[$(\mathsf{Lim})$] $\mathrm{C}(f) \to \mathrm{U}(\mathsf{l}f) \land f \sqsubset \mathsf{l}f$
\end{description}
\end{description}
\end{defn}

In the above definition, the predicate $\mathrm{U}(u)$ indicates that $u$ is a universe. 
The compositional axioms and the consistency axiom for $\mathrm{U}$ ensures that each universe satisfies the axioms of $\mathsf{KF}$; thus, universes work as quasi-truth predicates.
%This is expressed by relativising these axioms to elements of $\mathrm{U}$, as given in Definition~\ref{def:KFU}.
Of course, $\mathsf{KF}$-axioms relativised to universes do not have any logical strength unless we do not further postulate the existence of universes.
Thus, we added the axiom $(\mathsf{Lim})$ in Definition~\ref{def:KFU}, which generates an infinite series of universes:
$u_0 \sqsubset u_1 \sqsubset u_2 \sqsubset \cdots \sqsubset u_n \sqsubset \cdots$.
This roughly corresponds to a hierarchy of infinitely iterated truth predicates for $\mathsf{KF}$, whence the strength of $\mathsf{KFU}$ exceeds that of $\mathsf{KF}$.
Moreover, Kahle showed that $\mathsf{KFU}$ is a proper extension of $\mathsf{KF}$:

\begin{fact}[{\cite[Proposition~14]{kahle2003universes}}]\label{fact:KFU_implies_KF}
$\mathsf{KF}$ is a subtheory of $\mathsf{KFU}$.
\end{fact}

\begin{fact}[{\cite[Theorem~33]{kahle2003universes}}]\label{strength_of_KFU}
$\mathsf{KFU}$ and $\mathsf{EMU}$ have the same $\mathcal{L}$-theorems.
%In particular, they are proof-theoretically equivalent to 
%$\widehat{\mathsf{ID}}_{<\varepsilon_0}.$
\end{fact}

%%%%%%%%%%%%%%%%%%%%%%%%%%%%%%%%%%%%%%%%%%%%%%%%%%%%%%%%%%%%%%%%%%%%%%%%%%%%%%%%%%%%%%

\section{Extension by proposition induction}\label{sec:extension_by_PI}
While adding universes makes $\mathsf{KF}$ stronger, $\mathsf{KFU}$ falls within metapredicativity in strength,
so it is still much weaker than impredicative theories such as $\Pi^1_1 \mhyph \mathsf{CA}^{\mhyph \mhyph}_0$.
This section aims to propose a stronger principle that gives the same strength as $\mathsf{T}_0$.

For our first attempt, we shall borrow an idea from the system $\mathsf{NAI}$ of explicit mathematics in \cite{jager2001universes}. In the system $\mathsf{EMU}$,
names are inductively generated by several axioms, such as join and the limit axiom. The system $\mathsf{NAI}$ is then defined as an extension of $\mathsf{EMU}$ with induction principles requiring that the whole universe of names only contains those inductively generated by these axioms. 
Despite lacking inductive generations, $\mathsf{NAI}$ is shown to be proof-theoretically equivalent to $\mathsf{T}_0$ \cite[Conclusion~25]{jager2001universes}.
Since \emph{class} in Frege structure is an analogue of \emph{name} in explicit mathematics, we can expect to get a theory of Frege structure equivalent to $\mathsf{T}_0$ by considering an induction principle for classes or propositions.

\subsection{Proposition induction}\label{subsec:PI}
How, then, should we inductively characterise \emph{propositions} in $\mathsf{KF}$?
For example, let us consider the conjunction. In the Strong Kleene schema, $A \land B$ is a proposition when both $A$ and $B$ are propositions. In fact, we can derive $\mathrm{P}(x) \land \mathrm{P}(y) \to \mathrm{P}(x\dot{\land}y)$ in $\mathsf{KF}.$
On the other hand, once we establish the falsity of either $A$ or $B$,
the conjunction $A \land B$ is determined to be false regardless of the other conjunct.
This could be expressed, for example, as $([\mathrm{P}(x) \land \mathrm{T}(\dot{\neg}x)] \lor [\mathrm{P}(y) \land \mathrm{T}(\dot{\neg}y)]) \to \mathrm{P}(x \dot{\land} y),$ which is indeed provable in $\mathsf{KF}.$ So, in $\mathsf{KF},$ we can
characterise conjunctive propositions in two ways. 
Understanding the other logical symbols similarly, we can characterise propositions for $\mathsf{KF}$ as follows:

\begin{lem}\label{lem:prop_in_KF}
$\mathsf{KF}$ derives the following:
\begin{enumerate}
\item $\mathrm{P}(x \dot{=} y) \land \mathrm{P}(\dot{\mathsf{N}}x),$
\item $\mathrm{P}(\dot{\neg}x) \leftrightarrow \mathrm{P}(x),$
\item $\mathrm{P}(x \dot{\land} y) \leftrightarrow [\mathrm{P}(x) \land \mathrm{P}(y)] \lor [\mathrm{P}(x) \land \mathrm{T}(\dot{\neg}x)] \lor [\mathrm{P}(y) \land \mathrm{T}(\dot{\neg}y)],$
\item $\mathrm{P}(x \dot{\to} y) \leftrightarrow [\mathrm{P}(x) \land \mathrm{P}(y)] \lor [\mathrm{P}(x) \land \mathrm{T}(\dot{\neg}x)] \lor [\mathrm{P}(y) \land \mathrm{T}(y)],$
\item $\mathrm{P}(\dot{\forall}f) \leftrightarrow [\forall x. \ \mathrm{P}(fx)] \lor [\exists y. \ \mathrm{P}(fy) \land \mathrm{T}(\dot{\neg}(fy))].$
\end{enumerate}
\end{lem}

\begin{rmk}\label{rmk:relation_with_classical_determinateness}
This inductive characterisation of propositions is nearly the same as that of \emph{determinateness} in \cite{halbach2023classical}, where Halbach and Fujimoto advocate the definition by appealing to the function of truth as a generalising device 
\cite[p.~5]{halbach2023classical}.
We also remark that in $\mathsf{KF}$, it might also be natural to treat each proposition and its negation independently. For instance, in the Strong Kleene evaluation, $\neg(A\land B)$ is true when either $\neg A$ or $\neg B$ is true. So, whether $\neg (A \land B)$ is a proposition can be determined by $\neg A$ and $\neg B$, rather than by $A \land B$.
Formally speaking, $\mathsf{KF}$ derives: 
\[
\mathrm{P}(\dot{\neg}(x \dot{\land} y)) \leftrightarrow [\mathrm{P}(\dot{\neg}x) \land \mathrm{P}(\dot{\neg}y)] \lor [\mathrm{P}(\dot{\neg}x) \land \mathrm{T}(\dot{\neg}x)] \lor [\mathrm{P}(\dot{\neg}y) \land \mathrm{T}(\dot{\neg}y)].
\]
The same remark applies to the other negated propositions $\neg \neg A$, $\neg(A \to B)$, and $\neg \forall x. \ A$.
\end{rmk}

Finally, the axiom $(\mathsf{Lim})$ generates a new class $\mathsf{l}f$ from a given class $f$. In particular, we have:
\[
\mathsf{KFU} \vdash \mathrm{C}(f) \to \forall x. \ \mathrm{P}(\mathsf{l}fx).
\]

In summary, we can express the above construction of propositions in $\mathsf{KFU}$ by a single operator $\mathscr{A}$. For each $\mathcal{L}_{FS}$-formula $B$ and a free variable $x$, 
the formula $\mathscr{A} (B(\bullet), x)$
\footnote{Here, $\bullet$ means every occurrence of some fixed free variable (\emph{placeholder}).} is the disjunction of the following:
\begin{enumerate}
\item $\exists y, z. \ x = (y \dot{=} z) \land y = z$
%\item $\exists y, z. \ x = \dot{\neg} (y \dot{=} z) \land y \neq z$
\item $\exists y. \ x = (\dot{\mathsf{N}} y) \land \mathrm{N}(y)$
%\item $\exists y. \ x = \dot{\neg} (\dot{\mathsf{N}} y) \land \neg \mathrm{N}(y)$
\item $\exists y. \ x = \dot{\neg} y \land   B(y)$
\item $\exists y, z. \ x = y \dot{\land} z \land \{ [B(y) \land B(z)] \lor [B(y) \land \mathrm{T}(\dot{\neg}y)] \lor [B(z) \land \mathrm{T} (\dot{\neg}z)]\}$
\item $\exists y, z. \ x = y \dot{\to} z \land \{ [B(y) \land B (z)] \lor [B (y) \land \mathrm{T} (\dot{\neg}y)] \lor [B(z) \land \mathrm{T} (z)]\}$
\item $\exists g. \ x = \dot{\forall} g \land \{  [\forall y. \ B(gy) ] \lor \exists y. \ B(gy) \land \mathrm{T}(\dot{\neg}(gy))\}$
\item $\exists f, y. \ x=\mathsf{l}fy \land \forall z. \ B(fz)$
\end{enumerate}

Then, the proposition induction schema $(\mathsf{PI})$ is given by:
\begin{align}
[\forall x. \ \mathscr{A}(B(\bullet),x) \to B(x)] \to \forall x. \ \mathrm{P}(x) \to B(x), \tag{$\mathsf{PI}$}\label{PI}
\end{align}
for all $\mathcal{L}_{FS}$-formulas $B(x).$

The schema $(\mathsf{PI})$ assures that all propositions are obtained by the
inductive construction as explained above.

For the proof-theoretic purposes, we also need the following axioms.
\begin{defn}\label{def:UG}
The schema $(\mathsf{UG})$ consists of axioms asserting that each constant symbol in $\{ \dot{=}, \dot{\neg}, \dot{\land}, \dot{\to}, \dot{\forall}, \dot{\mathsf{N}}, \mathsf{l} \}$ can be uniquely decomposed. For example, $(\mathsf{UG})$ has the following:
\[
\dot{\land} \neq \dot{\neg},
\]
\[
\dot{\land}a = \dot{\land}b \to a=b.
\]
\end{defn}

With the help of $(\mathsf{UG})$, we can verify in $\mathsf{KFU}$ that $\mathrm{P}(x)$ is an $\mathscr{A}$-closure:
\begin{align}
\mathsf{KFU}+(\mathsf{UG}) \vdash \forall x. \ \mathscr{A}(\mathrm{P}(\bullet),x) \to \mathrm{P}(x). \notag
\end{align}

\begin{defn}\label{def:KFU+PI}
The $\mathcal{L}_{FS}$-theory $\mathsf{KFUPI}$ is $\mathsf{KFU}$ with the schemata $(\mathsf{UG})$ and $(\mathsf{PI})$.  
\end{defn}

\subsection{Lower bound of $\mathsf{KFUPI}$}\label{subsec:lower-bound_of_KFUPI}
The remainder of this section is devoted to the proof-theoretic analysis of $\mathsf{KFUPI}$.
In this subsection, we give a relative interpretation of $\mathsf{T}_0$ into $\mathsf{KFUPI}$ by extending the translation $'$ used in Fact~\ref{fact:strength_of_KF}. 
For each symbol except $\mathsf{i}$, we can employ the same interpretation as given in the proof of Fact~\ref{fact:strength_of_KF}. 
In particular, recall that the predicate $\mathrm{R}(f)$ is interpreted as $\mathrm{C}(f),$ namely the statement that $f$ is a \emph{class}.  
Moreover, $x \in y$ is interpreted as the formula $\mathrm{T}(yx).$
For readability, we sometimes refer to $\mathrm{T}(yx)$ as $x \dot{\in} y.$ Also, $\forall x \dot{\in}y. \ A(x)$ means $\forall x. \ x \dot{\in} y \to A(x)$.
Finally, we want to interpret the inductive generation $\mathsf{i}$
 as an appropriate term $\mathsf{acc}.$
In other words, we need to find a term $\mathsf{acc}$ such that $\mathsf{KFUPI}$ derives the translation of the axioms $(\mathsf{IG}.1)$ and $(\mathsf{IG}.2)$ in Definition~\ref{def:T_0}:
\begin{description}
\item[$(\mathsf{IG}.1)'$] $\mathrm{C}(a) \land \mathrm{C}(b) \to \mathrm{C}(\mathsf{acc}(a,b)) \land \mathsf{Closed}'(a,b, \mathrm{T}(\mathsf{acc}(a,b)(\bullet)))$;
\item[$(\mathsf{IG}.2)'$] $\mathrm{C}(a) \land \mathrm{C}(b) \land \mathsf{Closed}'(a,b,A(\bullet)) \to \forall x \dot{\in} \mathsf{acc}(a,b). \ A(x),$
\end{description}
where $A$ is any $\mathcal{L}_{FS}$-formula. We used the notations:
\begin{itemize}
\item $y \dot{<}_{b} x : \equiv (y,x) \dot{\in} b :\equiv \mathrm{T}(b(y,x))$;
\item $\mathsf{Closed}'(a,b,A(\bullet)) :\equiv \forall x \dot{\in} a. \ [\forall y \dot{\in} a. \ y \dot{<}_{b} x \to A(y)] \to A(x)$.
\end{itemize}

We essentially follow the proof of the lower bound of $\mathsf{LUN}$ 
(\cite[pp.~153-156]{jager2001universes}),
to derive the above formulas.
Roughly speaking, the interpretations of $(\mathsf{IG}.1)$ and $(\mathsf{IG}.2)$ are derivable by $(\mathsf{L}1)$ and $(\mathsf{L}2),$ respectively.
Let $\oplus$ be a term such that 
\[
\oplus(u, v) = \lambda z. [((z)_{0} \dot{=} 0 \dot{\land} u ((z)_{1})) \ \dot{\lor} \ (\dot{\neg}((z)_{0} \dot{=} 0) \dot{\land} v ((z)_{1}))],
\]
which represents the disjoint sum of $u$ and $v$.
By recursion, we can construct a term $t$ such that
\begin{center}
$t(u,v,w)  = \dot{\forall} y (\oplus(u, v) (0,y) \dot{\to} [\oplus(u, v) (1, (y,w)) \dot{\to} t(u,v,y)]).$
\end{center}
The term $t(u,v,w)$ informally says that in $u$, every $v$-predecessor $y$ of $w$ satisfies $t(u,v,y)$. Of course, $t(u,v,w)$ is not, in general, a proposition even if $u$ and $v$ are classes; hence, we cannot simply employ $t$ as the interpretation of $\mathsf{i}$.
Nevertheless, by attaching the limit constant $\mathsf{l}$ to $t$, we can treat $t$ as a proposition.
Thus, for the interpretation of $\mathsf{i},$ we further define a term $\mathsf{acc}$ such that 
\[
\mathsf{acc}(u,v)x = 
\bigl( \oplus(u, v) (0,x) \bigr) \dot{\land} \bigl( \dot{\forall} y (\oplus(u, v) (0,y) \dot{\to} [\oplus(u, v) (1,(y,x)) \dot{\to} (\mathsf{l}(\oplus(u, v)) (t(u,v,y))]) \bigr).
\]

By the following lemma, the term $\mathsf{acc}(a,b)$ is indeed a class whenever $a$ and $b$ are classes, so the first conjunct of $(\mathsf{IG}.1)'$ is derived. 
\begin{lem}\label{lem:acc_is_class}
$\mathsf{KFU} \vdash \mathrm{C}(a) \land \mathrm{C}(b) \to \mathrm{C}(\oplus(a,b)) \land \mathrm{C}(\mathsf{acc}(a,b)).$
\end{lem}

\begin{pf}
Assume $\mathrm{C}(a) \land \mathrm{C}(b).$ First, taking any $z$, we show that $\oplus(a,b)z$ is a proposition in order to show that $\oplus(a,b)$ is a class.
By the definition of $\oplus$, $\oplus(a,b)z$ is of the form: 
\begin{center}
$((z)_{0} \dot{=} 0 \dot{\land} a ((z)_{1})) \ \dot{\lor} \ (\dot{\neg}((z)_{0} \dot{=} 0) \dot{\land} b ((z)_{1}))$.
\end{center}
Here, $(z)_{0} \dot{=} 0$ and $\dot{\neg}((z)_{0} \dot{=} 0)$ are propositions by Lemma~\ref{lem:prop_in_KF}, and
$a ((z)_{1})$ and $b ((z)_{1})$ are also propositions by the assumption $\mathrm{C}(a) \land \mathrm{C}(b).$
Therefore, $\oplus(a, b) z$ is a proposition, again by Lemma~\ref{lem:prop_in_KF}. Since $z$ is arbitrary, $\oplus(a, b)$ is a class, and thus so is $\mathsf{l}(\oplus(a,b))$ by the axiom $(\mathsf{Lim})$.
Second, we deduce that $\mathsf{acc}(a,b)x$ is a proposition for any $x$.
Recall that $\mathsf{acc}(a,b)x$ is of the form 
\begin{center}
$\oplus(a, b) (0,x) \dot{\land} \dot{\forall} y (\oplus(a, b) (0,y) \dot{\to} [\oplus(a, b) (1,(y,x)) \dot{\to} (\mathsf{l}(\oplus(a, b)) (t(a,b,y)))])$.
\end{center}
Since we have shown that $\oplus(a, b)$ and $\mathsf{l}(\oplus(a, b))$ are classes,
Lemma~\ref{lem:prop_in_KF} implies that $\mathsf{acc}(a,b)x$ is a propositon, as required. \qed
\end{pf}

The following is the second half of $(\mathsf{IG}.1)'.$
\begin{lem}\label{lem:IG1_in_KFU}
Recall: 
$\mathsf{Closed}'(a,b,B(\bullet)):\equiv \forall x \dot{\in} a. \ [\forall y \dot{\in} a. \ y \dot{<}_{b} x \to A(y)] \to A(x).
$

Then, $\mathsf{KFU} \vdash \mathrm{C}(a) \land \mathrm{C}(b) \to \mathsf{Closed}'(a,b, \mathrm{T}(\mathsf{acc}(a,b)(\bullet))).$
\end{lem}
Informally, this lemma says that if every $b$-predecessor of $x$ is contained in the $b$-accessible part of $a$, then $x$ is also contained in the $b$-accessible part of $a$.

\begin{pf}
Assume $\mathrm{C}(a)$ and $\mathrm{C}(b)$. Taking any $c$, we further assume that $\mathrm{T}(ac)$ and $\forall y. \ \mathrm{T}(ay) \to \mathrm{T}(b(y,c)) \to \mathrm{T} (\mathsf{acc}(a,b)y)$,
then we have to derive $\mathrm{T}(\mathsf{acc}(a,b)c).$ The first conjunct $\mathrm{T}(\oplus(a, b) (0,c))$ of $\mathrm{T}(\mathsf{acc}(a,b)c)$ is equivalent to $\mathrm{T}(ac)$  in $\mathsf{KF}$ and is nothing but one of the assumptions.
To derive the second conjunct of $\mathrm{T}(\mathsf{acc}(a,b)c)$, we take any $d$, and we will derive: 
\begin{align}
\mathrm{T} ((\oplus(a, b) (0,d)) \dot{\to} (\oplus(a, b) (1, (d,c)) \dot{\to} (l(\oplus(a, b)) t(a,b,d)))). \label{pf_1:lem:IG1_in_KFUPI}
\end{align}
Then, the axioms of $\mathsf{KF}$ yield the desired conclusion. 

Since the term $\oplus(a, b)$ is a class by Lemma~\ref{lem:acc_is_class} and is reflected by $l(\oplus(a,b))$, the axioms of $\mathsf{KF}$ imply that the formula (\ref{pf_1:lem:IG1_in_KFUPI}) is equivalent to:
\begin{center}
$\mathrm{T}(\oplus(a, b) (0,d)) \to \mathrm{T}(\oplus(a, b) (1, (d,c))) \to \mathrm{T}(l(\oplus(a, b)) t(a,b,d)).$
\end{center}
Hence, we suppose that $\mathrm{T}(\oplus(a, b) (0,d))$ and $\mathrm{T}(\oplus(a, b) (1, (d,c)))$ in order to derive $\mathrm{T}(l(\oplus(a, b)) t(a,b,d))$. 
%$\mathrm{T}(l(\oplus(a, b)) t(a,b,d))$ is implied by suppositions
%$\mathrm{T}(\oplus(a, b) (0,d))$ and $\mathrm{T}(\oplus(a, b) (1, (d,c))).$ 
These suppositions are equivalent to $\mathrm{T}(ad)$ and $\mathrm{T}(b(d,c))$, respectively;
thus, letting $y:=d$, the initial assumption $\forall y. \ \mathrm{T}(ay) \to \mathrm{T}(b(y,c)) \to \mathrm{T} (\mathsf{acc}(a,b)y)$ implies $\mathrm{T}(\mathsf{acc}(a,b)d),$ and hence
\begin{center}
$\mathrm{T}(\dot{\forall} y (\oplus(a, b) (0,y) \dot{\to} (\oplus(a, b) (1, (y,d)) \dot{\to} (l(\oplus(a, b)) t(a,b,y)))))$
\end{center}
by the definition of $\mathsf{acc}.$
By the axioms $(\mathsf{K}_{\forall})$ and $(\mathsf{K}_{\to})$, we have 
\[
\forall y. \ \mathrm{T}  (\oplus(a, b) (0,y)) \to \mathrm{T} (\oplus(a, b) (1, (y,d))) \to \mathrm{T} (l(\oplus(a, b)) t(a,b,y)).
\]
As $\oplus(a, b)$ is a class, the axioms $(\mathsf{U} \mhyph \mathsf{True})$ and $(\mathsf{Lim})$ imply that
\[
\forall y. \ \mathrm{T}  (l(\oplus(a,b))(\dot{\neg}(\oplus(a, b) (0,y)))) \lor \mathrm{T} (l(\oplus(a,b))(\dot{\neg}(\oplus(a, b) (1, (y,d)))))
\]
\[
\lor \mathrm{T} (l(\oplus(a, b)) t(a,b,y)).
\]
Therefore, by using the axioms $(\mathsf{U}_{\to})$ and $(\mathsf{U}_{\forall}),$
we obtain 
\begin{center}
$\mathrm{T}(l (\oplus(a, b)) ( \dot{\forall} y (\oplus(a, b) (0,y) \dot{\to} (\oplus(a, b) (1, (y,d)) \dot{\to} t(a,b,y))))),$
\end{center}
which is nothing but the desired formula $\mathrm{T}(l(\oplus(a, b)) t(a,b,d))$.
Therefore, the formula~(\ref{pf_1:lem:IG1_in_KFUPI}) is derived. \qed
\end{pf}

Finally, we want to derive $(\mathsf{IG}.2)':$
\[
\mathrm{C}(a) \land \mathrm{C}(b) \land \mathsf{Closed}'(a,b,A(\bullet)) \to \forall x. \ \mathrm{T}(\mathsf{acc}(a,b)x) \to A(x).
\]
The author shall roughly explain the basic idea of the proof, which is essentially based on a standard technique found in, e.g., \cite{burgess2009friedman, cantini2016truth, jager2001universes, kahle2000theory}. Taking any $x$, we want to deduce $A(x)$ from the assumptions $\mathrm{C}(a), \mathrm{C}(b), \mathsf{Closed}'(a,b,A(\bullet))$ and $\mathrm{T}(\mathsf{acc}(a,b)x)$.
By the definition of $t$ and $\mathsf{acc},$ the last assumption implies that $t(a,b,x)$ is true, and hence is a proposition. 
In addition, by Lemma~\ref{lem:BA_is_A-closure} below, whether $A(x)$ holds or not can be reduced to whether $t(a,b,x)$ is a proposition. Therefore, $A(x)$ follows, as required.

We require the following lemma for the proof of Lemma~\ref{lem:BA_is_A-closure}.

\begin{lem}\label{lem:acc_is_closed_under_pred}
$\mathsf{KFU}$ derives the following:
\[
\mathrm{C}(u) \land \mathrm{C}(v) \land w \dot{\in} \mathsf{acc}(u,v) \to \forall x \dot{\in} u. \ x \dot{<}_{v} w \to x \dot{\in} \mathsf{acc}(u,v).
\]
\end{lem}

Informally, it says that if $w$ is contained in the $v$-accessible part of $u$,
then its every $v$-predecessor $x$ in $u$ also belongs in the $v$-accessible part of $u$.

\begin{pf}
Assume $\mathrm{C}(u), \mathrm{C}(v),$ and $w \dot{\in} \mathsf{acc}(u,v).$
In addition, we take any $x$ such that $x \dot{\in} u$ and $x \dot{<}_{v} w$.
Then, our purpose is to derive $x \dot{\in} \mathsf{acc}(u,v)$.
Recall that $w \dot{\in} \mathsf{acc}(u,v)$ is the following formula:
\[
\mathrm{T}(\oplus(u, v) (0,w)) \dot{\land} \dot{\forall} y (\oplus(u, v) (0,y) \dot{\to} [\oplus(u, v) (1,(y,w)) \dot{\to} (l(\oplus(u, v)) (t(u,v,y)))]).
\]
Thus, in a similar way as the proof of Lemma~\ref{lem:IG1_in_KFU}, we obtain:
\[
\forall y. \ \mathrm{T}(\oplus(u, v) (0,y)) \to \mathrm{T}(\oplus(u, v) (1,(y,w))) \to
\mathrm{T}(l(\oplus(u, v)) (t(u,v,y)))).
\]
For $y:=x,$ it follows from the assumptions that $\mathrm{T}(l(\oplus(u, v)) (t(u,v,x)))).$
Therefore, from the definition of $t(u,v,x)$, we similarly have: 
\[
\forall y. \ \mathrm{T}(\oplus(u, v) (0,y)) \to \mathrm{T}(\oplus(u, v) (1,(y,x))) \to
\mathrm{T}(l(\oplus(u, v)) (t(u,v,y)))).
\]
Since $\oplus(u, v)$ is a class, we obtain:
\[
\mathrm{T}(\dot{\forall} y (\oplus(u, v) (0,y) \dot{\to} [\oplus(u, v) (1,(y,x)) \dot{\to} (l(\oplus(u, v)) (t(u,v,y)))]).
\]
Moreover, by the assumption, we also have $\mathrm{T}(\oplus(u, v) (0,x))$; thus $(\mathsf{K}_{\land})$ yields $x \dot{\in} \mathsf{acc}(u,v),$ as required. \qed
\end{pf}

\begin{lem}\label{lem:BA_is_A-closure}
For an $\mathcal{L}_{FS}$-formula $A(x)$, define $B_{A}(u,v,w)$ to be the formula $\forall x. \
\mathrm{T}(\mathsf{acc}(u,v)x) \to F_A(u,v,w,x)$, 
where $F_A(u,v,w,x)$ is the conjunction of the following formulas:
\begin{enumerate}
\item $w = t(u,v,x) \ \to \ A(x),$

\item $\forall z. \ w = \{(\oplus(u, v) (1, (z,x))) \dot{\to} (t(u,v,z))\} \to \mathrm{T}(\oplus(u, v) (1, (z,x))) \to A(z),$

\item $\forall z. \ w = \{(\oplus(u, v) (0,z)) \dot{\to} (\oplus(u, v) (1, (z,x)) \dot{\to} t(u,v,z))\} \to
%$ \ \ \ \to 
\mathrm{T}(\oplus(u, v) (0,z)) \to \mathrm{T}(\oplus(u, v) (1, (z,x))) \to A(z).$

%\item $w = \dot{\forall} z (\oplus(u, v) (0,z) \dot{\to} (\oplus(u, v) (1, (z,x)) \dot{\to} t(u,v,z)))$

%$ \ \ \ \to \forall z. \mathrm{T}(\oplus(u, v) (0,z)) \to \mathrm{T}(\oplus(u, v) (1, (z,x))) \to A(z)$
\end{enumerate}

Then, $\mathsf{KFUPI}$ derives the following:
\[
\mathrm{C}(u) \land \mathrm{C}(v) \land \mathsf{Closed}' (u,v,A(\bullet)) \to \forall w. \ \mathscr{A}(B_{A}(u,v, \bullet), w) \to B_{A}(u,v,w). 
\]

\end{lem}

\begin{pf}
Assume $\mathrm{C}(u), \mathrm{C}(v), \mathsf{Closed}' (u,v,A(\bullet))$ and $\mathscr{A}(B_{A}(u,v, \bullet), w).$ Then, we need to prove $B_{A}(u,v,w)$. So, taking any $x$ such that $\mathrm{T}(\mathsf{acc}(u,v)x)$, we show $F_A(u,v,w,x).$ The proof is mainly divided into three cases according to the form of $w$.
\begin{enumerate}
\item Assume $w = t(u,v,x):=
\dot{\forall} z (\oplus(u, v) (0,z) \dot{\to} (\oplus(u, v) (1, (z,x)) \dot{\to} t(u,v,z))). $
Then, we have to show $A(x).$
Now, $w$ is of the universal form $\dot{\forall}f$, thus by the definition of $\mathscr{A}(B_{A}(u,v, \bullet), w)$ and $(\mathsf{UG})$, we have:
\[
[\forall z. \ B_A(u,v, fz)] \lor \exists z. \ B_A(u,v, fz) \land \mathrm{T}(\dot{\neg} (fz)).
\]
However, in a similar manner to the proof of Lemma~\ref{lem:IG1_in_KFU}, the assumption $\mathrm{T}(\mathsf{acc}(u,v)x)$ implies, in $\mathsf{KFU},$ $\forall z. \ \mathrm{T}(fz)$; that is,
\[
\forall z. \ \mathrm{T}(\oplus(u, v) (0,z) \dot{\to} (\oplus(u, v) (1, (z,x)) \dot{\to} t(u,v,z))).
\] 
Thus, $\exists z. \ \mathrm{T}(\dot{\neg} (fz))$ is false by $(\mathsf{T} \mhyph \mathsf{Cons})$. Therefore, it follows that $\forall z. \ B_A(u,v, fz).$ %Now we take any $z$ such that $z \dot{\in} u$ and $z \dot{<}_v x,$ then 
%by Lemma~\ref{lem:acc_is_closed_under_pred}, 
%we have $\mathsf{acc}(u,v)z.$ 
%$B_A(u,v, fz)$ and
Then, the assumption $\mathsf{acc}(u,v)x$ yields $\forall z. \ F_A(u,v,fz,x).$
Since $fz$ is of the conditional form, the third conjunct of 
$F_A(u,v, fz,x)$ is applied, and thus, for all $z$, we have:
\[
\mathrm{T}(\oplus(u, v) (0,z)) \to \mathrm{T}(\oplus(u, v) (1, (z,x))) \to A(z),
\]
which is equivalent to $\forall z \dot{\in} u. \ z \dot{<}_v x \to A(z)$.
Because $\mathrm{T}(\mathsf{acc}(u,v)x)$ implies $x \dot{\in}u,$
the assumption $\mathsf{Closed}' (u,v,A(\bullet))$ implies the formula $A(x),$ as required.

\item Taking any $z$, we assume $w = (\oplus(u, v) (1, (z,x))) \dot{\to} t(u,v,z)$.
Then, we will deduce $\mathrm{T}(\oplus(u, v) (1, (z,x))) \to A(z)$.
Thus, similar to the case of $w = t(u,v,x)$, the assumption $\mathscr{A}(B_{A}(u,v, \bullet), w)$ implies the following:
\begin{align}
\mathrm{T}(\oplus(u, v) (1, (z,x))) \to B_A(u,v,t(u,v,z)). \label{pf_1:lem:BA_is_A-closure}
\end{align}
By the assumption $\mathrm{T}(\mathsf{acc}(u,v)x)$ and Lemma~\ref{lem:acc_is_closed_under_pred}, we have $\mathrm{T}(\mathsf{acc}(u,v)z)$. Therefore, (\ref{pf_1:lem:BA_is_A-closure}) implies:
\[
\mathrm{T}(\oplus(u, v) (1, (z,x))) \to F_A(u,v,t(u,v,z),z).
\] 
Thus, the first conjunct of $F_A(u,v,t(u,v,z),z)$ is applied, and we get
\[
\mathrm{T}(\oplus(u, v) (1, (z,x))) \to A(z).
\] 

\item The case $w = \oplus(u, v) (0,z) \dot{\to} (\oplus(u, v) (1, (z,x)) \dot{\to} t(u,v,z))$ is similarly proved as the second case.

\item If $w$ is of the other form, we trivially have $B_{A}(u,v,w).$ \qed
\end{enumerate}
\end{pf}

Using the above lemma, we can derive $(\mathsf{IG}. 2)'$.
\begin{lem}\label{lem:IG2_in_KFUPI}
$\mathsf{KFUPI}$ derives the following:
\[
\mathrm{C}(a) \land \mathrm{C}(b) \land Closed(a,b,A) \to \forall x. \ \mathrm{T}(\mathsf{acc}(a,b)x) \to A(x).
\]
\end{lem}

\begin{pf}
We assume $\mathrm{C}(a), \mathrm{C}(b), Closed(a,b,A)$ and $\mathrm{T}(\mathsf{acc}(a,b)c)$ for an arbitrary $c$. Then, we have to derive $A(c)$.
From the assumptions and Lemma~\ref{lem:BA_is_A-closure}, we obtain $\forall w. \ \mathscr{A}(B_{A}(u,v, \bullet), w) \to B_{A}(u,v,w).$
Thus,  $(\mathsf{PI})$ yields:
\[
\forall x. \ \mathrm{P}(x) \to B_{A}(a,b,x).
\]
Letting $x := t(a,b,c),$ we want to show $\mathrm{P}(t(a,b,c)),$ %that is,
%\begin{center}
%$\mathrm{T}(l (\oplus(a, b)) (\dot{\forall} y (\oplus(a, b) (0,y) \dot{\to} [\oplus(a, b) (1, (y,c)) \dot{\to} t(a,b,y)])) ),$
%\end{center}
from which $B_{A}(a,b,t(a,b,c))$ follows. Then, we can immediately obtain the conclusion $A(c)$ by the assumption $\mathrm{T}(\mathsf{acc}(a,b)c)$ and the definition of $B_{A}$. 

In order to prove $\mathrm{P}(t(a,b,c)),$ it suffices to derive the following:
\begin{align}
%\item $\mathrm{T}(l (\oplus(a, b)) (\oplus(a, b) (0,c)));$
\forall y. \ \mathrm{T}(\oplus(a, b) (0,y)) \to \mathrm{T}(\oplus(a, b) (1, (y,c))) \to \mathrm{T}(t(a,b,y)). \label{pf_1:lem:IG2_in_KFUPI}
\end{align}
In fact, this formula implies, in $\mathsf{KFU}$, $\mathrm{T}(t(a,b,c) )$ and thus  $\mathrm{P}(t(a,b,c) )$ because $\oplus(a, b)$ is a class.
%\begin{center}
%$\mathrm{T}((l (\oplus (a,b))) (\dot{\forall} y ((\oplus(a, b) (0,y)) \dot{\to} [(\oplus(a, b) (1, (y,c))) \dot{\to} t(a,b,y)]))),$
%\end{center}
%which is the same as the desired formula $\mathrm{T}(l (\oplus(a, b)) t(a,b,c) ).$% that is,
%\begin{center}
%$\mathrm{T}((l (\oplus (a,b))) ((\oplus(a, b) (0,w)) \dot{\land} (\dot{\forall} y ((\oplus(a, b) (0,y)) \dot{\to} [(\oplus(a, b) (1, (y,w))) \dot{\to} t(a,b,y)])))).$
%\end{center}
\  
However,  the formula (\ref{pf_1:lem:IG2_in_KFUPI}) follows from the assumption $\mathrm{T}(\mathsf{acc}(a,b)c)$. \qed
\end{pf}

To summarise, we obtain the interpretation of $\mathsf{T}_0$ into $\mathsf{KFUPI},$ in which $\mathcal{L}$-vocabularies are preserved.

\begin{thm}\label{thm:lower-bound_of_KFUPI}
For each $\mathcal{L}_{EM}$-sentence $A$,
if $\mathsf{T}_{0} \vdash A$, then $\mathsf{KFUPI} \vdash A'.$

\end{thm}

\subsection{Kahle's model for $\mathsf{KFU}$}\label{subsec:Kahle's_model_for_KFU}
In this subsection, we introduce Kahle's model construction for $\mathsf{KFU}$
\cite[pp.~214-215]{kahle2003universes} and observe that this model also satisfies $\mathsf{KFUPI}$. The upper-bound proof of $\mathsf{KFUPI}$ in Section~\ref{subsec:upper-bound_of_KFUPI} is essentially based on this model.

First, %fixing some G\"{o}del-numbering, we identify every closed term of $\mathcal{L}_{FS}$ with its code. Then, 
the base theory $\mathsf{TON}$ is interpreted by the closed total term model $\mathcal{CTT}$ (e.g. \cite[p.~26]{cantini1996logical});
the domain of $\mathcal{CTT}$ is the set of all closed $\mathcal{L}$-terms;
the constant symbols are interpreted as themselves;
the application function $\mathsf{App}(x,y)$ is defined to be the juxtaposition of $x$ and $y$;
and it is defined that the equation $x=y$ holds when $x$ and $y$ are $\beta$-equivalent in the standard sense. 
Then, $\mathrm{N}(x)$ holds when $x$ reduces to a numeral.
%As for the vocabularies of $\mathcal{L}_{FS}$, the additional constants $\dot{=}, \dot{\neg}, \dot{\land}, \dot{\to}, \dot{\forall}, \dot{\mathsf{N}},  \mathsf{l}$ of $\mathcal{L}_{FS}$ can be interpreted as closed terms so that the schema $(\mathsf{UG})$ \ref{def:UG} is satisfied. 
To interpret the predicates $\mathrm{T}$ and $\mathrm{U},$ we define an operator $\Phi(X,Y,a,\alpha)$ for an ordinal number $\alpha:$
\begin{enumerate}
\item $a \in Y,$
\item $\exists b,c. \ \alpha = 0 \land a = (b \dot{=} c) \land b=c,$
\item $\exists b,c. \ \alpha = 0 \land a = \dot{\neg}(b \dot{=} c) \land b\neq c,$
\item $\exists b. \ \alpha = 0 \land a = \dot{\mathsf{N}}b \land \mathrm{N}(b), $
\item $\exists b. \ \alpha = 0 \land a = \dot{\neg}(\dot{\mathsf{N}}b) \land \neg\mathrm{N}(b), $
\item $\exists b. \ a = \dot{\neg}(\dot{\neg}b) \land a \notin Y \land b \in X,$
\item $\exists b,c. \ a = (b \dot{\land} c) \land a \notin Y \land b \in X \land c \in X,$
\item $\exists b,c. \ a = \dot{\neg}(b \dot{\land} c) \land a \notin Y \land [\dot{\neg}b \in X \lor \dot{\neg}c \in X],$
\item $\exists b,c. \ a = (b \dot{\to} c) \land a \notin Y \land [\dot{\neg}b \in X \lor c \in X],$
\item $\exists b,c. \ a = \dot{\neg}(b \dot{\to} c) \land a \notin Y \land b \in X \land \dot{\neg}c \in X,$
\item $\exists f. \ a = (\dot{\forall}f) \land a \notin Y \land \forall x. \ fx \in X,$
\item $\exists f. \ a = \dot{\neg}(\dot{\forall}f) \land a \notin Y \land \exists x. \ \dot{\neg}(fx) \in X,$
\item $\exists b,c. \ \alpha \in \mathrm{Suc} \land a=\mathsf{l}bc \land a \notin Y \land \dot{\neg}a \notin Y \land [\forall x. \ bx \in Y \lor \dot{\neg}(bx) \in Y] \land c \in Y,$
\item $\exists b,c. \ \alpha \in \mathrm{Suc} \land a=\dot{\neg}(\mathsf{l}bc) \land a \notin Y \land \dot{\neg}a \notin Y \land [\forall x. \ bx \in Y \lor \dot{\neg}(bx) \in Y] \land c \notin Y,$
\end{enumerate}
where $\alpha \in \mathrm{Suc}$ means that $\alpha$ is a successor ordinal.

A $\Phi$-sequence $(Z_{\alpha})$ of least fixed points $Z_{\alpha}$ is defined as follows:
\begin{enumerate}
\item $Z_0 := \emptyset.$
\item $Z_{\alpha+1}:=$ the least fixed point of $\Phi(X,Z_{\alpha},x,\alpha+1).$
\item $Z_{\alpha}:=$ the least fixed point of $\Phi(X, \bigcup_{\beta < \alpha}Z_{\beta}, x, \alpha)$, for a limit ordinal $\alpha$.
\end{enumerate} 
Since this sequence $(Z_{\alpha})$ is monotonically increasing, we eventually reach the least ordinal $\iota$ such that $Z_{\iota}= Z_{\iota +1}.$\footnote{In fact, this ordinal $\iota$ is identical to the first recursively inaccessible ordinal, for the definition of which, see, e.g. \cite{barwise1975admissible}. } 

Then, we interpret $\mathrm{T}(t)$ as $t \in Z_{\iota},$ and $\mathrm{U}(u)$ as $\exists f. \ u=\mathsf{l}f \land \forall x. \ fx \in Z_{\iota} \lor \dot{\neg}(fx) \in Z_{\iota}.$
Combining this with the above interpretation, we obtain an $\mathcal{L}_{FS}$-model $\mathcal{M}_{\mathsf{KF}}.$

\begin{prop}\label{prop:model_of_KFUPI}
$\mathcal{M}_{\mathsf{KF}} \models \mathsf{KFUPI}$.
\end{prop}

\begin{pf}\label{pf:prop:model_of_KFUPI}
Kahle \cite[p.~215]{kahle2003universes} already showed $\mathcal{M}_{\mathsf{KF}} \models \mathsf{KFU} + (\mathsf{UG})$; thus, we concentrate on the schema (\ref{PI}):
\[
[\forall x. \ \mathscr{A}(B(\bullet),x) \to B(x)] \to \forall x. \ \mathrm{P}(x) \to B(x).
\]
We take any $\mathcal{L}_{FS}$-formula $B$ such that $\mathcal{M}_{\mathsf{KF}} \models \forall x. \ \mathscr{A}(B(\bullet),x) \to B(x)$. 
%Since this $B(x)$ is closed under $\dot{\neg},$
%it suffices to show that $\mathcal{M}_{\mathsf{KF}} \models \forall x. \mathrm{T}(\dot{\neg}x) \to B(x).$ 
Then, by induction on $\alpha \leq \iota$, we prove that
if $x \in Z_{\alpha}$ or $\dot{\neg}x \in Z_{\alpha}$, then $\mathcal{M}_{\mathsf{KF}} \models B(x).$ 

For example, we consider the case $x = \mathsf{l}fy$ with $\dot{\neg}(\mathsf{l}fy) \in Z_{\alpha}$.
Then, we must show $\mathcal{M}_{\mathsf{KF}} \models B(\mathsf{l}fy).$
By the definition of the operator $\Phi,$ the crucial case is when $\alpha$ is the least successor ordinal such that $\forall z. \ fz \in \bigcup_{\beta< \alpha}Z_{\beta} \lor \dot{\neg}(fz) \in \bigcup_{\beta< \alpha}Z_{\beta}.$ Thus, by the induction hypothesis, we obtain
$\mathcal{M}_{\mathsf{KF}} \models \forall x. \ B(fz)$. 
Since $B$ is $\mathscr{A}$-closed, it follows that $\mathcal{M}_{\mathsf{KF}} \models B(\mathsf{l}fy).$
The other cases are similarly proved by using subinduction on the construction of $Z_{\alpha}.$ \qed
\end{pf}

\begin{rmk}\label{rmk:structure_of_M_KF}
As Kahle remarks \cite[p.~215]{kahle2003universes}, the model $\mathcal{M}_{\mathsf{KF}}$ also satisfies the following principles:
\begin{description}
\item[$(\mathsf{U} \mhyph \mathsf{Tran})$] $\mathrm{U}(u) \land \mathrm{U}(v) \to [\mathrm{T}(v(ux)) \to \mathrm{T}(vx)]$.
\item[$(\mathsf{U} \mhyph \mathsf{Dir})$] $\mathrm{U}(u) \land \mathrm{U}(v) \to \exists w. \ \mathrm{U}(w) \land u \sqsubset w \land v \sqsubset w.$
\item[$(\mathsf{U} \mhyph \mathsf{Nor})$] $\mathrm{U}(u) \to \exists f. \ \mathrm{C}(f) \land u = \mathsf{l}f.$
\item[$(\mathsf{U} \mhyph \mathsf{Lin})$] $\mathrm{U}(u) \land \mathrm{U}(v) \to u \sqsubset v \lor v \sqsubset u \lor \forall x. x \dot{\in} u \leftrightarrow x \dot{\in} v.$

\end{description}
Therefore, Proposition~\ref{prop:model_of_KFUPI} also verifies the consistency of $\mathsf{KFUPI} + (\mathsf{U} \mhyph \mathsf{Tran}) + (\mathsf{U} \mhyph \mathsf{Dir}) + (\mathsf{U} \mhyph \mathsf{Nor}) + (\mathsf{U} \mhyph \mathsf{Lin}).$
\end{rmk}

\subsection{Upper bound of $\mathsf{KFUPI}$}
\label{subsec:upper-bound_of_KFUPI}

In order to obtain the upper bound of $\mathsf{KFUPI},$ we interpret the truth predicate of $\mathsf{KFUPI}$ as the least fixed point of a non-monotone operator which is almost the same as $\Phi$ in Section~\ref{subsec:Kahle's_model_for_KFU}. In the operator $\Phi$, the truth condition of $\dot{\neg}(\mathsf{l}bc)$ for a class $b$ was characterised as the \emph{non}-truth of $c$, where non-monotonicity concerns. That is why monotone inductive operators are not enough to formalise $\mathcal{M}_{\mathsf{KF}}$.
Thus, we use the theory $\mathsf{FID([POS,QF])}$ of non-monotone inductive definition in \cite{jager2001first}.

Let $\mathcal{L}'$ be an language of Peano arithmetic $(\mathsf{PA})$ that contains symbols for all primitive recursive functions. Let $\mathcal{L}'(X)$ be the extension of $\mathcal{L}'$ with a new unary predicate symbol $x \in X$. An \emph{operator form} $\mathfrak{B}(X,u)$ is an $\mathcal{L}'(X)$-formula in which at most one variable $u$ occurs freely. We write $\mathfrak{B}(A,u)$ as the result of replacing each occurrence of $t \in X$ by a formula $A(t).$
An operator form $\mathfrak{B}(X, u)$ is in $\mathsf{POS}$ if $X$ occurs only positively in $\mathfrak{B};$
an operator form $\mathfrak{B}(X, u)$ is in $\mathsf{QF}$ if it contains no quantifier.
%For each pair of operator forms 
We consider the particular operators $\mathfrak{A}_{0} \in \mathsf{POS}$ and $\mathfrak{A}_{1} \in \mathsf{QF}$, which are defined below.
Then, the operator $\mathfrak{A}$ is defined to be the following formula:
%The set $[\mathsf{POS,QF}]$ is the collection of all operator forms of the following form:
\begin{center}
%$\mathfrak{A}(\mathrm{P},\vec{u}) : = 
$\mathfrak{A}(X,u):= \mathfrak{A}_{0}(X, u) \lor ([\forall x. \ \mathfrak{A}_{0}(X, x) \to x \in X] \land \mathfrak{A}_{1}(X,u)).$
\end{center}
%where $\mathfrak{A}_{0}(\mathrm{P}, \vec{u}) \in \mathsf{POS}$ and $\mathfrak{A}_{1}(\mathrm{P}, \vec{u}) \in \mathsf{QF}.$

Let the two-sorted language $\mathcal{L}_{\mathcal{K}}$ be the extension of $\mathcal{L}'$ with ordinal variables $\alpha, \beta, \dots,$ a binary relation symbol $<$ on the ordinals, and the binary relation symbol $\mathrm{P}_{\mathfrak{A}}$ for the above operator $\mathfrak{A}$.

We shall use the following notations:
\begin{itemize}
\item $\mathrm{P}^{\alpha}_{\mathfrak{A}}(s) := \mathrm{P}_{\mathfrak{A}}(\alpha, s),$
\item $\mathrm{P}^{< \alpha}_{\mathfrak{A}}(s) := \exists \xi < \alpha. \ \mathrm{P}^{\xi}_{\mathfrak{A}}(s),$
\item $\mathrm{P}_{\mathfrak{A}}(s) := \exists \alpha. \ \mathrm{P}^{\alpha}_{\mathfrak{A}}(s).$
\end{itemize}

\begin{defn}\label{def:FID}
The $\mathcal{L}_{\mathcal{K}}$-theory $\mathsf{FID}([\mathsf{POS,QF}])$ consists of $\mathsf{PA}$ with the full induction schema for $\mathcal{L}_{\mathcal{K}}$ and the following axioms:\footnote{In contrast to J\"{a}ger and Studer's formulation of $\mathsf{FID}([\mathsf{POS,QF}])$ in \cite{jager2001first}, we now consider only one operator form $\mathfrak{A}$ for simplicity.}
\begin{enumerate}
\item $\alpha \nless \alpha \land (\alpha < \beta \land \beta < \gamma \to \alpha < \gamma) \land (\alpha < \beta \lor \alpha = \beta \lor \alpha > \beta).$
\item %For all operator forms $\mathfrak{A}(\mathrm{P},\vec{u}) \in [\mathsf{POS,QF}],$
\begin{description}
\item[$(\mathsf{OP}.1)$] $\mathrm{P}^{\alpha}_{\mathfrak{A}}(s) \leftrightarrow \mathrm{P}^{< \alpha}_{\mathfrak{A}}(s) \lor \mathfrak{A}(\mathrm{P}^{< \alpha}_{\mathfrak{A}}, s),$
\item[$(\mathsf{OP}.2)$] $\mathfrak{A}(\mathrm{P}_{\mathfrak{A}}, s) \to \mathrm{P}_{\mathfrak{A}}(s).$
\end{description}
\item $[\forall \xi. \ \{ \forall \eta < \xi. \ A(\eta)\} \to A(\xi)] \to \forall \xi. \ A(\xi),$ for all $\mathcal{L}_{\mathcal{K}}$-formulae $A(\alpha)$.

\end{enumerate}
\end{defn}

Based on \cite[Section~5.2.2]{kahle2003universes}, we now formalise the closed total term model of $\mathsf{KFUPI}$ within $\mathsf{FID}([\mathsf{POS,QF}])$.
First, we define a translation $^{\star}: \mathcal{L} \to \mathcal{L}'$.
Fixing some G\"{o}del-numbering, each closed term $t$ is assigned the corresponding natural number $\ulcorner t \urcorner$; thus, the domain of $\mathcal{CTT}$ can be seen as a subset of $\mathbb{N}$, which is expressed by some arithmetical formula. Thus, the quantifier symbols are relativised to this formula.
The translation $\mathsf{App}^{\star}$ of the application function is defined to be the standard Kleene bracket $\{ x \} (y)$, which returns the result of the partial recursive function $\{ x \}$ at $y$ if it has a value.
Each constant symbol $\mathsf{c}$ is interpreted as a code $e$ such that $\{ e \} (\ulcorner t \urcorner) \simeq \ulcorner \mathsf{c} t \urcorner$ for each closed term $t$.
Similarly, $\mathrm{N}^{\star}(x)$ is an arithmetical formula  expressing that $x$ is the code of a numeral, and an arithmetical formula $x =^{\star} y$ is true when $x$ and $y$  have the same reduct. 
Under this interpretation, every theorem of $\mathsf{TON}$ is derivable in $\mathsf{PA}$ (cf. \cite[Theorem~4.13]{cantini1996logical}).
Next, to expand $^{\star}$ to the language $\mathcal{L}_{FS}$, we need to define the extension of the truth predicate $\mathrm{T}$.

We define the operator form $\mathfrak{A}_0(X,a) \in \mathsf{POS}$ to be the disjunction of the following:
\begin{enumerate}
\item $\exists b,c. \ a = (b \dot{=} c) \land b = c$
\item $\exists b,c. \ a = \dot{\neg} (b \dot{=} c) \land b \neq c$
\item $\exists b. \ a = (\dot{\mathsf{N}} b) \land \mathrm{N}(b)$
\item $\exists b. \ a = \dot{\neg} (\dot{\mathsf{N}} b) \land \neg \mathrm{N}(b)$
\item $\exists b. \ a = \dot{\neg}(\dot{\neg} b) \land b \in X$
\item $\exists b,c. \ a = (b \dot{\land} c) \land b \in X \land c \in X$
\item $\exists b,c. \ a = \neg (b \dot{\land} c) \land [\dot{\neg}b \in X \lor \dot{\neg}c \in X]$
\item $\exists b,c. \ a = (b \dot{\to} c) \land [\dot{\neg}b \in X \lor c \in X]$
\item $\exists b,c. \ a = \dot{\neg} (b \dot{\to} c) \land b \in X \land \dot{\neg} c \in X$
\item $\exists f. \ a = \dot{\forall} f \land \forall b. \ fb \in X$
\item $\exists f. \ a = \dot{\neg}(\dot{\forall} f) \land \exists b. \ \dot{\neg}(fb) \in X$
\item $\exists f. \ a = \mathsf{l}f(\mathsf{0} \dot{=} \mathsf{0}) \land \forall b. \ fb \in X \lor \dot{\neg}fb \in X $
\end{enumerate}
Note that the last clause is given to record that $f$ is a class at the stage.
Next, let an operator form $\mathfrak{A}_{1}(X, a) \in \mathsf{QF}$ be the disjunction of the following:
\begin{enumerate}
\item $\exists f,b. \ a = \mathsf{l}fb \land \mathsf{l}f(\mathsf{0} \dot{=}\mathsf{0}) \in X \land \mathsf{l}f(\dot{\mathsf{N}}\mathsf{0}) \notin X \land b \in X.$
\item $\exists f,b. \ a = \dot{\neg} (\mathsf{l}fb) \land \mathsf{l}f(\mathsf{0} \dot{=}\mathsf{0}) \in X \land \mathsf{l}f(\dot{\mathsf{N}}\mathsf{0}) \notin X \land b \notin X.$
\end{enumerate}
It should be clear that $\mathfrak{A}_1$ can be given as a quantifier-free formula.
The condition that $f$ is a class is now given by the quantifier-free formula $\mathsf{l}f(0 \dot{=}0) \in X,$ so $\mathfrak{A}_{1}(X, a)$ can dispense with quantifier symbols.
The condition $\mathsf{l}f(\dot{\mathsf{N}}\mathsf{0}) \notin X$ is used to ensure that 
every term of the form $\mathsf{l}fu$ or $\dot{\neg}(\mathsf{l}fu)$ is simultaneously determined to be true or false at some stage (see Lemma~\ref{lem:simultaneousness_in_U}).
Finally, the operator form $\mathfrak{A}$ is defined, as explained above:
\[
\mathfrak{A}(X,a) : = \mathfrak{A}_{0}(X, a) \lor ([\forall x. \ \mathfrak{A}_{0}(X, x) \to x \in X] \land \mathfrak{A}_{1}(X,a)).
\]
To complete the definition of the translation $^{\star}$, let $\mathrm{T}^{\star}(t)$ be $\mathrm{P}_{\mathfrak{A}}(t),$
and $\mathrm{U}^{\star}(u)$ be $\exists f. \ u = \mathsf{l}f \land \forall x. \mathrm{P}_{\mathfrak{A}}(fx) \lor \mathrm{P}_{\mathfrak{A}}(\dot{\neg}(fx))$.
Under the interpretation $^\star$, we want to derive all the theorems of $\mathsf{KFUPI}$ in $\mathsf{FID}([\mathsf{POS,QF}]).$ For that purpose, we use the following lemma, which shows that when $f$ is a class,  $\mathsf{l}fb$ is determined to be true or false at some stage $\alpha$, and the truth value never changes at the later stages. 

\begin{lem}\label{lem:simultaneousness_in_U}
We work in $\mathsf{FID}([\mathsf{POS,QF}])$.

Assume $\forall x. \ \mathrm{P}_{\mathfrak{A}}(fx) \lor \mathrm{P}_{\mathfrak{A}}(\dot{\neg}(fx)).$ Then, there exists an ordinal $\alpha$ such that $\mathrm{P}^{\alpha}_{\mathfrak{A}}(\mathsf{l}f(\dot{\mathsf{N}\mathsf{0}})) \land \neg \mathrm{P}^{<\alpha}_{\mathfrak{A}}(\mathsf{l}f(\dot{\mathsf{N}\mathsf{0}})).$
Furthermore, for this $\alpha,$ we have the following:
\begin{enumerate}
\item $\forall b. \ \mathrm{P}_{\mathfrak{A}}(\mathsf{l}fb) \leftrightarrow \mathrm{P}^{\alpha}_{\mathfrak{A}}(\mathsf{l}fb),$
\item $\forall b. \ \mathrm{P}_{\mathfrak{A}}(\dot{\neg}(\mathsf{l}fb)) \leftrightarrow \mathrm{P}^{\alpha}_{\mathfrak{A}}(\dot{\neg}(\mathsf{l}fb)).$
\end{enumerate}
\end{lem}

\begin{pf}\label{pf:lem:simultaneousness_in_U}
We assume that $\forall x. \ \mathrm{P}_{\mathfrak{A}}(fx) \lor \mathrm{P}_{\mathfrak{A}}(\dot{\neg}(fx))$. Firstly, we prove $\mathrm{P}_{\mathfrak{A}}(\mathsf{l}f(\dot{\mathsf{N}}\mathsf{0})).$ For a contradiction, we suppose $\neg \mathrm{P}_{\mathfrak{A}}(\mathsf{l}f(\dot{\mathsf{N}}\mathsf{0})).$ Since $(\mathsf{OP}.2)$ implies that $\mathrm{P}_{\mathfrak{A}}$ is closed under $\mathfrak{A}_0$, we have $\mathrm{P}_{\mathfrak{A}}(\mathsf{l}f(\mathsf{0} \dot{=} \mathsf{0})).$ Again by $(\mathsf{OP}.2)$, it follows that $\forall a. \ \mathfrak{A}_1(\mathrm{P}_{\mathfrak{A}},a) \to \mathrm{P}_{\mathfrak{A}}(a).$ Therefore, by the supposition $\neg \mathrm{P}_{\mathfrak{A}}(\mathsf{l}f(\dot{\mathsf{N}}\mathsf{0})),$ we obtain $\mathrm{P}_{\mathfrak{A}}(\mathsf{l}f(\dot{\mathsf{N}}\mathsf{0})),$ a contradiction. Thus, $\mathrm{P}_{\mathfrak{A}}(\mathsf{l}f(\dot{\mathsf{N}}\mathsf{0}))$ is proved.
By the transfinite induction schema, we also have an ordinal $\alpha$ such that 
$\mathrm{P}^{\alpha}_{\mathfrak{A}}(\mathsf{l}f(\dot{\mathsf{N}}\mathsf{0})) \land \neg \mathrm{P}^{<\alpha}_{\mathfrak{A}}(\mathsf{l}f(\dot{\mathsf{N}}\mathsf{0})).$

Next, we show $\mathrm{P}_{\mathfrak{A}}(\mathsf{l}fb) \to \mathrm{P}^{\alpha}_{\mathfrak{A}}(\mathsf{l}fb).$ Note that the right-to-left direction is obvious.
Since $\mathrm{P}^{<\alpha}_{\mathfrak{A}}(\mathsf{0} \dot{=} \mathsf{0}),$
it suffices to consider the case $b \neq (\mathsf{0} \dot{=} \mathsf{0}).$
Now, we assume $\mathrm{P}_{\mathfrak{A}}(\mathsf{l}fb)$, and hence we can take a least $\beta$ such that $\mathrm{P}^{\beta}_{\mathfrak{A}}(\mathsf{l}fb)$.
Then, $\beta$ is clearly equal to $\alpha,$ and thus we have $\mathrm{P}^{\alpha}_{\mathfrak{A}}(\mathsf{l}fb),$ as required. We can similarly verify the second item $\mathrm{P}_{\mathfrak{A}}(\dot{\neg}(\mathsf{l}fb)) \leftrightarrow \mathrm{P}^{\alpha}_{\mathfrak{A}}(\dot{\neg}(\mathsf{l}fb)).$ \qed
\end{pf}

\begin{lem}\label{lem:interpretation_of_KFUPI}
For every $\mathcal{L}_{FS}$-formula $A,$ if $\mathsf{KFUPI} \vdash A$, then $\mathsf{FID}([\mathsf{POS,QF}]) \vdash A^{\star}$.
\end{lem}

\begin{pf}\label{pf:lem:interpretation_of_KFUPI}
If $A$ is an axiom of $\mathsf{TON}$ or $(\mathsf{UG})$, we already have $\mathsf{PA} \vdash A^{\star}$. 
Thus, it suffices to deal with the axioms displayed in Definition~\ref{def:KFU} and the schema (\ref{PI}).
\begin{description}
\item[Compositional axioms in $\mathrm{U}$.]
We consider the axiom $(\mathsf{U}_{\forall}):$
\[
\mathrm{U}^{\star}(u) \to \forall g. \ \mathrm{T}^{\star}(u(\dot{\forall} g)) \leftrightarrow \forall x. \ \mathrm{T}^{\star}(u(gx)).
\]

Thus, we assume $\mathrm{U}^{\star}(u)$, so we can take a closed term $f$ such that $u=\mathsf{l}f \land \forall x. \ \mathrm{P}_{\mathfrak{A}}(fx) \lor \mathrm{P}_{\mathfrak{A}}(\dot{\neg}(fx)).$ By Lemma~\ref{lem:simultaneousness_in_U}, there exists an ordinal $\alpha$ such that $\mathrm{P}^{\alpha}_{\mathfrak{A}}(\mathsf{l}f(\dot{\mathsf{N}}\mathsf{0})) \land \neg \mathrm{P}^{<\alpha}_{\mathfrak{A}}(\mathsf{l}f(\dot{\mathsf{N}}\mathsf{0})).$ Moreover, for any closed term $g,$ we have:
\begin{align}
\mathrm{P}_{\mathfrak{A}}(u(\dot{\forall} g)) &\leftrightarrow \mathrm{P}^{\alpha}_{\mathfrak{A}}(u(\dot{\forall} g)) \tag{$\because$ Lemma~\ref{lem:simultaneousness_in_U}} \\
&\leftrightarrow \mathrm{P}^{<\alpha}_{\mathfrak{A}}(\dot{\forall} g) \tag{$\because$ def. of $\mathfrak{A}$} \\
&\leftrightarrow \forall x. \ \mathrm{P}^{<\alpha}_{\mathfrak{A}}(gx) \tag{$\because$ $\mathrm{P}^{<\alpha}$ is $\mathfrak{A}_0$-closed} \\
&\leftrightarrow \forall x. \ \mathrm{P}^{\alpha}_{\mathfrak{A}}(u(gx)) \tag{$\because$ def. of $\mathfrak{A}$} \\
&\leftrightarrow \forall x. \ \mathrm{P}_{\mathfrak{A}}(u(gx)). \tag{$\because$ Lemma~\ref{lem:simultaneousness_in_U}}
\end{align}
Therefore, we obtain $\forall g. \ \mathrm{T}^{\star}(u(\dot{\forall} g)) \leftrightarrow \forall x. \ \mathrm{T}^{\star}(u(gx))$. The other compositional axioms are similarly treated.

\item[Consistency in $\mathrm{U}$.]
We consider the axiom $(\mathsf{U} \mhyph \mathsf{Cons})$:
\[
\mathrm{U}^{\star}(u) \to \forall x. \ \neg [\mathrm{T}^{\star}(ux) \land \mathrm{T}^{\star}(u(\dot{\neg}x))].
\]
As above, from the assumption $\mathrm{U}^{\star}(u)$, we take a closed term $f$ such that $u=\mathsf{l}f \land \forall x. \mathrm{P}_{\mathfrak{A}}(fx) \lor \mathrm{P}_{\mathfrak{A}}(\dot{\neg}(fx)).$
Also, we get an ordinal $\alpha$ such that $\mathrm{P}^{\alpha}_{\mathfrak{A}}(\mathsf{l}f(\dot{\mathsf{N}}\mathsf{0})) \land \neg \mathrm{P}^{<\alpha}_{\mathfrak{A}}(\mathsf{l}f(\dot{\mathsf{N}}\mathsf{0})).$
In order to get a contradiction, we further take any closed term $x$ such that $\mathrm{P}_{\mathfrak{A}}(ux)$ and $\mathrm{P}_{\mathfrak{A}}(u(\dot{\neg}x)).$
Then, Lemma~\ref{lem:simultaneousness_in_U} implies 
$\mathrm{P}^{\alpha}_{\mathfrak{A}}(ux)$ and $\mathrm{P}^{\alpha}_{\mathfrak{A}}(u(\dot{\neg}x)),$ and thus we also have the inconsistency $\mathrm{P}^{<\alpha}_{\mathfrak{A}}(x) \land \neg \mathrm{P}^{<\alpha}_{\mathfrak{A}}(x)$ by the definition of $\mathfrak{A}_1.$

\item[Structural properties of $\mathrm{U}$.]
We consider the axiom $(\mathsf{Lim})$:
\[
\mathrm{C}^{\star}(f) \to \mathrm{U}^{\star}(\mathsf{l}f) \land (f \sqsubset \mathsf{l}f)^{\star}.
\]
We assume $\mathrm{C}^{\star}(f),$ then $\mathrm{U}^{\star}(\mathsf{l}f)$ is clear from the definition. Thus, we show the translation of $f \sqsubset \mathsf{l}f$:
\begin{enumerate}
\item $\forall x. \ \mathrm{P}_{\mathfrak{A}}(fx) \to \mathrm{P}_{\mathfrak{A}}(\mathsf{l}f(fx)),$
\item $\forall x. \ \mathrm{P}_{\mathfrak{A}}(\dot{\neg}(fx)) \to \mathrm{P}_{\mathfrak{A}}(\mathsf{l}f(\dot{\neg}(fx))).$
\end{enumerate}
As to the first item, we take any $x$ satisfying $\mathrm{P}_{\mathfrak{A}}(fx)$, and show $\mathrm{P}_{\mathfrak{A}}(\mathsf{l}f(fx))$.
By the assumption $\mathrm{C}^{\star}(f),$ Lemma~\ref{lem:simultaneousness_in_U} implies that there exists a least $\alpha$ such that $\mathrm{P}^{\alpha}_{\mathfrak{A}}(\mathsf{l}f(\dot{\mathsf{N}}\mathsf{0}))$. Thus, by the definition of $\mathfrak{A},$ it follows that $\forall y. \ \mathrm{P}^{<\alpha}_{\mathfrak{A}}(fy) \lor \mathrm{P}^{<\alpha}_{\mathfrak{A}}(\dot{\neg}(fy)).$ By the consistency of $\mathrm{P}_{\mathfrak{A}}$, we obtain $\mathrm{P}^{<\alpha}_{\mathfrak{A}}(fx).$ Therefore, we get $\mathrm{P}^{\alpha}_{\mathfrak{A}}(\mathsf{l}f(fx))$, and thus $\mathrm{P}_{\mathfrak{A}}(\mathsf{l}f(fx))$. The second item is similar.

The other structural axioms are similarly proved.

\item[Propositon induction.]
The proof is almost the same as that of Proposition~\ref{prop:model_of_KFUPI}. \qed
\end{description}
\end{pf}

\begin{thm}\label{thm:strength_of_KFUPI}
$\mathsf{KFUPI}$ and $\mathsf{T}_0$ are proof-theoretically equivalent.
\end{thm}

\begin{pf}\label{pf:thm:strength_of_KFUPI}
The lower bound of $\mathsf{KFUPI}$ is given by Theorem~\ref{thm:lower-bound_of_KFUPI}.
As for the upper bound, Lemma~\ref{lem:interpretation_of_KFUPI} showed that $\mathsf{KFUPI}$ is proof-theoretically reducible to $\mathsf{FID}([\mathsf{POS,QF}])$, which is known to be proof-theoretically equivalent to $\mathsf{T}_0$ \cite{jager2001first}. \qed
\end{pf}

%%%%%%%%%%%%%%%%%%%%%%%%%%%%%%%%%%%%%%%%%%%%%%%%%%%%%%%%%%%%%%%%%%%%%%%%%%%%%%%%%%%%%%

\section{Frege structure by the Aczel-Feferman schema}\label{sec:FS_by_AF}
In the present section, we consider another theory, $\mathsf{PT}$, based on Aczel's original Frege structure \cite{aczel1980frege}.  
Since the notion of propositions in $\mathsf{PT}$ is different from that in $\mathsf{KF}$, we need to change the definition of proposition induction slightly.
Nevertheless, almost the same proof-theoretic analysis as for $\mathsf{KFUPI}$ can be applied to $\mathsf{PT}$, and hence we obtain the theory $\mathsf{PTUPI}$, which is proof-theoretically equivalent to $\mathsf{T}_0$ (Theorem~\ref{thm:strength_of_PTUPI}).
\subsection{System $\mathsf{PT}$ and universes}\label{subsec:system_PT}
Whereas $\mathsf{KF}$ is based on Strong Kleene logic,
Aczel's original Frege structure \cite{aczel1980frege} is essentially based on \emph{Aczel--Feferman} logic, a variant of Weak Kleene logic which has the following truth tables:

\begin{center}
\begin{tabular}{|c|c|c|c|}
\hline
$\neg$        &     \\ \hline
$\mathrm{T}$  & $\mathrm{F}$    \\ \hline
$\mathrm{U}$ & $\mathrm{U}$  \\ \hline
$\mathrm{F}$  & $\mathrm{T}$    \\ \hline
\end{tabular} \ \ \ \ 
\begin{tabular}{|c|c|c|c|}
\hline
$\lor$        & $\mathrm{T}$  & $\mathrm{U}$ & $\mathrm{F}$  \\ \hline
$\mathrm{T}$  & $\mathrm{T}$  & $\mathrm{U}$ & $\mathrm{T}$  \\ \hline
$\mathrm{U}$ & $\mathrm{U}$ & $\mathrm{U}$ & $\mathrm{U}$ \\ \hline
$\mathrm{F}$  & $\mathrm{T}$  & $\mathrm{U}$ & $\mathrm{F}$  \\ \hline
\end{tabular} \ \ \ \ 
\begin{tabular}{|c|c|c|c|}
\hline
$\land$        & $\mathrm{T}$  & $\mathrm{U}$ & $\mathrm{F}$  \\ \hline
$\mathrm{T}$  & $\mathrm{T}$  & $\mathrm{U}$ & $\mathrm{F}$  \\ \hline
$\mathrm{U}$ & $\mathrm{U}$ & $\mathrm{U}$ & $\mathrm{U}$ \\ \hline
$\mathrm{F}$  & $\mathrm{F}$  & $\mathrm{U}$ & $\mathrm{F}$  \\ \hline
\end{tabular} \ \ \ \ 
%\end{table}
%\begin{table}[]
\begin{tabular}{|c|c|c|c|}
\hline
$\to$        & $\mathrm{T}$  & $\mathrm{U}$ & $\mathrm{F}$  \\ \hline
$\mathrm{T}$  & $\mathrm{T}$  & $\mathrm{U}$ & $\mathrm{F}$  \\ \hline
$\mathrm{U}$ & $\mathrm{U}$ & $\mathrm{U}$ & $\mathrm{U}$ \\ \hline
$\mathrm{F}$  & $\mathrm{T}$  & $\mathrm{T}$  & $\mathrm{T}$  \\ \hline
\end{tabular}
\end{center}
Note that while $A \lor B$ is definable as $\neg (\neg A \land \neg B) $, $\to$ cannot be defined by $\neg$ and $\land$.

Based on the truth table, the theory $\mathsf{PT}$ (proposition and truth) is defined as follows.
\begin{defn}[cf. \cite{feferman2008axioms}]\label{def:PT}
The $\mathcal{L}_{FS}$-theory $\mathsf{PT}$ has $\mathsf{TON}$ and (the universal closure of) the following axioms:
\begin{description}
\item[Compositional axioms.] \
\begin{itemize}
\item $\mathrm{T}(x \dot{=} y) \leftrightarrow x = y$

\item $\mathrm{T}(\dot{\neg}(x \dot{=} y)) \leftrightarrow x \neq y$

\item $\mathrm{T}(\dot{\mathsf{N}}x) \leftrightarrow \mathrm{N}(x)$

\item $\mathrm{T}(\dot{\neg}(\dot{\mathsf{N}}x)) \leftrightarrow \neg \mathrm{N}(x)$

\item $\mathrm{T} (\dot{\neg} (\dot{\neg} x)) \leftrightarrow \mathrm{T}(x)$

\item $\mathrm{T} (x \dot{\land} y) \leftrightarrow \mathrm{T}(x) \land \mathrm{T} (y)$

\item $\mathrm{T} (\dot{\neg} (x \dot{\land} y)) \leftrightarrow [\mathrm{T}(\dot{\neg}x) \land \mathrm{T} (\dot{\neg}y)] \lor [\mathrm{T}(\dot{\neg}x) \land \mathrm{T} (y)] \lor [\mathrm{T}(x) \land \mathrm{T} (\dot{\neg}y)]$

\item $\mathrm{T} (x \dot{\to} y) \leftrightarrow [\mathrm{T}(x) \land \mathrm{T}(y)] \lor \mathrm{T}(\dot{\neg}x)$

\item $\mathrm{T}(\dot{\neg} (x \dot{\to} y)) \leftrightarrow [\mathrm{T}(x) \land \mathrm{T}(\dot{\neg}y)]$

\item $\mathrm{T}(\dot{\forall} f) \leftrightarrow \forall x. \ \mathrm{T}(fx)$

\item $\mathrm{T}(\dot{\neg}(\dot{\forall} f)) \leftrightarrow [\forall x. \ \mathrm{T}(fx) \lor \mathrm{T}(\dot{\neg}(fx))]  \land \exists x. \ \mathrm{T}(\dot{\neg}(fx))$
\end{itemize}
\item[Consistency.] \
\begin{itemize}
\item $\neg [\mathrm{T}(x) \land \mathrm{T}(\dot{\neg}x)]$
\end{itemize}
\end{description}
\end{defn}

In Aczel's Frege structure, propositions can be naturally characterised inductively. %\footnote{In fact, Aczel introduced $\mathrm{P}(x)$ as a primitive symbol, and thus Lemma~\ref{lem:prop_in_PT} was given as axioms, not derived.}
%We remark that Aczel's originally considers the notion of $proposition$ in addition to truth; similarly to the truth condition, whether a given sentence is a proposition or not is computed from those of simpler sentences. 
Recall the notation $\mathrm{P}(x) := \mathrm{T}(x) \lor \mathrm{T}(\dot{\neg}x)$. Then, the following holds.

\begin{lem}[cf. \cite{aczel1980frege}]\label{lem:prop_in_PT} 
$\mathsf{PT}$ derives the following:
\begin{enumerate}
\item $\mathrm{P}(x \dot{=} y) \land \mathrm{P}(\dot{\mathsf{N}}x),$
\item $\mathrm{P}(\dot{\neg}x) \leftrightarrow \mathrm{P}(x),$
\item $\mathrm{P}(x \dot{\land} y) \leftrightarrow \mathrm{P}(x) \land \mathrm{P}(y),$
\item $\mathrm{P}(x \dot{\to} y) \leftrightarrow \mathrm{P}(x) \land (\mathrm{T}(x) \to \mathrm{P}(y)),$
%\item $[ \mathrm{P}(x) \land ( \mathrm{T}(x) \to \mathrm{T}(y) ) \to \mathrm{T}(x \dot{\to} y)] \land [ \mathrm{T}(x) \land \mathrm{T}(x \dot{\to} y) \to \mathrm{T}(y) ],$
\item $\mathrm{P}(\dot{\forall}f) \leftrightarrow \forall x. \ \mathrm{P}(fx).$
\end{enumerate}
\end{lem}

The proof-theoretic analysis of $\mathsf{PT}$ is found in, e.g., \cite{beeson1985foundations,cantini2016truth,hayashi1995new,fujimoto2010relative}:
\begin{fact}\label{fact:strength_of_PT}
$\mathsf{PT}$ and $\mathsf{EM}$ have the same $\mathcal{L}$-theorems.
\end{fact}

\begin{defn}\label{def:PTU}
The $\mathcal{L}_{FS}$-system $\mathsf{PTU}$ consists of $\mathsf{TON}$ and the following axioms:
\begin{description}
\item[Compositional axioms in $\mathrm{U}$.] \
\begin{itemize}
\item $\mathrm{U}(u) \to \forall x, y. \ \mathrm{T}(x \dot{=} y) \leftrightarrow x = y$
\item $\mathrm{U}(u) \to \forall x, y. \ \mathrm{T}(x \dot{\neq} y) \leftrightarrow x \neq y$
\item $\mathrm{U}(u) \to \forall x. \ \mathrm{T}(\dot{\mathsf{N}}x) \leftrightarrow \mathrm{N}(x)$
\item $\mathrm{U}(u) \to \forall x. \ \mathrm{T}(\dot{\neg}(\dot{\mathsf{N}}x)) \leftrightarrow \neg \mathrm{N}(x)$

\item $\mathrm{U}(u) \to \forall x. \ \mathrm{T} (\dot{\neg} (\dot{\neg} x)) \leftrightarrow \mathrm{T}(x)$

%\item $\mathrm{U}(u) \to
%\forall x, y. \mathrm{T} (x \dot{\lor} y) \leftrightarrow [\mathrm{T}(x) \land \mathrm{T} (y)] \lor [\mathrm{T}(\dot{\neg}x) \land \mathrm{T} (y)] \lor [\mathrm{T}(x) \land \mathrm{T} (\dot{\neg}y)]
%$

%\item $\mathrm{U}(u) \to \forall x, y. \mathrm{T} (\dot{\neg}(x \dot{\lor} y)) \leftrightarrow \mathrm{T}(\dot{\neg}x) \land \mathrm{T} (\dot{\neg}y)$

\item $\mathrm{U}(u) \to \forall x, y. \ \mathrm{T} (x \dot{\land} y) \leftrightarrow \mathrm{T}(x) \land \mathrm{T} (y)$

\item $\mathrm{U}(u) \to
\forall x, y. \ \mathrm{T} (\dot{\neg} (x \dot{\land} y)) \leftrightarrow [\mathrm{T}(\dot{\neg}x) \land \mathrm{T} (\dot{\neg}y)] \lor [\mathrm{T}(\dot{\neg}x) \land \mathrm{T} (y)] \lor [\mathrm{T}(x) \land \mathrm{T} (\dot{\neg}y)]
$

\item $\mathrm{U}(u) \to \forall x, y. \ \mathrm{T} (x \dot{\to} y) \leftrightarrow [\mathrm{T}(x) \land \mathrm{T}(y)] \lor \mathrm{T}(\dot{\neg}x)$

\item $\mathrm{U}(u) \to \forall x, y. \ \mathrm{T}(\dot{\neg} (x \dot{\to} y)) \leftrightarrow \mathrm{T}(x) \land \mathrm{T}(\dot{\neg}y)$

\item $\mathrm{U}(u) \to \forall f. \ \mathrm{T}(\dot{\forall} f) \leftrightarrow \forall x. \ \mathrm{T}(fx)$

\item $\mathrm{U}(u) \to \forall f. \ \mathrm{T}(\dot{\neg}(\dot{\forall} f)) \leftrightarrow [\forall x. \  \mathrm{T}(fx) \lor \mathrm{T}(\dot{\neg}(fx)) ] \land \exists x. \ \mathrm{T}(\dot{\neg}(fx))$
\end{itemize}

\item[Consistency in $\mathrm{U}$.] \
\begin{itemize}
\item $\mathrm{U}(u) \to \forall x. \ \neg [\mathrm{T}(x) \land \mathrm{T}(\dot{\neg}x)]$
\end{itemize}

\item[Structural properties of $\mathrm{U}$.] \
\begin{itemize}
\item $\mathrm{U}(u) \to \mathrm{C}(u)$
\item $\mathrm{U}(u) \to \forall x. \ \mathrm{T}(ux) \to \mathrm{T}(x)$
\item $\mathrm{C}(f) \to \mathrm{U}(\mathsf{l}f) \land f \sqsubset \mathsf{l}f$
\end{itemize}
\end{description}
\end{defn}

As far as the author knows, the system $\mathsf{PTU}$ has not been explicitly formulated in the literature, though Fujimoto \cite[pp.~933-935]{fujimoto_2011} formulated and analysed transfinite iterations of the truth theory $\mathsf{DT}$, which is essentially $\mathsf{PT}$ formulated over Peano arithmetic. Fujimoto's proof can be adapted to the proof-theoretic analysis of $\mathsf{PTU}.$
As a result, we obtain the following:
\begin{corr}\label{corr:strength_of_PTU}
$\mathsf{PTU}$ and $\mathsf{EMU}$ have the same $\mathcal{L}$-theorems.
\end{corr}

\subsection{$\mathsf{PTU}$ with proposition induction}\label{subsec:PTU_with_PI}
In order to define the proposition induction for $\mathsf{PT},$
we now use another operator $\mathscr{A}^{\mathsf{PT}}$, which is based on the inductive characterisation of propositions in Lemma~\ref{lem:prop_in_PT}. For each $\mathcal{L}_{FS}$-formula $B$ and a free variable $x$, 
the formula $\mathscr{A}^{\mathsf{PT}} (B(\bullet), x)$ is the disjunction of the following:
\begin{enumerate}
\item $\exists y, z. \ x = (y \dot{=} z) \land y = z$
%\item $\exists y, z. \ x = \dot{\neg} (y \dot{=} z) \land y \neq z$
\item $\exists y. \ x = (\dot{\mathsf{N}} y) \land \mathrm{N}(y)$
%\item $\exists y. \ x = \dot{\neg} (\dot{\mathsf{N}} y) \land \neg \mathrm{N}(y)$
\item $\exists y. \ x = \dot{\neg} y \land   B(y)$
\item $\exists y, z. \ x = y \dot{\land} z \land B(y) \land B(z)$
\item $\exists y, z. \ x = y \dot{\to} z \land B(y) \land [\mathrm{T} (y) \to B(z)]$
\item $\exists g. \ x = \dot{\forall} g \land \forall y. \ B(gy)$
\item $\exists f,y. \ x=\mathsf{l}fy \land \forall z. \ B(fz)$
\end{enumerate}

Then, the proposition induction schema $(\mathsf{PI}^{\mathsf{PT}})$ is given by:
\[
[\forall x. \ \mathscr{A}^{\mathsf{PT}}(B(\bullet),x) \to B(x)] \to \forall x. \ \mathrm{P}(x) \to B(x),
\]
for all $\mathcal{L}_{FS}$-formulas $B(x).$

\begin{defn}\label{def:PTUPI}
The $\mathcal{L}_{FS}$-theory $\mathsf{PTUPI}$ is $\mathsf{PTU}$ with the schemata $(\mathsf{UG})$ and $(\mathsf{PI}^{\mathsf{PT}})$.  
\end{defn}

As for the proof-theoretic strength of $\mathsf{PTUPI},$ exactly the same lower-bound proof of Section~\ref{subsec:lower-bound_of_KFUPI} can be given in $\mathsf{PTUPI}$; thus, we have $\mathsf{T}_0 \leq \mathsf{PTUPI}.$
Moreover, Kahle's model construction in Section~\ref{subsec:Kahle's_model_for_KFU} and, thus, the upper-bound proof of $\mathsf{KFUPI}$ in Section~\ref{subsec:upper-bound_of_KFUPI} are easily modified for $\mathsf{PTUPI}.$ As a result, we also obtain the upper bound.

\begin{thm}\label{thm:strength_of_PTUPI}
$\mathsf{PTUPI}$ and $\mathsf{T}_0$ are proof-theoretically equivalent.
\end{thm}

%%%%%%%%%%%%%%%%%%%%%%%%%%%%%%%%%%%%%%%%%%%%%%%%%%%%%%%%%%%%%%%%%%%%%%%%%%%%%%%%%%%%%%

\section{Extension by least universes}\label{sec:extension_by_LU}
In Sections~\ref{sec:extension_by_PI} and~\ref{sec:FS_by_AF}, we extended theories of Frege structure by using the induction principle on propositions, in analogy with name induction in \cite{jager2001universes}. In the same paper, J\"{a}ger, Kahle and Studer further proposed the idea of \emph{least universes}, based on which they formulated the theory $\mathsf{LUN}$. As $\mathsf{LUN}$ is proof-theoretically equivalent to $\mathsf{T}_0$ \cite[Conclusion~25]{jager2001universes}, we can expect to obtain a strong theory of Frege structure by requiring some kind of leastness on universes.
In fact, Burgess \cite{burgess2009friedman} proposed a truth theory $\mathsf{KF}_{\mu}$ with impredicative strength, in which the truth predicate is intended to denote the \emph{least} fixed point of the Kripke operator. Therefore, in this section, we aim to extend $\mathsf{KF}$ by least universes.% which work like the truth predicate of $\mathsf{KF}_{\mu}$.

In Kahle's model $\mathcal{M}_{\mathsf{KF}}$ in Section~\ref{subsec:Kahle's_model_for_KFU}, each universe of the form $\mathsf{l}f$ reflects a class $f$ and is a fixed point of the Kripke operator.
Thus, one could naturally introduce leastness by defining a universe $\mathsf{l}f$ as the least set that reflects $f$ and is closed under the Kripke operator.\footnote{In $\mathsf{LUN}$, the least universe $\mathsf{lt}(a)$ for a name $a$ is defined to be the least set that contains $a$ and is closed under the set constructions of $\mathsf{EMU}.$ } 
However, we can easily observe that universes so defined violate natural structural properties in Remark~\ref{rmk:structure_of_M_KF}, similar to 
\cite[Theorem~14]{jager2001universes}. 
This means that  such universes are incompatible with $\mathcal{M}_{\mathsf{KF}}$.
To make matters worse, these universes are not generally fixed points of the Kripke operator, so they do not sufficiently serve as truth predicates. This definition of leastness does not capture universes in $\mathcal{M}_{\mathsf{KF}}$  because universes in $\mathcal{M}_{\mathsf{KF}}$ may contain other small universes. 
Therefore, in order to define a least universe $\mathsf{l}f$ compatible with $\mathcal{M}_{\mathsf{KF}}$, we also need to accommodate smaller universes $\mathsf{l}g$ such that $\mathsf{l}g \sqsubset \mathsf{l}f$ holds. 
Taking these into consideration, we shall define the least universe $\mathsf{l}f$ for a class $f$ as the least set that reflects $f$ and other smaller universes, and is closed under the Kripke operator. In addition, to get a strong system, we characterise least universes in terms of \emph{proposition} as in Lemma~\ref{lem:prop_in_KF}, rather than \emph{truth}.

Let $\mathrm{C}_{\mathsf{l}f}(x) :=  \forall y. \ \mathrm{P}_{\mathsf{l}f}(xy) :=\forall y. \ \mathrm{T}(\mathsf{l}f(xy)) \lor \mathrm{T}(\mathsf{l}f(\dot{\neg}(xy))).$ Thus, this means that $x$ is a class within $\mathsf{l}f.$
Then, for each term $f$, $\mathcal{L}_{FS}$-formula $B$ and a free variable $x$, 
the formula $\mathscr{A}^{\mathsf{LU}} (f, B(\bullet), x)$ is defined to be the disjunction of the following:
\begin{enumerate}
\item $\exists y, z. \ x = (y \dot{=} z) \land y = z$
\item $\exists y. \ x = (\dot{\mathsf{N}} y) \land \mathrm{N}(y)$

\item $\exists y. \ x = \dot{\neg} y \land B(y)$

\item $\exists y, z. \ x = (y \dot{\land} z) \land \{[B(y) \land B(z)] \lor [B(y) \land \mathrm{T}(\mathsf{l}f(\dot{\neg}y))] \lor [B(z) \land \mathrm{T}(\mathsf{l}f(\dot{\neg}z))] \}$

\item $\exists y, z. \ x = (y \dot{\to} z) \land \{[B(y) \land B(z)] \lor [B(y) \land \mathrm{T}(\mathsf{l}f(\dot{\neg}y))] \lor [B(z) \land \mathrm{T}(\mathsf{l}fz)] \}$

\item $\exists g. \ x = (\dot{\forall} g) \land  \{[\forall y. \ B(gy)] \lor \exists y. \ B(gy) \land \mathrm{T}(\mathsf{l}f(\dot{\neg}(gy)) \}$

\item $\exists y. \ x = fy$
\item $\exists g,y. \ x = \mathsf{l}gy \land \mathrm{P}_{\mathsf{l}f}(x)$
\end{enumerate}

The least universe schema $(\mathsf{LU})$ is given by:
\begin{align}
[\mathrm{C}(f) \land \forall x. \ \mathscr{A}^{\mathsf{LU}}(f, B(\bullet),x) \to B(x)] \to \forall x. \ \mathrm{P}_{\mathsf{l}f}(x) \to B(x), \tag{$\mathsf{LU}$}\label{LU}
\end{align}
for all $\mathcal{L}_{FS}$-formulas $B(x).$

\begin{defn}\label{def:KFLU}
The $\mathcal{L}_{FS}$-theory $\mathsf{KFLU}$ is $\mathsf{KFU}$ with the schemata $(\mathsf{UG})$ and $(\mathsf{LU})$.  
\end{defn}

As an immediate consequence of the definition, we can show that each least universe is closed under the operator $\mathscr{A}^{\mathsf{LU}}:$
\begin{corr}\label{corr:KFLU_is_A'-closed}
$\mathsf{KFLU} \vdash \mathrm{C}(f) \to \forall x. \ \mathscr{A}^{\mathsf{LU}}(f,
\mathrm{P}_{\mathsf{l}f}(\bullet),x) \to \mathrm{P}_{\mathsf{l}f}(x).$
\end{corr}

For the lower bound of $\mathsf{KFLU},$ we can interpret $\mathsf{T}_0$ into $\mathsf{KFLU}$ using the same translation as for $\mathsf{KFUPI}$. Moreover, the proof is almost the same as for $\mathsf{KFUPI},$ except that Lemmata~\ref{lem:BA_is_A-closure} and \ref{lem:IG2_in_KFUPI} for the derivation of $(\mathsf{IG}.2)'$ are now replaced by the following lemmata:
\begin{lem}\label{lem:BA_is_A-closure_for_KFLU}
For any $\mathcal{L}_{FS}$-formula $A$, recall that $B_{A}(u,v,w)$ is the formula $\forall x. \
\mathrm{T}(\mathsf{acc}(u,v)x) \to F_A(u,v,w,x)$, 
where $F_A(u,v,w,x)$ is the conjunction of the following formulas:
\begin{enumerate}
\item $w = t(u,v,x) \ \to \ A(x),$

\item $\forall z. \ w = \{(\oplus(u, v) (1, (z,x))) \dot{\to} t(u,v,z)\} \to \mathrm{T}(\oplus(u, v) (1, (z,x))) \to A(z),$

\item $\forall z. \ w = \{\oplus(u, v) (0,z) \dot{\to} (\oplus(u, v) (1, (z,x)) \dot{\to} t(u,v,z))\} 
\to \mathrm{T}(\oplus(u, v) (0,z)) \to \mathrm{T}(\oplus(u, v) (1, (z,x))) \to A(z).$
\end{enumerate}

Then, $\mathsf{KFLU}$ derives the following:
\[
\mathrm{C}(u) \land \mathrm{C}(v) \land \mathsf{Closed}' (u,v,A(\bullet)) \to \forall w. \ \mathscr{A}^{\mathsf{LU}}(\oplus(u,v),B_{A}(u,v, \bullet), w) \to B_{A}(u,v,w). 
\]
\end{lem}

\begin{lem}\label{lem:IG2_in_KFLU}
$\mathsf{KFUPI}$ derives the following:
\[
\mathrm{C}(a) \land \mathrm{C}(b) \land Closed(a,b,A) \to \forall x. \ \mathrm{T}(\mathsf{acc}(a,b)x) \to A(x).
\]
\end{lem}

\begin{pf}
We assume $\mathrm{C}(a), \mathrm{C}(b), Closed(a,b,A)$ and $\mathrm{T}(\mathsf{acc}(a,b)c)$ for an arbitrary $c$. Then, we have to derive $A(c)$.
From the assumptions and Lemma~\ref{lem:BA_is_A-closure_for_KFLU}, we obtain $\mathrm{C}(\mathsf{l}(\oplus(a,b)))$ and $\forall w. \ \mathscr{A}^{\mathsf{LU}}(\oplus(u, v), B_{A}(u,v, \bullet), w) \to B_{A}(u,v,w).$
Thus,  $(\mathsf{LU})$ yields:
\[
\forall x. \ \mathrm{P}_{\mathsf{l}(\oplus(a,b))}(x) \to B_{A}(a,b,x).
\]
Letting $x := t(a,b,c),$ we want to show $\mathrm{P}_{\mathsf{l}(\oplus(a,b))}(t(a,b,c)),$ 
from which we have $B_{A}(a,b,t(a,b,c))$, and thus we can immediately obtain the conclusion $A(c)$ by the assumption $\mathrm{T}(\mathsf{acc}(a,b)c)$ and the definition of $B_{A}$. 

In order to prove $\mathrm{P}_{\mathsf{l}(\oplus(a,b))}(t(a,b,c)),$ it suffices to derive the following:
\begin{align}
\forall y. \ \mathrm{T}(\oplus(a, b) (0,y)) \to \mathrm{T}(\oplus(a, b) (1, (y,c))) \to \mathrm{T}(\mathsf{l}(\oplus(a,b))(t(a,b,y))). \label{pf_1:lem:IG2_in_KFLU}
\end{align}
In fact, this formula implies $\mathrm{T}(\mathsf{l}(\oplus(a,b))(t(a,b,c)) )$ and $\mathrm{P}_{\mathsf{l}(\oplus(a,b))}(t(a,b,c) ),$ because $\oplus(a, b)$ is a class within $\mathsf{l}(\oplus(a,b))$.

\  However,  the formula (\ref{pf_1:lem:IG2_in_KFLU}) follows from the assumption $\mathrm{T}(\mathsf{acc}(a,b)c)$. \qed
\end{pf}

As for the upper bound, we can interpret $\mathsf{KFLU}$ into $\mathsf{FID}([\mathsf{POS,QF}])$ in a similar manner as for $\mathsf{KFUPI}.$
So, we only observe that the schema $(\mathsf{KFLU})$ is satisfied in $\mathcal{M}_{\mathsf{KF}}.$

\begin{prop}\label{prop:model_of_KFLU}
$\mathcal{M}_{\mathsf{KF}} \models \mathsf{KFLU}$.
\end{prop}

\begin{pf}\label{pf:prop:model_of_KFLU}
By Proposition~\ref{prop:model_of_KFUPI}, we concentrate on the schema (\ref{LU}):
\[
[\mathrm{C}(f) \land \forall x. \ \mathscr{A}^{\mathsf{LU}}(f, B(\bullet),x) \to B(x)] \to \forall x. \ \mathrm{P}_{\mathsf{l}f}(x) \to B(x).
\]
So, we take any class $f$ and assume that $\mathcal{M}_{\mathsf{KF}} \models \forall x. \ \mathscr{A}^{\mathsf{LU}}(f,B(\bullet),x) \to B(x)$. 
Now, we prove by induction on $\alpha$ that if $x \in Z_{\alpha}$ or $\dot{\neg}x \in Z_{\alpha}$, then $B(x)$.
%Then, taking any $x$, we need to show
%$\mathcal{M}_{\mathsf{KF}} \models \mathrm{P}_{\mathsf{l}f}(x) \to B(x).$ 
%If $\mathcal{M}_{\mathsf{KF}} \models \mathrm{P}_{\mathsf{l}f}(x)$, 
%then by the definition of $\mathcal{M}_{\mathsf{KF}},$ there exists a least $\alpha$ 
%such that $x \in Z_{\alpha}$ or $\dot{\neg}x \in Z_{\alpha}$ holds.
%Now, we define the rank $\mathrm{rk}(x)$ of $x$ to be this $\alpha$.
%Then, we shall prove, by induction on $\alpha,$ that if $\mathcal{M}_{\mathsf{KF}} \models \mathrm{P}_{\mathsf{l}f}(x)$ and $\mathrm{rk}(x)= \alpha$, then $\mathcal{M}_{\mathsf{KF}} \models B(x).$ 

For example, consider the case where $x$ is of the form $\mathsf{l}gy.$ 
%Then, we can easily verify from the definition of $\mathcal{M}_{\mathsf{KF}}$ that 
%$\mathsf{l}g$ is a class within $\mathsf{l}f,$ that is, $\mathcal{M}_{\mathsf{KF}} \models \mathrm{C}_{\mathsf{l}f}(\mathsf{l}g).$ 
Then, since $B$ is $\mathscr{A}^{\mathsf{LU}}$-closed, it follows that $B(\mathsf{l}gy).$ The other cases are shown by  the subinduction on the construction of $Z_{\alpha}$ (cf. Proposition~\ref{prop:model_of_KFUPI}).
\qed
\end{pf}

In conclusion, we obtain the proof-theoretic strength of $\mathsf{KFLU}:$
\begin{thm}\label{thm:strength_of_KFLU}
$\mathsf{KFLU}$ and $\mathsf{T}_0$ are proof-theoretically equivalent.
\end{thm}

Similarly to $\mathsf{KFLU}$, we can also consider least universes for $\mathsf{PT}$. As for the Aczel-Feferman schema, we can naturally characterise least universes in terms of  
$truth$ only; thus, we define an operator $\mathscr{A}^{\mathsf{PTLU}}$ as follows:
 For each term $f$, $\mathcal{L}_{\mathrm{T}}$-formula $B$ and a free variable $x$, 
the formula $\mathscr{A}^{\mathsf{PTLU}} (f, B(\bullet), x)$ is the disjunction of the following:
\begin{enumerate}
\item $\exists y, z. \ x = (y \dot{=} z) \land y = z$
\item $\exists y, z. \ x = (\dot{\neg}(y \dot{=} z)) \land y \neq z$

\item $\exists y. \ x = (\dot{\mathsf{N}} y) \land \mathrm{N}(y)$
\item $\exists y. \ x = (\dot{\neg}(\dot{\mathsf{N}} y)) \land \neg\mathrm{N}(y)$

\item $\exists y. \ x = \dot{\neg}(\dot{\neg} y) \land B(y)$

\item $\exists y, z. \ x = (y \dot{\land} z) \land B(y) \land B(z)$
\item $\exists y, z. \ x = (\dot{\neg}(y \dot{\land} z)) \land \{ [B(\dot{\neg}y) \land B(\dot{\neg}z)] \lor [B(y) \land B(\dot{\neg}z)] \lor [B(\dot{\neg}y) \land B(z)] \}$

\item $\exists y, z. \ x = (y \dot{\to} z) \land [\mathrm{T}(\mathsf{l}fy) \lor \mathrm{T}(\mathsf{l}f(\dot{\neg}y))] \land [\mathrm{T}(\mathsf{l}f(y)) \to B(z)]$
\item $\exists y, z. \ x = (\dot{\neg}(y \dot{\to} z)) \land B(y) \land B(\dot{\neg}z)$

\item $\exists g. \ x = (\dot{\forall} g) \land \forall y. B(gy)$
\item $\exists g. \ x = (\dot{\neg}(\dot{\forall} g)) \land [\forall y. B(gy) \lor B(\dot{\neg}(gy))] \land \exists y. B(\dot{\neg}(gy))$

\item $\exists y. \ [x = fy \lor x = \dot{\neg}(fy)] \land \mathrm{T}(x)$
\item $\exists g,y. \ [x = \mathsf{l}gy \lor x = \dot{\neg}(\mathsf{l}gy) ] \land \mathrm{T}(\mathsf{l}fx)$
\end{enumerate}

The least universe schema $(\mathsf{LU}^{\mathsf{PT}})$ is given by:
\begin{align}
[\mathrm{C}(f) \land \forall x. \ \mathscr{A}^{\mathsf{PTLU}}(f, B(\bullet),x) \to B(x)] \to \forall x. \ \mathrm{T}(\mathsf{l}fx) \to B(x), \tag{$\mathsf{LU}^{\mathsf{PT}}$}\label{LU_PT}
\end{align}
for all $\mathcal{L}_{FS}$-formulas $B(x).$

Thus, the schema $(\mathsf{LU}^{\mathsf{PT}})$ informally says that $\mathrm{T}(\mathsf{l}fx)$ is the least truth set satisfying each clause of $\mathscr{A}^{\mathsf{PTLU}}$.

\begin{defn}\label{def:PTLU}
The $\mathcal{L}_{FS}$-theory $\mathsf{PTLU}$ is $\mathsf{PTU}$ with the schemata $(\mathsf{UG})$ and $(\mathsf{LU}^{\mathsf{PT}})$.  
\end{defn}

In exactly the same way as for $\mathsf{KFLU}$, we can determine the proof-theoretic strength of $\mathsf{PTLU}$:
\begin{thm}\label{thm:strength_of_PTLU}
$\mathsf{PTLU}$ and $\mathsf{T}_0$ are proof-theoretically equivalent.
\end{thm}

%%%%%%%%%%%%%%%%%%%%%%%%%%%%%%%%%%%%%%%%%%%%%%%%%%%%%%%%%%%%%%%%%%%%%%%%%%%%%%%%%%%%%%

\section{Frege structure by the supervaluation schema}\label{sec:FS_by_supervaluation}
In this section, we consider the supervaluational Frege structure. 
Kripke initially sketched a semantic theory of truth based on the supervaluation schema \cite[p.~711]{kripke1976outline}.
Based on Kripke's semantic definition, Cantini \cite{cantini1990theory} defined and studied a formal theory $\mathsf{VF}$ (van Fraassen) over Peano arithmetic. 
In particular, Cantini \cite{cantini1990theory} proved that $\mathsf{VF}$ is proof-theoretically equivalent to the theory $\mathsf{ID}_1$ of arithmetical monotone inductive definition. 
In \cite{cantini1996logical,kahle2001truth}, the system $\mathsf{VF}$ as a theory of Frege structure is formulated and found to be proof-theoretically equivalent to the one over Peano arithmetic.
For future research, Kahle \cite[p.~124]{kahle2001truth} suggested extending $\mathsf{VF}$ by adding universes, similar to
$\mathsf{KFU}$ and $\mathsf{PTU}$. 
Given that $\mathsf{KFU}$ is roughly a transfinite iteration of $\mathsf{KF}$, it is natural to expect $\mathsf{VF}$ with universes to have the strength of a transfinite iteration of $\mathsf{VF}$. 
Thus, Kahle conjectured that such a theory would have at least the strength of $\mathsf{ID}_{\alpha}$ for some ordinal $\alpha$
 (cf. \cite{kahle2001truth, kahle2015sets}).
 The purpose of this section is to implement the idea of universes for $\mathsf{VF}$ and to verify Kahle's conjecture. In particular, we will show that $\mathsf{VF}$ with universes is proof-theoretically equivalent to $\mathsf{T}_0.$

\subsection{System $\mathsf{VF}$ and universes}\label{subsec:system_VF}
In this section, we add a constant symbol $\dot{\mathrm{T}}$ to $\mathcal{L}_{FS}$ for technical reasons (see Remark~\ref{rmk:necessity_of_dotT}).
Following \cite{kahle2001truth, kahle2003universes}, we inductively define a corresponding term $\dot{A}$ for each formula $A$.
\begin{itemize}
\item $\dot{\overbrace{s = t}} := (s \dot{=} t),$ \ $\dot{\overbrace{\mathrm{N}(s)}} := \dot{\mathsf{N}}s,$ \ $\dot{\overbrace{\mathrm{T}(s)}} := \dot{\mathsf{T}}s;$
\item $\dot{\overbrace{\neg A}} := \dot{\neg} \dot{A}$,  \ $\dot{\overbrace{A \land B}} := \dot{A} \dot{\land} \dot{B}$, \ 
$\dot{\overbrace{A \to B}} := \dot{A} \dot{\to} \dot{B};$
\item $\dot{\overbrace{\forall x. \ A}} := \dot{\forall} (\lambda x. \dot{A}).$
\end{itemize}

Then, we define the $\mathcal{L}_{FS}$-theory $\mathsf{VF},$ the formulation of which is essentially based on \cite{cantini1990theory,cantini1996logical,kahle2001truth}.

\begin{defn}\label{def:VF}
The $\mathcal{L}_{FS}$-theory $\mathsf{VF}$ consists of $\mathsf{TON}$ and the following axioms:
\begin{description}
\item[$(\mathsf{T} \mhyph {\mathsf{Out}})$] $\mathrm{T} (\dot{A}) \to A$;
\item[$(\mathsf{T} \mhyph {\mathsf{Elem}})$] $P \to \mathrm{T}(\dot{P}),$ for any $\mathcal{L}$-literal formula $P$;
\item[$(\mathsf{T} \mhyph {\mathsf{Imp}})$] $\mathrm{T}(\dot{\overbrace{A \to B}}) \to (\mathrm{T}(\dot{A}) \to \mathrm{T}(\dot{B}))$;
\item[$(\mathsf{T} \mhyph {\mathsf{Univ}})$] $(\forall x. \ \mathrm{T}(\dot{A})) \to \mathrm{T} (\dot{\overbrace{\forall x. \ A}})$;\footnote{One might prefer the single axiom $(\forall x. \mathrm{T}(fx)) \to \mathrm{T}(\dot{\forall} f)$ similarly to $(\mathsf{K}_{\forall})$ in Definition~\ref{def:KF}.
The reason the author choses the schematic form is just to simplify the upper-bound proof in Section~\ref{subsec:upper-bound_of_VFU}. Nevertheless, the author believes that this single axiom does not affect the proof-theoretic strength of both $\mathsf{VF}$ and $\mathsf{VFU}$ (Definition~\ref{def:VFU}).
}
\item[$(\mathsf{T} \mhyph {\mathsf{Log}})$] $\mathrm{T}(\dot{C})$ for any logical theorem $C$;
\item[$(\mathsf{T}\mhyph {\mathsf{Cons}})$] $\neg [\mathrm{T}(x) \land \mathrm{T}(\dot{\neg}x)];$
\item[$(\mathsf{T} \mhyph \mathsf{Self})$] $[\mathrm{T}(x) \leftrightarrow \mathrm{T}(\dot{\mathsf{T}}x)] \ \land \ [\mathrm{T}(\dot{\neg}x) \leftrightarrow \mathrm{T}(\dot{\neg}(\dot{\mathsf{T}}x))].$
\end{description}
Here, $A$ and $B$ are any $\mathcal{L}_{FS}$-formulas.
\end{defn}

\begin{rmk}\label{rmk:necessity_of_dotT}
In Kahle's formulation of $\mathsf{VF}$ (called $\mathsf{SON}$ in \cite{kahle2001truth}), the constant $\dot{\mathrm{T}}$ is defined to be an identity function $\lambda x.x,$
so $\dot{\mathrm{T}}s$ is $\beta$-equivalent to $s$ itself.
Consequently, the axiom $(\mathsf{T} \mhyph \mathsf{Self})$ becomes trivial \cite[pp.~111-112]{kahle2001truth}. 
While this definition is not problematic in $\mathsf{KF}$ (Definition~\ref{def:KF}) and $\mathsf{PT}$ (Definition~\ref{def:PT}), it causes a contradiction in $\mathsf{VF}.$ In fact, for every $\mathcal{L}_{FS}$-sentence $A$, we have $\mathrm{T}(\dot{A} \dot{\lor} (\dot{\neg} \dot{A}))$ by $(\mathsf{T} \mhyph \mathsf{Log})$,
which is, according to Kahle's definition, equivalent to $\mathrm{T}((\dot{\mathsf{T}}\dot{A}) \dot{\lor} (\dot{\mathsf{T}}(\dot{\neg} \dot{A}))).$ Therefore, $(\mathsf{T} \mhyph \mathsf{Out})$ implies $\mathrm{T}(\dot{A}) \lor \mathrm{T}(\dot{\neg} \dot{A}).$
However, it is well known that $\mathsf{VF}$ with the schema $\mathrm{T}(\dot{A}) \lor \mathrm{T}(\dot{\neg} \dot{A})$ is inconsistent (see, e.g., \cite{friedman1987axiomatic}). That is the reason we explicitly introduced the constant $\dot{\mathsf{T}}$ and instead required the axiom $(\mathsf{T} \mhyph \mathsf{Self})$.
\end{rmk}

As is often said, $\mathsf{VF}$ has a non-compositional nature, and thus, it is not suitable for the inductive characterisation of proposition or truth, unlike $\mathsf{KF}$ and $\mathsf{PT}$.\footnote{Note that Stern's supervaluation-style truth \cite{stern2018supervaluation} is an attempt to overcome this difficulty of $\mathsf{VF}$.}
Alternatively, we show in the next subsection that simply adding the axiom $(\mathsf{Lim})$ (Definition~\ref{def:KFU}) to the universe-relative version of $\mathsf{VF}$ gives the same strength as $\mathsf{T}_0.$
So, analogously to $\mathsf{KFU}$ and $\mathsf{PTU},$ we propose the system $\mathsf{VFU}.$

Recall that $\mathrm{C}(f):\equiv \forall x. \ \mathrm{T}(fx) \lor \mathrm{T}(\dot{\neg}(fx)).$
\begin{defn}\label{def:VFU}
The $\mathcal{L}_{FS}$-theory $\mathsf{VFU}$ consists of $\mathsf{TON}$ and the following axioms:
\begin{description}
\item[$\mathsf{VF}$-axioms in $\mathrm{U}$.] \
\begin{description}
\item[$(\mathsf{U} \mhyph \mathsf{Out})$] $\mathrm{U}(u) \to [\mathrm{T} (u(\dot{A})) \to A]$;
\item[$(\mathsf{U} \mhyph \mathsf{Elem})$] $\mathrm{U}(u) \to [P \to \mathrm{T}(u(\dot{P}))],$ for any $\mathcal{L}$-literal formula $P$;
\item[$(\mathsf{U} \mhyph \mathsf{Imp})$] $\mathrm{U}(u) \to [\mathrm{T}(u(\dot{\overbrace{A \to B}})) \to \{ \mathrm{T}(u\dot{A}) \to \mathrm{T}(u\dot{B})\} ]$;
\item[$(\mathsf{U} \mhyph \mathsf{Univ})$] $\mathrm{U}(u) \to [\{ \forall x. \ \mathrm{T}(u(\dot{A}))\} \to \mathrm{T} (u(\dot{\overbrace{\forall x. \ A}}))]$;
\item[$(\mathsf{U} \mhyph \mathsf{Log})$] $\mathrm{U}(u) \to \mathrm{T}(u(\dot{C}))$ for any logical theorem $C$;
\item[$(\mathsf{U} \mhyph \mathsf{Cons})$] $\mathrm{U}(u) \to \forall x. \ \neg [\mathrm{T}(ux) \land \mathrm{T}(u(\dot{\neg}x))];$
\item[$(\mathsf{U} \mhyph \mathsf{Self})$] $\mathrm{U}(u) \to [\mathrm{T}(ux) \leftrightarrow \mathrm{T}(u(\dot{\mathsf{T}}x))] \ \land \ [\mathrm{T}(u(\dot{\neg}x)) \leftrightarrow \mathrm{T}(u(\dot{\neg}(\dot{\mathsf{T}}x)))];$
\end{description}
Here, $A$ and $B$ are any $\mathcal{L}_{FS}$-formulas.

\item[Structural properties of $\mathrm{U}$.] \
\begin{description}
\item[$(\mathsf{U} \mhyph \mathsf{Class})$] $\mathrm{U}(u) \to \mathrm{C}(u);$
\item[$(\mathsf{U} \mhyph \mathsf{True})$] $\mathrm{U}(u) \to \forall x. \ \mathrm{T}(ux) \to \mathrm{T}(x);$
\item[$(\mathsf{Lim})$] $\mathrm{C}(f) \to \mathrm{U}(\mathsf{l}f) \land f \sqsubset \mathsf{l}f.$
\end{description}
\end{description}
\end{defn}

Similarly to Fact~\ref{fact:KFU_implies_KF}, we can prove the following:
\begin{lem}\label{lem:VFU_implies_VF}
$\mathsf{VF}$ is a subtheory of $\mathsf{VFU}$.

\end{lem}

\subsection{Lower bound of $\mathsf{VFU}$}\label{subsec:lower-bound_of_VFU}
In this subsection, we determine the lower bound of $\mathsf{VFU}$.
The proof proceeds in a similar way as for $\mathsf{KFUPI}$ and $\mathsf{PTUPI}$. 
Thus, we define a translation $': \mathcal{L}_{EM} \to \mathcal{L}_{FS}$.
The interpretation of $\mathcal{L}$ is exactly the same as in Theorem~\ref{thm:lower-bound_of_KFUPI}.
Each generator of $\mathsf{T}_{0}$ is interpreted as a term of $\mathcal{L}_{FS}.$
In particular, the interpretation of the inductive generation $\mathsf{i}$ is essentially a generalisation of Cantini's original lower-bound proof of $\mathsf{VF}$ \cite{cantini1990theory, cantini1996logical}.

The following two lemmata are essentially by way of \cite[Lemma~59.2]{cantini1996logical}.

\begin{lem}\label{lem:T-POS_in_VFU}
An $\mathcal{L}_{FS}$-formula $A$ is $\mathrm{T}$-$positive$ if each truth predicate $\mathrm{T}$ in $A$ occurs only positively. 
Then, for every $\mathrm{T}$-positive $\mathcal{L}_{FS}$-formula $A$,
\begin{center}
$\mathsf{VFU} \vdash \mathrm{T}(\dot{A}) \leftrightarrow A.$
\end{center}
\end{lem}

\begin{lem}\label{lem:compositionality_in_VFU}
In $\mathsf{VFU},$ the following are derivable:
\begin{enumerate}
\item $\mathrm{T}(\dot{A}) \land \mathrm{T}(\dot{B}) \leftrightarrow \mathrm{T}(\dot{\overbrace{A \land B}});$
\item $\mathrm{T}(\dot{A}) \lor \mathrm{T}(\dot{B}) \to \mathrm{T}(\dot{\overbrace{A \lor B}});$
\item $(\mathrm{T}(\dot{A}) \lor \mathrm{T}(\dot{\overbrace{\neg A}})) \to [\{ \mathrm{T}(\dot{A}) \to \mathrm{T}(\dot{B}) \} \leftrightarrow \mathrm{T}(\dot{\overbrace{A \to B}})];$
\item $[\forall x. \ \mathrm{T}(\dot{A})] \leftrightarrow \mathrm{T}(\dot{\overbrace{\forall x A}}) ;$
\item $[\exists x. \ \mathrm{T}(\dot{A})] \to \mathrm{T}(\dot{\overbrace{\exists x A}});$
\item $[\mathrm{T}(a) \lor \mathrm{T}(\dot{\neg}a)] \to [\neg \mathrm{T}(\dot{\mathrm{T}}a) \leftrightarrow \mathrm{T}(\dot{\neg}(\dot{\mathrm{T}}a))].$
\end{enumerate}
\end{lem}

Using these lemmata, $\mathsf{VFU}$ can interpret each generator except $\mathsf{i}$. 
Here, as an example, we only deal with join $\mathsf{j}$, but the other generators can be similarly treated.

\begin{lem}\label{lem:interpretation_of_join_VFU}
Let a term $\underline{\mathsf{j}}$ be such that $\underline{\mathsf{j}}(x,f) = \lambda z. \dot{\overbrace{\exists v, w. \ z = (v,w) \land \mathrm{T}  \mathit{(x v)}  \land \mathrm{T} \mathit{(fvw)}}}.$
Then, $\mathsf{VFU} \vdash (\mathrm{join})'$, that is,
\[
\mathsf{VFU} \vdash \mathrm{C}(x) \land [\forall y. \ \mathrm{T}(xy) \to \mathrm{C}(fy)] \to \mathrm{C}(\underline{\mathsf{j}}(x,f)) \land \Sigma^{'}(x,f, \underline{\mathsf{j}} (x,f)),
\]
where, $\Sigma^{'}(x,f, \underline{\mathsf{j}} (x,f))$ is the following:
\[
\forall z. \ \mathrm{T}( \underline{\mathsf{j}} (x,f)z) \leftrightarrow \exists v, w. \ z = (v,w) \land \mathrm{T}(xv) \land \mathrm{T}(fvw).
\]
Therefore, the term $\underline{\mathsf{j}}$ interprets the join axiom in $\mathsf{VFU},$
where recall that $\mathrm{R}(x)$ and $x \in y$ are interpreted as $\mathrm{C}(x)$ and $\mathrm{T}(yx),$ respectively.
\end{lem}

\begin{pf}
Suppose that $\mathrm{C}(x)$ and $\forall y (\mathrm{T}(xy) \to \mathrm{C}(fy)).$
As the second conjunct is obvious from Lemma~\ref{lem:T-POS_in_VFU}, we show that $\mathrm{C}(\underline{\mathsf{j}}(x,f)).$
Take any $z$ and assume that $\neg \mathrm{T}(\underline{\mathsf{j}}(x,f)z),$ then $\mathrm{T}(\dot{\neg}(\underline{\mathsf{j}}(x,f)z))$ is derived with the help of Lemma~\ref{lem:compositionality_in_VFU}:

\begin {align}
        &\neg \mathrm{T}(\dot{\overbrace{\exists v, w. \ z = (v,w) \land \mathrm{T}  \mathit{(x v)}  \land \mathrm{T} \mathit{(fvw)}}}) \nonumber \\ 
        &\Longrightarrow \ 
        \neg \exists v, w. \ \mathrm{T}(\dot{\overbrace{z = (v,w) \land \mathrm{T}  \mathit{(x v)}  \land \mathrm{T} \mathit{(fvw)}}}) \nonumber \\
       &\Longleftrightarrow \neg \exists v, w. \ \mathrm{T}(z = (v,w)) \land \mathrm{T}  \mathit{(x v)}  \land \mathrm{T} \mathit{(fvw)} \nonumber \\
       &\Longleftrightarrow \forall v, w. \ \neg \mathrm{T}(z = (v,w)) \lor \neg \mathrm{T}  \mathit{(x v)} \lor \neg \mathrm{T} \mathit{(fvw)} \nonumber \\
       &\Longleftrightarrow \forall v, w. \ \mathrm{T}(\dot{\neg} (z \dot{=} (v,w))) \lor \mathrm{T} (\dot{\neg} \mathit{(x v)}) \lor \mathrm{T} (\dot{\neg}\mathit{(fvw))} \nonumber \\
       &\Longrightarrow \forall v, w. \ \mathrm{T}(\dot{\neg} (\dot{\mathrm{T}}(z \dot{=} (v,w)))) \lor \mathrm{T} (\dot{\neg} (\dot{\mathrm{T}}\mathit{(x v)})) \lor \mathrm{T} (\dot{\neg}(\dot{\mathrm{T}}(\mathit{fvw))}) \nonumber \\
       &\Longrightarrow \forall v, w. \ \mathrm{T}(\dot{\overbrace{\neg(\mathrm{T}(\mathit{z = (v,w)}) \land \mathrm{T} \mathit{(x v)} \land \mathrm{T} \mathit{(fvw)})}}) \nonumber \\
       &\Longleftrightarrow \mathrm{T}(\dot{\neg}(\underline{\mathsf{j}}(x,f)z)). \nonumber
\end {align} \qed
\end{pf}

To complete the interpretation of $\mathsf{T}_{0}$, we have to give the interpretation of inductive generation $\mathsf{i}$.
We define $\underline{\mathsf{i}}$ as the following term $\mathsf{acc}:$

$\mathsf{acc}(a,b) := \lambda z. \dot{\overbrace{TI[a,b,z]}},$ where 

%\begin{center}
$WF[a,b,f] := \forall x. \ \mathrm{T}(ax) \to [ \forall y. \ \mathrm{T}(ay) \to \mathrm{T}(b(y,x)) \to \mathrm{T}(\mathsf{l}_{a, b} (fy))] \to \mathrm{T}(\mathsf{l}_{a, b}(fx));$

$TI[a,b,z] := \mathrm{T}(az) \land \forall f. \ WF[a,b,f] \to \mathrm{T}(\mathsf{l}_{a, b}(fz)) .$

\

The term  $\mathsf{l}_{a, b}$ is defined below.
Informally speaking, $\mathsf{acc}(a,b)$ is the intersection of every set $f$ that includes the $<_b$-accessible part of $a$.
%In short, $\mathsf{acc}(a,b)$ is the accessible part of $a$ under the relation $b$.

\begin{lem}\label{lem:convergence_in_VFU}
For each class $u$ and $v$, we can take a universe $\mathsf{l}_{u, v}$ that reflects on both $u$ and $v$.
That is, 
\begin{center}
$\mathsf{VFU} \vdash \forall u, v. \ \mathrm{C}(u) \land \mathrm{C}(v) \to [\mathrm{U}(\mathsf{l}_{u, v}) \land u \sqsubset \mathsf{l}_{u, v} \land v \sqsubset \mathsf{l}_{u, v}].$
\end{center}
\end{lem}

\begin{pf}
Take any classes $u$ and $v$. We define a term $\oplus$ to be such that $u\oplus v = \lambda x. \dot{\overbrace{\mathrm{N}((\mathit{x})_{\mathrm{0}}) \to \mathrm{T}(\mathsf{d}_{\mathsf{N}}(\mathit{u(x)}_{\mathrm{1}})(v(x)_{\mathrm{1}})((x)_{\mathrm{0}})\mathsf{0})}}$
and let $\mathsf{l}_{u, v} := \mathsf{l}(u \oplus v)$. First, to show that $\mathsf{l}_{u, v}$ is a class, we prove that $u \oplus v$ is a class. Thus, taking any object $x$, we prove that $(u \oplus v)x$ is a proposition.
If $\neg \mathrm{N}((x)_0),$ we have $\mathrm{T}(\dot{\overbrace{\mathrm{N}((\mathit{x})_{\mathrm{0}}) \to \mathrm{T}(\mathsf{d}_{\mathsf{N}}(\mathit{u(x)}_{\mathrm{1}})(v(x)_{\mathrm{1}})((x)_{\mathrm{0}})\mathsf{0})}}),$ hence $\mathrm{P}((u \oplus v)x)$. 
Thus, we can assume $\mathrm{N}((x)_0)$.
If $(x)_{0} = \mathsf{0}$, then 
we have $\mathsf{d}_{\mathsf{N}}(u(x)_{1})(v(x)_{1})((x)_{0})\mathsf{0} = u(x)_{1}.$ As $u$ is a class, $(u \oplus v)x$ is a proposition, as required. Similarly, if $(x)_{0} \neq \mathsf{0},$ then $\mathsf{d}_{\mathsf{N}}(u(x)_{1})(v(x)_{1})((x)_{0})\mathsf{0} = v(x)_{1},$ and thus $(u \oplus v)x$ is a proposition.
Thus, in any case, $(u \oplus v)x$ is a proposition. Therefore, $u \oplus v$ is a class, and thus $(\mathsf{Lim})$ yields that $\mathsf{l}_{u, v}$ is a class.
Second, we show that $u \sqsubset \mathsf{l}_{u, v}.$ For an arbitrary object $x$, assume $\mathrm{T}(ux).$ Then, since $(u \oplus v)(\mathsf{0},x) = \dot{\overbrace{\mathrm{N}(\mathsf{0}) \to \mathrm{T}(\mathit{ux})}},$ it follows that $\mathrm{T}((u \oplus v)(\mathsf{0},x))$; thus, we obtain
$\mathrm{T}(\mathsf{l}_{u,v}((u \oplus v)(\mathsf{0},x))),$ which, by $(\mathsf{U} \mhyph \mathsf{Log})$, $(\mathsf{U} \mhyph \mathsf{Imp})$, and $(\mathsf{U} \mhyph \mathsf{Self})$, implies $\mathrm{T}(\mathsf{l}_{u, v}(ux))$. Similarly, $\mathrm{T}(\dot{\neg}(ux))$ implies $\mathrm{T}(\mathsf{l}_{u, v}(\dot{\neg}(ux)))$.
Thus, the conclusion $u \sqsubset \mathsf{l}_{u, v}$ is obtained. In the same way, we also have $v \sqsubset \mathsf{l}_{u, v}$. \qed
\end{pf}

\begin{lem}\label{lem:IG_in_VFU}
\begin{enumerate}
\item $\mathsf{VFU} \vdash \mathrm{C}(a) \land \mathrm{C}(b) \to \mathrm{C}(\mathsf{acc}(a,b));$
\item $\mathsf{VFU} \vdash \mathrm{C}(a) \land \mathrm{C}(b) \to \mathsf{Closed}' (a,b, \mathsf{acc}(a,b));$
\item $\mathsf{VFU} \vdash \mathrm{C}(a) \land \mathrm{C}(b) \land \mathsf{Closed}' (a,b, A(\bullet)) \to \forall x. \ \mathrm{T}(\mathsf{acc}(a,b)x) \to A(x)$, for each $\mathcal{L}_{FS}$-formula $A$.
\end{enumerate}
\end{lem}

\begin{pf}
We assume that $\mathrm{C}(a)$ and $\mathrm{C}(b)$.
\begin{enumerate}
\item Take any $z$; then we have to prove $\mathrm{P}(\mathsf{acc}(a,b)z)$. Therefore, supposing $\neg \mathrm{T}(\mathsf{acc}(a,b)z),$ we show $\mathrm{T}(\dot{\neg}(\mathsf{acc}(a,b)z)).$
For that purpose, we prove that $WF[a,b,f]$ is a proposition for any $f$.
First, by repeated use of Lemma~\ref{lem:compositionality_in_VFU}, we observe that $\dot{\overbrace{\forall y. \ \mathrm{T}(\mathit{ay}) \to \mathrm{T}(\mathit{b(y,x)}) \to \mathrm{T}(\mathsf{l}_{\mathit{a, b}} (\mathit{fy}))}}$ is a proposition:

\begin{align}
& \neg \mathrm{T}(\dot{\overbrace{\forall y. \ \mathrm{T}(\mathit{ay}) \to \mathrm{T}(\mathit{b(y,x)}) \to \mathrm{T}(\mathsf{l}_{\mathit{a, b}} (\mathit{fy}))}}) \nonumber \\
&\Longleftrightarrow \neg \forall y. \ \mathrm{T}(\dot{\overbrace{\mathrm{T}(\mathit{ay}) \to\mathrm{T}(\mathit{b(y,x)}) \to \mathrm{T}(\mathsf{l}_{\mathit{a, b}} (\mathit{fy}))}}) \nonumber \\
&\Longleftrightarrow \neg \forall y. \ \mathrm{T}(\dot{\overbrace{\mathrm{T}(\mathit{ay})}}) \to \mathrm{T}(\dot{\overbrace{\mathrm{T}(\mathit{b(y,x)}) \to \mathrm{T}(\mathsf{l}_{\mathit{a, b}} (\mathit{fy}))}}) \nonumber \\
&\Longleftrightarrow \neg \forall y. \ \mathrm{T}(\dot{\overbrace{\mathrm{T}(\mathit{ay})}}) \to \mathrm{T}(\dot{\overbrace{\mathrm{T}(\mathit{b(y,x)})}}) \to \mathrm{T}(\dot{\overbrace{\mathrm{T}(\mathsf{l}_{\mathit{a, b}} (\mathit{fy}))}})  \nonumber \\
&\Longleftrightarrow \exists y. \ \mathrm{T}(\dot{\overbrace{\mathrm{T}(\mathit{ay})}}) \land \mathrm{T}(\dot{\overbrace{\mathrm{T}(\mathit{b(y,x)})}}) \land \neg \mathrm{T}(\dot{\overbrace{\mathrm{T}(\mathsf{l}_{\mathit{a, b}} (\mathit{fy}))}}) \nonumber \\
&\Longleftrightarrow \exists y. \ \mathrm{T}(\dot{\overbrace{\mathrm{T}(\mathit{ay})}}) \land \mathrm{T}(\dot{\overbrace{\mathrm{T}(\mathit{b(y,x)})}}) \land \mathrm{T}(\dot{\overbrace{\neg \mathrm{T}(\mathsf{l}_{\mathit{a, b}} (\mathit{fy}))}}) \nonumber \\
&\Longleftrightarrow \exists y. \ \mathrm{T}(\dot{\overbrace{\neg (\mathrm{T}(\mathit{ay}) \to \mathrm{T}(\mathit{b(y,x)}) \to \mathrm{T}(\mathsf{l}_{\mathit{a, b}} (\mathit{fy})))}}) \nonumber \\
&\Longrightarrow \mathrm{T}(\dot{\overbrace{\neg \forall y. \ \mathrm{T}(\mathit{ay}) \to \mathrm{T}(\mathit{b(y,x)}) \to \mathrm{T}(\mathsf{l}_{\mathit{a, b}} (\mathit{fy}))}}). \nonumber
\end{align}
Therefore, it follows that $WF[a,b,f]$ is a proposition:
\begin{align}
& \neg \mathrm{T}(\dot{\overbrace{WF[a,b,f]}}) \nonumber \\
&\Longleftrightarrow \neg \forall x. \ \mathrm{T}(\dot{\overbrace{\mathrm{T}(\mathit{ax}) \to [\forall y. \ \mathrm{T}(\mathit{ay}) \to \mathrm{T}(\mathit{b(y,x)}) \to \mathrm{T}(\mathsf{l}_{\mathit{a, b}} (\mathit{fy}))] \to \mathrm{T}(\mathsf{l}_{\mathit{a, b}}(\mathit{fx}))}}) \nonumber \\
&\Longleftrightarrow \neg \forall x. \ \mathrm{T}(\dot{\overbrace{\mathrm{T}(\mathit{ax})}}) \to \mathrm{T}(\dot{\overbrace{[\forall y. \ \mathrm{T}(\mathit{ay}) \to \mathrm{T}(\mathit{b(y,x)}) \to \mathrm{T}(\mathsf{l}_{\mathit{a, b}} (\mathit{fy}))] \to \mathrm{T}(\mathsf{l}_{\mathit{a, b}}(\mathit{fx}))}}) \nonumber \\
&\Longleftrightarrow \exists x. \ \mathrm{T}(\dot{\overbrace{\mathrm{T}(\mathit{ax})}}) \land \neg \mathrm{T}(\dot{\overbrace{[\forall y. \ \mathrm{T}(\mathit{ay}) \to \mathrm{T}(\mathit{b(y,x)}) \to \mathrm{T}(\mathsf{l}_{\mathit{a, b}} (\mathit{fy}))] \to \mathrm{T}(\mathsf{l}_{\mathit{a, b}}(\mathit{fx}))}}) \nonumber \\
&\Longleftrightarrow \exists x. \ \mathrm{T}(\dot{\overbrace{\mathrm{T}(\mathit{ax})}}) \land \mathrm{T}(\dot{\overbrace{[\forall y. \ \mathrm{T}(\mathit{ay}) \to \mathrm{T}(\mathit{b(y,x)}) \to \mathrm{T}(\mathsf{l}_{\mathit{a, b}} (\mathit{fy}))] \land \neg \mathrm{T}(\mathsf{l}_{\mathit{a, b}}(\mathit{fx}))}}) \nonumber \\
&\Longleftrightarrow \exists x. \ \mathrm{T}(\dot{\overbrace{\mathrm{T}(\mathit{ax}) \land [\forall y. \ \mathrm{T}(\mathit{ay}) \to \mathrm{T}(\mathit{b(y,x)}) \to \mathrm{T}(\mathsf{l}_{\mathit{a, b}} (\mathit{fy}))] \land \neg \mathrm{T}(\mathsf{l}_{\mathit{a, b}}(\mathit{fx}))}}) \nonumber \\
&\Longrightarrow \mathrm{T}(\dot{\overbrace{\neg WF[a,b,f]}}). \nonumber
\end{align}
Using this, we can similarly get $\mathrm{T}(\dot{\neg}(\mathsf{acc}(a,b)z))$ from $\neg \mathrm{T}(\mathsf{acc}(a,b)z).$

\item
Assume that $\mathrm{T}(ax)$ and $\forall y [\mathrm{T}(ay) \to (\mathrm{T}(b(y,x)) \to \mathrm{T}(\mathsf{acc}(a,b)y))]$; then we want to derive $\mathrm{T}(\mathsf{acc}(a,b)x)$:

\begin{align}
& \mathrm{T}(\dot{\overbrace{\forall f. \ WF[a,b,f] \to \mathrm{T}(\mathsf{l}_{\mathit{a, b}}(\mathit{fx})) }}) \nonumber \\
&\Longleftrightarrow \forall f. \ \mathrm{T}(\dot{\overbrace{ WF[a,b,f] \to \mathrm{T}(\mathsf{l}_{\mathit{a, b}}(\mathit{fx})) }}) \nonumber \\
&\Longleftrightarrow \forall f. \ \mathrm{T}(\dot{\overbrace{ WF[a,b,f]}}) \to \mathrm{T}(\dot{\overbrace{\mathrm{T}(\mathit{l}_{\mathit{a, b}}(\mathit{fx})) }}). \nonumber
\end{align}
To show the last formula, we take any $f$ and suppose $\mathrm{T}(\dot{\overbrace{ WF[a,b,f]}})$. Then, we need to derive $\mathrm{T}(\dot{\overbrace{\mathrm{T}(\mathsf{l}_{\mathit{a, b}}(\mathit{fx})) }})$.
By the assumption, for any $y$ such that $\mathrm{T}(ay)$ and $\mathrm{T}(b(y,x))$ we have $\mathrm{T}(\mathsf{acc}(a,b)y).$
Thus, the supposition $\mathrm{T}(\dot{\overbrace{ WF[a,b,f]}})$ implies that $\mathrm{T}(\dot{\overbrace{\mathrm{T}(\mathsf{l}_{\mathit{a, b}}(\mathit{fy})) }})$.
As $y$ is arbitrary, this implies that 
$\mathrm{T}(\dot{\overbrace{\forall y. \ \mathrm{T}(\mathit{ay}) \to \mathrm{T}(\mathit{b(y,x)}) \to \mathrm{T}(\mathsf{l}_{\mathit{a, b}} (\mathit{fy}))}}).$
Combining this with the assumption $\mathrm{T}(ax),$ the desired conclusion  $\mathrm{T}(\dot{\overbrace{\mathrm{T}(\mathsf{l}_{\mathit{a, b}}(\mathit{fx})) }})$ follows from $\mathrm{T}(\dot{\overbrace{ WF[a,b,f]}})$.

\item Take any $x$ and assume that $\mathsf{Closed}' (a,b, A(\bullet))$ and $\mathrm{T}(\mathsf{acc}(a,b)x)$; we show $A(x)$.
Let $A'(x) := \mathsf{Closed}' (a,b,A(\bullet)) \to A(x)$, then we easily have $\mathsf{Closed}'(a,b,A'(\bullet))$ by logic.
Thus, in $\mathsf{VFU}$, we also obtain $\mathrm{T}(\mathsf{l}_{a, b} (\dot{\overbrace{\mathsf{Closed}'(\mathit{a,b,A'(\bullet)})}})).$ 
Next, we can derive $\mathsf{Closed}'(a,b,\mathrm{T}(\mathsf{l}_{a, b}(A'(\bullet))))$ in the following way:
\begin{align}
& \mathrm{T}(\mathsf{l}_{a, b} (\dot{\overbrace{\mathsf{Closed}'(\mathit{a,b,A'}(\bullet))}})) \nonumber \\
&\Longrightarrow \forall x. \ \mathrm{T}(\mathsf{l}_{a, b}(\dot{\overbrace{\mathrm{T}(\mathit{ax})}})) \to \mathrm{T}(\mathsf{l}_{a, b}(\dot{\overbrace{\forall y. \ \mathrm{T}(\mathit{ay}) \to \mathrm{T}(\mathit{b(y,x)}) \to A'(y)}})) \to \mathrm{T}(\mathsf{l}_{a, b}(\dot{\overbrace{A'(x)}})) \nonumber \\
&\Longleftrightarrow \forall x. \ \mathrm{T}(\mathit{ax}) \to \mathrm{T}(\mathsf{l}_{\mathit{a, b}}(\dot{\overbrace{\forall y. \ \mathrm{T}(\mathit{ay}) \to \mathrm{T}(\mathit{b(y,x)}) \to A'(y)}})) \to \mathrm{T}(\mathsf{l}_{a, b}(\dot{\overbrace{A'(x)}})) \nonumber \\
&\Longleftrightarrow \forall x. \ \mathrm{T}(ax) \to [\forall y. \ \mathrm{T}(\mathsf{l}_{\mathit{a, b}}(\dot{\overbrace{\mathrm{T}(\mathit{ay}) \to \mathrm{T}(\mathit{b(y,x)}) \to A'(y)}}))] \to \mathrm{T}(\mathsf{l}_{\mathit{a, b}}(\dot{\overbrace{A'(x)}})) \nonumber \\
&\Longleftrightarrow \forall x. \ \mathrm{T}(ax) \to [ \forall y. \ \mathrm{T}(ay) \to \mathrm{T}(b(y,x)) \to \mathrm{T}(\mathsf{l}_{a, b}(\dot{\overbrace{A'(y)}}))]
 \to \mathrm{T}(\mathsf{l}_{a, b}(\dot{\overbrace{A'(x)}})). \nonumber
\end{align}
Letting $f := \lambda x. \dot{\overbrace{A'(x)}}$, the assumption $\mathrm{T}(\mathsf{acc}(a,b)x)$ implies in $\mathsf{VFU}$ the formula
$\mathsf{Closed}'(a,b,\mathrm{T}(\mathsf{l}_{a, b}(A'(\bullet)))) \to \mathrm{T}(\mathsf{l}_{a, b} (\dot{\overbrace{A'(x)}}))$.
Therefore, we obtain $\mathrm{T}(\mathsf{l}_{a, b} (\dot{\overbrace{A'(x)}})),$ which yields $A'(x)$ by $(\mathsf{U}\mhyph \mathsf{Out})$.
Finally, combining this with the assumption  $\mathsf{Closed}' (a,b, A(\bullet))$, the conclusion $A(x)$ follows. \qed
\end{enumerate}
\end{pf}

\begin{thm}\label{thm:lower-bound_of_VFU}
For each $\mathcal{L}_{FS}$-sentence $A$,
if $\mathsf{T}_{0} \vdash A$, then $\mathsf{VFU} \vdash A'$.
In particular, every $\mathcal{L}$-theorem of $\mathsf{T}_{0}$ is derivable in $\mathsf{VFU}$. 
\end{thm}

\subsection{Truth-as-provability interpretation of $\mathsf{VFU}$}\label{subsec:model_of_VFU}
In this subsection, we give a model of $\mathsf{VFU}$ by generalising Cantini's \emph{truth}-\emph{as}-\emph{provability} interpretation for $\mathsf{VF}$ \cite{cantini1990theory, cantini1996logical},
which is formalisable in a suitable set theory, and thus the upper-bound of $\mathsf{VFU}$ is obtained (see subsection~\ref{subsec:upper-bound_of_VFU}).
The idea of our truth-as-provability interpretation is that the truth predicate $\mathrm{T}(x)$ is intepreted as the derivability of $x$ in the indexed infinitary sequent calculus, as is displayed on the table below. Then, each axiom of $\mathsf{VFU}$ is shown to be true under this interpretation.
Here, a \emph{sequent} $\varGamma$ is a finite set of closed terms, each of which is $\beta$-equivalent to $\dot{A}$ for some sentence $A$. 
To present the system in the form of Tait calculus, we consider only negation normal sentences. Therefore, the negation symbol $\neg$ may come only in front of atomic sentences, and then the global negation $\neg A$ becomes a defined expression with the help of De Morgan's law.
Note that the conditional $A \to B$ is defined by $\neg A \lor B$.
For simplicity, we do not distinguish terms which have the same reduct, thus we can suppose that every term is of the form $\dot{A}$ for some negation normal sentence $A$. 
%For example, if closed terms $s$ and $\dot{A}$ have the same reduct, then the sequent $\{ s \}$ is assumed to be identical with $\{ \dot{A} \}$.
For readability, we often simply write $A$ instead of $\dot{A}$.
Moreover, since we also consider terms of the form $\mathsf{l}ab,$ it is useful to treat them as if they were sentences.
Thus, we introduce a new binary predicate symbol $\mathrm{L}_x(y)$ and we let $\dot{\overbrace{\mathrm{L}_a(b)}} := \mathsf{l}ab.$ Similarly, let $\dot{\overbrace{\neg \mathrm{L}_a(b)}} := \dot{\neg}(\mathsf{l}ab).$

Next, we explain the calculus in more detail.
In the calculus, the predicate $\sststile{}{\alpha,\beta,\gamma} \varGamma$ means that $\varGamma$ is derived in the system with the $\mathrm{U}$-rank $\alpha,$
the $\mathrm{T}$-rank, and the derivation length $\gamma.$
We introduce several notations: 
\begin{itemize}
\item $\sststile{}{\alpha} \varGamma$ means $\sststile{}{\alpha, \beta, \gamma} \varGamma$ for some $\alpha, \beta, \gamma$;
\item $\sststile{}{<\alpha} \varGamma$ means  $\sststile{}{\alpha_0} \varGamma$ for some $\alpha_0<\alpha;$
\item $\sststile{}{\alpha,\beta} \varGamma$ means $\sststile{}{\alpha, \beta, \gamma} \varGamma$ for some $\gamma$;
\item $\sststile{}{\alpha,\beta,<\gamma} \varGamma$ means
$\sststile{}{\alpha,\beta,\gamma_0} \varGamma$ for some $\gamma_0 < \gamma$;
\item $\sststile{}{\alpha,<\beta,<\gamma} \varGamma$ means
$\sststile{}{\alpha,\beta_0,\gamma_0} \varGamma$ for some $\beta_0 < \beta$ and $\gamma_0 < \gamma$;
\item $\sststile{}{\leq\alpha,\leq\beta,\leq\gamma} \varGamma$ means
$\sststile{}{\alpha_0,\beta_0,\gamma_0} \varGamma$ for some $\alpha_0 \leq \alpha,$ $\beta_0 \leq \beta$ and $\gamma_0 \leq \gamma$;
\item $\sststile{}{} \varGamma$ means  $\sststile{}{\alpha} \varGamma$ for some $\alpha;$
\item if $\neg$ comes to the left of one of the above expressions, it negates the whole expression.
For example, $\neg\sststile{}{} \varGamma$ means that it is not the case that $\sststile{}{} \varGamma$.
\end{itemize}

Under these conventions, each rule is explained as follows. 
The rule $(\text{Lit})$ says that an $\mathcal{L}$-literal $P$ that is true in the closed term model $\mathcal{CTT}$ is derivable. The rules $(\text{Log}), (\land), (\lor)_i, (\exists), (\forall)$ are given similarly to the standard sequent calculus.
In particular, $(\forall)$ has infinitely many premises for each closed term $a$.
The rules $(\mathrm{T})$ and $(\neg\mathrm{T})$ respectively introduce $\mathrm{T}$ and $\neg\mathrm{T},$ with an increase in the $\mathrm{T}$-rank. 
%To avoid inconsistency, the premise does not have a context.
Note that the context of the premise of $(\mathrm{T})$ and $(\neg\mathrm{T})$ must be empty; otherwise, the system would be inconsistent according to the liar paradox.
The rule $(\text{Weak})$ assures the monotonicity of the derivability with respect to the $\mathrm{U}$-rank. 
Similarly to the operator $\Phi$ in $\mathcal{M}_{\mathsf{KF}},$ the rules $(\mathrm{U})$ and $(\neg \mathrm{U})$ have the side condition $\dagger$, which consists of the following conditions. 
First, $\alpha$ needs to be a successor ordinal (cf. the operator $\Phi$ in $\mathcal{M}_{\mathsf{KF}}$). Thanks to this condition, we can assure that $\sststile{}{<\alpha}$ is closed under $(\text{Lit}),$ $(\text{Log}),$ $(\land),$ $(\lor)_i,$ $(\exists)$ and $(\forall).$
In particular, if $\sststile{}{<\alpha}\varGamma,A(b)$ for all $b$, then $\sststile{}{<\alpha}\varGamma,\forall x. \ A(x)$ holds.
Second, $a$ must satisfy the following:
\[
\sststile{}{<\alpha} ac \ \text{or} \ \sststile{}{<\alpha} \dot{\neg}(ac) \ \text{holds for all closed terms} \ c.
\]
Thus, it roughly says that $a$ is a class, provably in $\sststile{}{<\alpha}.$
We express this property as $\sststile{}{<\alpha}a:\text{Class}$.
The third condition is that neither $\mathrm{L}_a(b)$ nor $\neg \mathrm{L}_a(b)$ are derived in $\sststile{}{<\alpha}:$
\[
\neg \sststile{}{<\alpha}\mathrm{L}_a(b), \ \text{and} \ \neg \sststile{}{<\alpha}\neg \mathrm{L}_a(b).
\]
Thus, it roughly says that $\mathrm{L}_a(b)$ is not a proposition in $\sststile{}{<\alpha}.$
We express this as $\neg \sststile{}{<\alpha}\mathrm{L}_a(b):\text{Prop}.$

\begin{table}[H] \label{tb:sequent}
\caption{Sequent system $\sststile{}{\alpha,\beta,\gamma}$}
\centering
\begin{tabular}{ll}
\hline

\\

\infer[(\text{Lit})]{\sststile{}{\alpha,\beta,\gamma} \varGamma,P}{\mathcal{CTT} \models P} & \infer[(\mathrm{Log})]{\sststile{}{\alpha,\beta,\gamma} \varGamma,\mathrm{T}(a), \neg \mathrm{T}(a)}{} \\

\\

\infer[(\land)]{\sststile{}{\alpha,\beta,\gamma} \varGamma,A_0 \land A_1}{\sststile{}{\alpha,\beta,<\gamma} \varGamma,A_0 \land A_1,A_0 \ \ \sststile{}{\alpha,\beta,<\gamma} \varGamma,A_0 \land A_1,A_1} & \infer[(\lor)_{i}]{\sststile{}{\alpha,\beta,\gamma} \varGamma,A_{0} \lor A_{1}}{\sststile{}{\alpha,\beta,<\gamma} \varGamma,A_{0} \lor A_{1},A_{i} \ (i \leq 1)} \\

\\

\infer[(\exists)]{\sststile{}{\alpha,\beta,\gamma} \varGamma,\exists x. \ A(x)}{\sststile{}{\alpha,\beta,<\gamma} \varGamma,\exists x. \ A(x),A(a)} & \infer[(\forall)]{\sststile{}{\alpha,\beta,\gamma} \varGamma,\forall x. \ A(x)}{\sststile{}{\alpha,\beta,<\gamma} \varGamma,\forall x. \ A(x),A(a), \ \mathrm{for \ all} \ a}  \\

\\

\infer[(\mathrm{T})]{\sststile{}{\alpha,\beta,\gamma} \varGamma,\mathrm{T}(\dot{A})}{
\sststile{}{\alpha,<\beta,<\gamma} A} & \infer[(\neg \mathrm{T})]{\sststile{}{\alpha,\beta,\gamma} \varGamma,\neg \mathrm{T}(\dot{A})}{\sststile{}{\alpha,<\beta,<\gamma} \neg A} \\

\\

\infer[(\text{Weak})]{\sststile{}{\alpha,\beta,\gamma} \varGamma}{\sststile{}{<\alpha} \varGamma} \\

\\

\infer[(\mathrm{U})^{\dagger}]{\sststile{}{\alpha,\beta,\gamma} \varGamma,\mathrm{L}_{a}(b)}{\sststile{}{<\alpha}b} &
\infer[(\neg\mathrm{U})^{\dagger}]{\sststile{}{\alpha,\beta,\gamma} \varGamma,\neg\mathrm{L}_{a}(b)}{\neg \sststile{}{<\alpha}b} \\

\\

\hline 
 
\end{tabular}
\end{table}

By transfinite induction, we easily obtain the following:
\begin{lem}\label{lem:consistency_weakening}
\begin{description}
\item[(Consistency)] The empty sequent is not derivable: $\neg \sststile{}{}\emptyset.$
\item[(Weakening)] If $\sststile{}{\leq\alpha,\leq\beta,\leq\gamma}\varGamma,$ then $\sststile{}{\alpha,\beta,\gamma}\varGamma, \varDelta.$
\end{description}
\end{lem}

We now show the cut-admissibility of the calculus.
For a formula $A$,
the logical complexity $\mathrm{co}(A)$ is defined as usual: if $P$ is any literal of $\mathcal{L} \cup \{ \mathrm{T}(x),\mathrm{L}_x(y) \}$, then $\mathrm{co}(P):=0;$
$\mathrm{co}(A \land B):=\mathrm{co}(A \lor B) := \max(\mathrm{co}(A),\mathrm{co}(B)) + 1;$ $\mathrm{co}(\forall x. A(x)) := \mathrm{co}(\exists x. A(x)) := \mathrm{co}(A(x))+1.$

\begin{lem}[Cut-admissibility]\label{lem:cut-admissibility}
If $\sststile{}{\alpha,\beta,\gamma}\varGamma,A$ and $\sststile{}{\delta,\varepsilon,\zeta}\varDelta,\neg A,$
then $\sststile{}{\max(\alpha,\delta)}\varGamma,\varDelta$.
\end{lem}

\begin{pf}\label{pf:lem:cut-admissibility}
We show the claim by septuple induction on $\alpha$, $\delta$, $\beta$, $\varepsilon$, $\mathrm{co}(A)$, $\gamma$, and $\zeta$.
The case where either $\varGamma, A$ or $\varDelta, \neg A$ is obtained by $(\text{Weak})$ is clear by the induction hypothesis. Thus, we can rule out such a case. 
If $A$ or $\neg A$ is not principal in the last rule, then the conclusion follows by the induction hypothesis. 
For example, assume that $\varDelta, \neg A$ is derived by $(\forall)$ from the premises $\sststile{}{\delta, \varepsilon, \zeta_{a}} \varDelta_{a}, \neg A$ with $\zeta_{a} < \zeta$ for all closed terms $a,$
then the induction hypothesis yields $\sststile{}{\max(\alpha,\delta)} \varGamma, \varDelta_a.$
Then, $(\forall)$ derives $\sststile{}{\max(\alpha,\delta)}\varGamma, \varDelta,$ as required.

Finally, we consider the case where both $A$ and $\neg A$ are principal. 
The inductive case $\text{co}(A) > 0$ is proved by a standard cut-elimination argument (cf. \cite[Theorem~62.1]{cantini1996logical}). 
Thus, we confine ourselves to the base cases $A \equiv \mathrm{T}(\dot{B})$ and $A \equiv \mathrm{L}_{a}(\dot{B})$.
\begin{description}
\item[$A \equiv \mathrm{T}(\dot{B})$] Firstly, if $\varGamma, A$ is an instance of $(\mathrm{Log})$, then $\neg \mathrm{T}(\dot{B})$ is contained in $\varGamma,$ hence
we have $\neg \mathrm{T}(\dot{B}), \varDelta \subseteq \varGamma, \varDelta$. Therefore, Lemma~\ref{lem:consistency_weakening} implies the conclusion. 
The case where $\varDelta, \neg A$ is $(\mathrm{Log})$ is similar.
Secondly, we assume that $\varGamma, A$ and $\varDelta, \neg A$ are respectively obtained by $(\mathrm{T})$ and $(\neg \mathrm{T})$:
\begin{center} \ 
\infer[(\mathrm{T})]{\sststile{}{\alpha, \beta, \gamma} \varGamma, \mathrm{T}(\dot{B})}{\sststile{}{\alpha,\beta', \gamma'} B}
\ \ \ \ \ \ \ \ \ \ \infer[(\neg\mathrm{T})]{\sststile{}{\delta,\varepsilon, \zeta} \varDelta, \neg\mathrm{T}(\dot{B})}{\sststile{}{\delta,\varepsilon', \zeta'} \neg B},
 \end{center}
where $\beta' < \beta, \gamma' < \gamma, \varepsilon' < \varepsilon$ and $\zeta' < \zeta.$
Since $\beta' < \beta$, the induction hypothesis for the premises yields that $\sststile{}{\max(\alpha,\delta)} \emptyset,$ which contradicts Lemma~\ref{lem:consistency_weakening}. Thus, this case cannot occur.

\item[$A \equiv \mathrm{L}_{a}(b)$]  As a crucial case, we suppose that $\varGamma, \mathrm{L}_{a}(b)$ and $\varDelta, \neg \mathrm{L}_{a}(b)$ are obtained by $(\mathrm{U})$ and $(\neg \mathrm{U})$, respectively:
\begin{center} \ 
\infer[(\mathrm{U})^{\dagger}]{\sststile{}{\alpha, \beta, \gamma} \varGamma, \mathrm{L}_{a}(b)}{\sststile{}{<\alpha} b}
\ \ \ \ \ \ \ \ \ \ \infer[(\neg\mathrm{U})^{\dagger}]{\sststile{}{\delta,\varepsilon,\zeta} \varDelta, \neg\mathrm{L}_{a}(b)}{\neg\sststile{}{<\delta} b}.
\end{center}
Here, we have the following side conditions:
\begin{gather}
%\sststile{}{<\alpha} a:\mathrm{Class} \label{pf_1:lem:cut-admissibility} \\
\neg \sststile{}{<\alpha} \mathrm{L}_{a}(b):\mathrm{Prop} \label{pf_1:lem:cut-admissibility} \\
%\sststile{}{<\delta} a:\mathrm{Class} \label{pf_3:lem:cut-admissibility} \\
\neg \sststile{}{<\delta} \mathrm{L}_{a}(b):\mathrm{Prop} \label{pf_2:lem:cut-admissibility}
\end{gather}

By (\ref{pf_1:lem:cut-admissibility}) and (\ref{pf_2:lem:cut-admissibility}), we clearly have $\alpha=\delta.$
Therefore, the premise $\sststile{}{<\alpha}b$ is identical with $\sststile{}{<\delta}b,$ which contradicts the other premise $\neg\sststile{}{<\delta}b.$
Thus, this case cannot occur. \qed
\end{description}
\end{pf}

Using the cut-admissibility, we can give a model of $\mathsf{VFU}.$
The $\mathcal{L}_{FS}$-model $\mathcal{M}_{\mathsf{VFU}}$ is an expansion of $\mathcal{CTT},$ in which the vocabularies of $\mathcal{L}$ and the additional constant symbols of $\mathcal{L}_{FS}$ are interpreted in the same way as for $\mathcal{M}_{\mathsf{KF}}$ in Section~\ref{subsec:Kahle's_model_for_KFU}.
Then, $\mathrm{T}(x)$ is interpreted as $\sststile{}{}\{x'\}$,
where $\{x'\}$ is the singleton of the negation normal form $x'$ of $x$.
%\footnote{For the definition of negation normal form, see, e.g.~\cite{cantini1990theory}. Indeed, we can primitive recursively translate a given formula into its negation normal form.} 
For simplicity, we write $\sststile{}{} x$ instead of $\sststile{}{}\{x'\}$.
Similarly, $\mathrm{U}(x)$ is interpreted as the statement that for some closed term $a,$ $x$ is of the form $\mathsf{l}a$ and $\sststile{}{}a:\text{Class}$ holds.

The next lemma verifies the $\mathsf{VF}$-axioms in $\mathrm{U}$ of Definition~\ref{def:VFU}.

\begin{lem}\label{lem_1:model_of_VFU}
Let $a$ and $b$ be any closed terms and suppose $\sststile{}{}a:\mathrm{Class}.$ 
Then, the following hold in $\mathcal{M}_{\mathsf{VFU}}$: 
\begin{description}
\item[$(\mathsf{U}\mhyph \mathsf{Elem})$] If $\mathcal{CTT} \models P$, then  $\sststile{}{} \mathrm{L}_{a}(\dot{P}),$ for each $\mathcal{L}$-literal sentence $P.$

\item[$(\mathsf{U}\mhyph \mathsf{Imp})$] If $\sststile{}{} \mathrm{L}_{a}(\dot{\overbrace{A \to B}})$ and  $\sststile{}{}\mathrm{L}_{a}(\dot{A})$, then $\sststile{}{}\mathrm{L}_{a}(\dot{B}).$

\item[$(\mathsf{U}\mhyph \mathsf{Univ})$] If $\sststile{}{}\mathrm{L}_{a}(\dot{\overbrace{A(c)}})$ for every closed term $c$, then $\sststile{}{}\mathrm{L}_{a}(\dot{\overbrace{\forall x. \ A(x)}}).$

\item[$(\mathsf{U}\mhyph \mathsf{Log})$] $\sststile{}{}\mathrm{L}_{a}(\dot{A})$ for each logical theorem $A$.

\item[$(\mathsf{U}\mhyph \mathsf{Cons})$] It is not the case that both $\sststile{}{}\mathrm{L}_{a}(b)$ and $\sststile{}{}\mathrm{L}_{a}(\dot{\neg}b).$

\item[$(\mathsf{U}\mhyph \mathsf{Self})$] $\sststile{}{}\mathrm{L}_{a}(b)$ if and only if $\sststile{}{}\mathrm{L}_{a}(\dot{\overbrace{\mathrm{T}(\mathit{b})}}).$ Similarly, $\sststile{}{}\mathrm{L}_{a}(\dot{\neg}b)$ if and only if $\sststile{}{}\mathrm{L}_{a}(\dot{\overbrace{\neg \mathrm{T}(\mathit{b})}}).$
\end{description}
\end{lem}

\begin{pf}
\begin{description}
\item[$(\mathsf{U}\mhyph \mathsf{Elem})$] Assuming $\mathcal{M}_{\mathsf{VFU}} \models P,$ we show that $\sststile{}{} \mathrm{L}_{a}(\dot{P}).$ From the supposition, we can take the least successor ordinal $\alpha$ such that $\sststile{}{<\alpha}a:\mathrm{Class}.$ Then, we obviously have $\neg\sststile{}{<\alpha}\mathrm{L}_a(\dot{P}):\text{Prop}.$
In addition, we have $\sststile{}{0}P$ by $(\mathrm{Lit})$; thus, by $(\mathrm{U})$, we can deduce $\sststile{}{\alpha} \mathrm{L}_{a}(\dot{P}),$ as required.

\item[$(\mathsf{U}\mhyph \mathsf{Imp})$] Assuming $\sststile{}{} \mathrm{L}_{a}(\dot{\overbrace{A \to B}})$ and  $\sststile{}{}\mathrm{L}_{a}(\dot{A}),$ we have to prove $\sststile{}{}\mathrm{L}_{a}(\dot{B}).$ Similarly to the above, we take the least successor ordinal $\alpha$ such that $\sststile{}{<\alpha}a:\mathrm{Class}.$ Then, we clearly have $\neg\sststile{}{<\alpha}\mathrm{L}_a(\dot{A}):\text{Prop}$ and $\sststile{}{\alpha}\mathrm{L}_a(\dot{A}):\text{Prop}.$
Here, if $\sststile{}{\alpha}\neg\mathrm{L}_a(\dot{A})$, then the assumption $\sststile{}{}\mathrm{L}_{a}(\dot{A})$ implies $\sststile{}{} \emptyset$ by Lemma~\ref{lem:cut-admissibility}, which contradicts Lemma~\ref{lem:consistency_weakening}.
Thus, $\sststile{}{\alpha}\mathrm{L}_a(\dot{A}),$ which yields $\sststile{}{<\alpha}A.$ 
Similarly, we also have $\sststile{}{<\alpha} A \to B.$ Therefore, again by Lemma~\ref{lem:cut-admissibility}, it follows that 
$\sststile{}{<\alpha}B,$ and thus we obtain $\sststile{}{\alpha}\mathrm{L}_{a}(\dot{B})$ by $(\mathrm{U})$.

\item[$(\mathsf{U}\mhyph \mathsf{Univ})$] Assuming $\sststile{}{}\mathrm{L}_{a}(\dot{\overbrace{A(c)}})$ for any closed term $c$, we show $\sststile{}{}\mathrm{L}_{a}(\dot{\overbrace{\forall x. \ A(x)}}).$ Similarly to the above, we take the least successor ordinal $\alpha$ such that $\sststile{}{<\alpha}a:\mathrm{Class}.$ Then, for all $c$, we clearly have $\sststile{}{\alpha}\mathrm{L}_a(\dot{\overbrace{A(c)}})$, and thus 
$\sststile{}{<\alpha}A(c)$. Since $\alpha$ is a successor ordinal, we also have $\sststile{}{\alpha-1}A(c)$ for all $c$, which implies $\sststile{}{\alpha-1}\forall x. \ A(x)$
by the rule $(\forall)$. Thus, it follows that $\sststile{}{}\mathrm{L}_{a}(\dot{\overbrace{\forall x. \ A(x)}}).$
\end{description}

The other cases are similarly proved. \qed
\end{pf}

Similarly to Lemma~\ref{lem_1:model_of_VFU}, the structural properties of $\mathsf{VFU}$ are satisfied in $\mathcal{M}_{\mathsf{VFU}}:$
\begin{lem}\label{lem_2:model_of_VFU}
Let $a$ and $b$ be any closed terms and suppose $\sststile{}{}a:\mathrm{Class}.$ 
Then, the following hold in $\mathcal{M}_{\mathsf{VFU}}$: 
\begin{description}
\item[$(\mathsf{U}\mhyph \mathsf{Class})$] $\sststile{}{} \mathrm{L}_{a}(b)$ or $\sststile{}{}\neg \mathrm{L}_{a}(b).$

\item[$(\mathsf{U}\mhyph \mathsf{True})$] If $\sststile{}{} \mathrm{L}_{a}(\dot{A}),$ then $\sststile{}{}A.$

\item[$(\mathsf{Lim})$] If $\sststile{}{} ab,$ then $\sststile{}{}\mathrm{L}_{a}(ab).$ 
If $\sststile{}{}\dot{\neg}(ab),$ then $\sststile{}{}\mathrm{L}_{a}(\dot{\neg}(ab)).$
\end{description}
\end{lem}

Finally, we verify the axiom $(\mathsf{U}\mhyph \mathsf{Out}):$

\begin{lem}\label{lem_3:model_of_VFU}
Let $a$ be any closed term and suppose $\sststile{}{}a:\mathrm{Class}.$
Then, the following is satisfied in $\mathcal{M}_{\mathsf{VFU}}:$
\begin{description}
\item[$(\mathsf{U} \mhyph \mathsf{Out})$] $\sststile{}{}\mathrm{L}_a(\dot{A})$ implies $\mathcal{M}_{\mathsf{VFU}} \models A$ for any $\mathcal{L}_{FS}$-sentence $A$.
\end{description}
\end{lem}

\begin{pf}
Since we observed that $(\mathsf{U} \mhyph \mathsf{True})$ is satisfied in $\mathcal{M}_{\mathsf{VFU}}$, it suffices to show that $\sststile{}{}A$ implies $\mathcal{M}_{\mathsf{VFU}} \models A$. 
Then, in exactly the same way as for \cite[Theorem~63.4]{cantini1996logical},
we can prove, by transfinite induction, that if $\sststile{}{}\varGamma$ for a sequent $\varGamma$ which consists only of $\mathcal{L}_{FS}$-sentences, then at least one sentence of $\varGamma$ is true in $\mathcal{M}_{\mathsf{VFU}}.$  \qed
\end{pf}

In conclusion, every theorem of $\mathsf{VFU}$ is satisfied in $\mathcal{M}_{\mathsf{VFU}}.$
\begin{thm}\label{thm:model_of_VFU}
$\mathcal{M}_{\mathsf{VFU}} \models \mathsf{VFU}.$
\end{thm}

\begin{rmk}\label{rmk:structure_of_M_VFU}
Similar to Remark~\ref{rmk:structure_of_M_KF}, we can also verify the axioms $(\mathsf{U} \mhyph \mathsf{Tran})$, $(\mathsf{U} \mhyph \mathsf{Dir})$,  $(\mathsf{U} \mhyph \mathsf{Nor})$ and $(\mathsf{U} \mhyph \mathsf{Lin})$ in $\mathcal{M}_{\mathsf{VFU}}$.
\end{rmk}

\subsection{Upper bound of $\mathsf{VFU}$}\label{subsec:upper-bound_of_VFU}
To determine the upper bound of $\mathsf{VFU}$, we want to formalise the model $\mathcal{M}_{\mathsf{VFU}}$ of the previous subsection. However, the theory $\mathsf{FID}([\mathsf{POS,QF}])$ of Section~\ref{subsec:upper-bound_of_KFUPI} is not expressive enough to formalise the cut-admissibility argument of Lemma~\ref{lem:cut-admissibility}.
Thus, we will construct $\mathcal{M}_{\mathsf{VFU}}$ within the Kripke-Platek set theory $\mathsf{KPi},$ which is proof-theoretically equivalent to $\mathsf{FID}([\mathsf{POS,QF}]),$ and hence it follows that $\mathsf{VFU} \leq \mathsf{T}_0.$

For the formulation of $\mathsf{KPi},$ we follow \cite{jager2001first}.
For the language $\mathcal{L}'$ of first-order Peano arithmetic, let $\mathcal{L}^{*} := \mathcal{L}' \cup \{ \in, \mathsf{N}, \mathrm{S}, \mathrm{Ad} \}$,
where $\in$ is the membership relation symbol; $\mathsf{N}$ is the set constant for the natural numbers; $\mathrm{S}$ is the unary predicate symbol, expressing that a given object is a set; and the unary predicate symbol $\mathsf{Ad}$ says that an object is an admissible set. 
Moreover, we assume that $\mathcal{L}^{*}$ contains restricted quantifiers $\forall x \in y$ and $\exists x \in y$ as primitive symbols.
In $\mathcal{L}^{*},$ the equality symbol $=$ is defined as the following formula:
$(a = b) := [a \in \mathsf{N} \land b \in \mathsf{N} \land a =_{\mathsf{N}} b] \lor [ \mathrm{S}(a) \land \mathrm{S}(b) \land ( \forall x \in a. \ x \in b) \land (\forall x \in b. \ x \in a)],$
where $=_{\mathsf{N}}$ is a primitive recursive equality on natural numbers, and thus is contained in $\mathcal{L}'.$
An $\mathcal{L}^{*}$-formula $A$ is $\Delta_{0}$ if $A$ contains no unrestricted quantifiers.
Let $\mathrm{Tran}(x)$ be a defined $\Delta_{0}$-predicate that expresses that $x$ is a transitive set.
For the $\mathcal{L}^{*}$-formula $A$, let $A^{a}$ be the result of replacing each unrestricted quantifier $\exists x. ()$ and $\forall x. ()$ in $A$ by $\exists x \in a. ()$ and $\forall x \in a. ()$, respectively.

\begin{defn}[cf. \cite{jager2001first}]\label{def:KPi}
The $\mathcal{L}^{*}$-theory $\mathsf{KPi}$ consists of the following axioms:
\begin{description}
\item[$\mathsf{N}$-Induction and foundation.] For all $\mathcal{L}^{*}$-formulas $A$,
\begin{itemize}
\item $A(0) \land [\forall x \in \mathsf{N}. \ A(x) \to A(x + 1)] \to \forall x \in \mathsf{N}. \ A(x)$; 
\item $[ \forall x. \ (\forall y \in x. \ A(y)) \to A(x) ] \to \forall x. \ A(x)$.
\end{itemize}
\item[Ontological axioms.] For all terms $a,b$ and $\vec{c}$ of $\mathcal{L}^{*}$, all function symbols $h$ and relation symbols $R$ of $\mathcal{L}'$ and all axioms $A(\vec{x})$ of Set-theoretic axioms whose free variables belong to $\vec{x}$,
\begin{itemize}
\item $a \in \mathsf{N} \leftrightarrow \neg \mathrm{S}(a);$
\item $\vec{c} \in \mathsf{N} \to h(\vec{c}) \in \mathsf{N};$
\item $R(\vec{c}) \to \vec{c} \in \mathsf{N};$
\item $a \in b \to \mathrm{S}(b);$
\item $\mathrm{Ad}(a) \to [\mathsf{N} \in a \land \mathrm{Tran}(a)];$
\item $\mathrm{Ad}(a) \to \forall \vec{x} \in a. \ A^{a}(\vec{x})$.
\end{itemize}
\item[Number-theoretic axioms.] For all axioms $A(\vec{x})$ of Peano arithmetic whose free variables belong to $\vec{x}$,
\begin{itemize}
\item $\forall \vec{x} \in \mathsf{N}. \ A^{\mathsf{N}}(\vec{x}).$
\end{itemize}
\item[Set-theoretic axioms.] For all terms $a$ and $b$ and all $\Delta_{0}$-formulas $A(x)$ and $B(x,y)$ of $\mathcal{L}^{*},$
\begin{description}
\item[Pair.] $\exists x. \ a \in x \land b \in x;$
\item[Transitive Hull.] $\exists x. \ a \subset x \land \mathrm{Tran}(x);$
\item[$\Delta_{0}$-Separation.] $\exists y. \ \mathrm{S}(y) \land y = \{ x \in a :  A(x) \};$
\item[$\Delta_{0}$-Collection.] $[\forall x \in a. \exists y. B(x,y)] \to \exists z. \ \forall x \in a. \ \exists y \in z. \ B(x,y);$
\end{description} 
\item[Limit axiom.] $\forall x. \ \exists y. \ x \in y \land \mathrm{Ad}(y).$
\end{description}
\end{defn}

Now, we describe how $\mathcal{M}_{\mathsf{VFU}}$ is formalised in $\mathsf{KPi}.$
Since $\mathsf{KPi}$ contains Peano arithmetic, the interpretation of $\mathcal{L}$ and the constant symbols of $\mathcal{L}_{FS}$ can be given by using fixed G\"{o}del numbering.
Thus, the remaining task is to formalise the sequent calculus $\sststile{}{}$ in the previous subsection for the interpretation of $\mathrm{T}$ and $\mathrm{U}.$
For that purpose, $\sststile{}{}$ in which $\mathrm{T}$-ranks and the derivation lengths are omitted is firstly formalised via the  operator $\Phi(X,\ulcorner \varGamma \urcorner)$ defined below, where $\varGamma$ is (the code of) a sequent  in the sense of the previous subsection. Since expressions of $\mathcal{L}_{FS} \cup \{ \mathrm{L}_x(y) \}$ are coded as natural numbers, $\Phi$ can be given as a formula of $\mathcal{L}^{*}$:

The $\Delta_0$-formula $\Phi_0(X,\varGamma)$ is defined to be the disjunction of the following:
\begin{itemize}
\item $P \in \varGamma$, for an $\mathcal{L}$-literal $P$ true in $\mathcal{CTT}$,
\item $A \land B \in \varGamma$, for some $A,B$ such that $\varGamma,A \in X$ and $\varGamma,B \in X$;
\item $A \lor B \in \varGamma$, for some $A,B$ such that either $\varGamma,A \in X$ or $\varGamma,B \in X$;

\item $\forall x. \ A(x) \in \varGamma$, for some $\forall x. \ A(x)$ such that $\varGamma,A(a) \in X$ for all closed terms $a$;
\item $\exists x. \ A(x) \in \varGamma$, for some $\exists x. \ A(x)$ such that $\varGamma,A(a) \in X$ for some closed term $a$;
\end{itemize}

Similarly, the $\Delta_0$-formula $\Phi_1(X,\varGamma)$ is defined to be the disjunction of the following:
\begin{itemize}
\item $\mathrm{T}(\dot{A}) \in \varGamma,$ for some $A \in X;$
\item $\neg \mathrm{T}(\dot{A}) \in \varGamma,$ for some $\neg  A \in X;$
\end{itemize}

Finally, the $\Delta_0$-formula $\Phi_2(X,\varGamma)$ is defined to be the disjunction of the following:
\begin{itemize}
\item $\mathrm{L}_a(b) \in \varGamma,$ for some $a,b$ such that $b \in X$ and $(\star)$ holds;
\item $\neg \mathrm{L}_a(b) \in \varGamma,$ for some $a,b$ such that $b \notin X$ and $(\star)$ holds,
\end{itemize}

where the condition $(\star)$ consists of the following:
\begin{enumerate}
\item $\forall c. \ ac \in X$ or $\dot{\neg}(ac) \in X$;
\item $\mathrm{L}_a(b) \notin X$ and $\neg \mathrm{L}_a(b) \notin X.$
\end{enumerate}

Thus, the operator $\Phi(X,\varGamma,\alpha)$ roughly means that $X$ contains the premise of one of the rules of $\sststile{}{}$.
Next, we want to characterise the set of derivable sequents, i.e., the set $\{ \varGamma \ : \ \sststile{}{} \varGamma \}$. Let $\mathrm{Fun}(f)$ be a $\Delta$-predicate meaning that $f$ is a function; 
let a unary $\Delta$-predicate $\mathrm{On}(x)$ express that $x$ is an ordinal number; we use $\alpha, \beta, \gamma, \alpha_{0}, \beta_{0}, \gamma_{0}, \dots$ as variables ranging over $\mathrm{On}$; a $\Sigma$-operation $\mathrm{Dom}(f)$ denotes the domain of $f$.
Then, we define a $\Delta_0$-predicate $\mathscr{H}(s, f),$ which roughly means that for each $\alpha,\beta,\gamma \in s$,
the value of $f$ at $(\alpha,\beta,\gamma)$ is the set $\{ \varGamma \ : \ \sststile{}{\alpha,\beta,\gamma} \varGamma \}$.
%$f$ is an inductive construction of $\sststile{}{<\delta,<\delta,<\delta}$ for a sufficiently large $\delta.$
\[
\mathscr{H}(s, f) := \mathrm{Ad}(s) \land \mathrm{Fun}(f) \land \forall \alpha, \beta,\gamma \in s. \ f(\alpha,\beta, \gamma)=S,
\]
where the set $S$ is the union of the following:
\begin{enumerate}
\item $f(<\alpha,\in s, \in s) \cup f(\alpha, \leq\beta,\leq\gamma),$
\item $\{ \varGamma \ : \ \Phi_0(f(\alpha,\beta,<\gamma), \varGamma) \lor \Phi_1(f(\alpha,<\beta,<\gamma), \varGamma) \},$
\item
$\{ \varGamma \ : \  \alpha \in \mathrm{Suc} \land \Phi_2(f(\alpha -1, \in s, \in s) ,\varGamma) \}.$
\end{enumerate}
Here, we used the following notations:
\begin{itemize}
\item $f(<\alpha, \in s, \in s) := \bigcup_{\alpha_0<\alpha, \beta \in s, \gamma \in s}f(\alpha_0,\beta,\gamma)$;
\item $f(\alpha,\leq\beta,\leq\gamma) := \bigcup_{\beta_0 \leq \beta,\gamma_0\leq\gamma}f(\alpha,\beta_0,\gamma_0)$;
\item $f(\alpha,\beta,<\gamma)$ and $f(\alpha,<\beta,<\gamma)$ are similarly defined;
\item $Cl_{\Phi_i}(X) : \leftrightarrow \forall x. \ \Phi_i(X ,x) \to x \in X$, for $i \in \{0,1\}$.
\end{itemize}

The next lemma shows that $\mathscr{H}(s,f)$ determines the set $f(\alpha,\beta,\gamma)$ for each $\alpha,\beta,\gamma \in s$ regardless of the particular choice of $s$ and $f$.
\begin{lem}[in $\mathsf{KPi}$]\label{lem:monotonicity_and_invariance_of_H}
Assume $\mathscr{H}(s,f)$ and $\mathscr{H}(s',g)$.
Then the following hold for all $\alpha, \beta, \gamma \in s \cap s'$:
\begin{enumerate}
\item $Cl_{\Phi_0}(f(\alpha,\in s,\in s))$ and $Cl_{\Phi_1}(f(\alpha,\in s,\in s))$,
\item $f(\alpha,\in s, \in s)=g(\alpha,\in s',\in s').$
\item $f(\alpha,\beta,\gamma)=g(\alpha,\beta,\gamma)$,
\end{enumerate}
\end{lem}

\begin{pf}
For the item~1, we show $Cl_{\Phi_0}(f(\alpha,\in s,\in s)).$ Thus, taking any sequent $\varGamma$ and assuming $\Phi_0(f(\alpha,\in s,\in s), \varGamma),$
we show $\varGamma \in f(\alpha,\in s,\in s).$ The proof is divided into cases according to the clauses of $\Phi_0.$ As the crucial case, suppose that a sentence $\forall x. \ A(x)$ is contained in $\varGamma$ and $\forall a \in \mathrm{Term}. \ \varGamma, A(a) \in f(\alpha,\in s,\in s).$ Since $s$ is admissible, by $\Sigma$-reflection within $s$, we can take an ordinal $\delta \in s$ such that $\forall a \in \mathrm{Term}. \ \varGamma, A(a) \in f(\alpha,<\delta,<\delta).$
Thus, by applying $\Phi_0$ we obtain $\varGamma \in f(\alpha,<\delta,\delta),$
and hence $\varGamma \in f(\alpha,\in s,\in s).$ The other cases are similar. Moreover, $Cl_{\Phi_1}(f(\alpha,\in s,\in s))$ is similarly proved.

As to item~2: $f(\alpha,\in s, \in s)=g(\alpha,\in s',\in s')$, we show $f(\alpha,\beta, \gamma) \subseteq g(\alpha, \in s', \in s')$ by induction on $\alpha$, $\beta$ and $\gamma.$
Therefore, assuming $\varGamma \in f(\alpha,\beta, \gamma),$ we have to show $\varGamma \in g(\alpha, \in s', \in s').$ By $\mathscr{H}(s,f),$ the proof is divided by cases according to the construction of $f(\alpha,\beta, \gamma)$. 
If $\varGamma \in f(<\alpha,\in s, \in s)$ or $\varGamma \in f(\alpha,\leq\beta, \leq\gamma)$, then the claim is obvious by the side-induction hypothesis.
If $\Phi_0(f(\alpha,\beta,<\gamma), \varGamma)$, then since $\Phi_0$ is positive, we have $\Phi_0(g(\alpha,\in s',\in s'), \varGamma)$ by the side-induction hypothesis.
As $g(\alpha,\in s',\in s')$ is $\Phi_0$-closed by the item~1, it follows that $\varGamma \in g(\alpha,\in s',\in s'),$ as required. The case $\Phi_1(f(\alpha,\beta,<\gamma), \varGamma)$ is similar.
The last case is where $\alpha \in \mathrm{Suc}$ and $\Phi_2(f(\alpha-1,\in s, \in s), \varGamma).$ Then, the main-induction hypothesis yields $\Phi_2(g(\alpha-1,\in s', \in s'), \varGamma)$, thus $\varGamma \in g(\alpha,\in s', \in s').$
In summary, we have $f(\alpha,\in s, \in s) \subseteq g(\alpha,\in s',\in s').$
The converse direction is similar, and thus the proof of item~3 is complete. 

Item~3: $f(\alpha,\beta,\gamma)=g(\alpha,\beta,\gamma)$ is easily proved by induction on $\alpha$, $\beta$ and $\gamma$, with the help of item~2. 
\qed
\end{pf}

We further introduce the following notations:
\begin{itemize}
\item $I^{\alpha,\beta,\gamma}(x) := \exists s, f. \ \{ \alpha,\beta,\gamma \} \subseteq s \land \mathscr{H}(s, f) \land x \in f (\alpha,\beta,\gamma),$
\item $I^{\alpha, \leq\beta,<\gamma}(x) := \exists \beta_{0} \leq \beta. \ \exists \gamma_0 < \gamma. \  I^{\alpha, \beta_0, \gamma_0}(x)$,
\item $I^{\alpha}(x) := \exists \beta, \gamma. \ I^{\alpha,\beta,\gamma}$,
\item $I^{<\alpha}(x)$, $I^{\alpha, \leq\beta,\leq\gamma}(x)$, $I^{\alpha, \beta,<\gamma}(x)$ and $I^{\alpha, <\beta,<\gamma}(x)$ are similarly defined. 
\end{itemize}
We now show that the $\Sigma$-predicate $I^{\alpha,\beta,\gamma}(x)$ expresses the required class $\{ \varGamma \ : \ \sststile{}{\alpha,\beta,\gamma}\varGamma \}$.

\begin{lem}\label{lem:I_is_a_required_derivability}
The following are derivable in $\mathsf{KPi}.$
\begin{enumerate}
\item $\forall \alpha. \ \exists s, f. \  \alpha \in s \land \mathscr{H}(s, f).$

\item $I^{\alpha, \beta, \gamma}(\varGamma)$ if and only if one of the following holds:
\begin{enumerate}
\item $I^{<\alpha}(\varGamma) \lor I^{\alpha,\leq\beta,\leq\gamma}(\varGamma),$
\item $\Phi_{0}(I^{\alpha, \beta, <\gamma}, \varGamma) \lor \Phi_{1}(I^{\alpha, <\beta, <\gamma}, \varGamma),$
\item $\alpha \in \mathrm{Suc} \land \Phi_2(I^{\alpha-1}, \varGamma).$
\end{enumerate}
\end{enumerate}
\end{lem}

\begin{pf}
\begin{enumerate}
\item By the limit axiom of $\mathsf{KPi}$, let $s$ be an admissible set that contains $\alpha$.
Then, from the definition of the $\Delta_0$-predicate $\mathscr{H}(s,f)$, we can construct a required function $f$ by $\Delta$-recursion, available in $\mathsf{KPi}$ (cf. \cite[p.~256]{pohlers2008proof}).

\item For the left-to-right direction, we assume $I^{\alpha, \beta, \gamma}(\varGamma)$; thus, we take sets $s$ and $f$ such that $\{ \alpha,\beta,\gamma \} \subseteq s \land \mathscr{H}(s, f) \land \varGamma \in f (\alpha,\beta,\gamma)$.
For instance, we consider the case where $\varGamma \in f (\alpha,\beta,\gamma)$ is obtained from $\Phi_2$, then
we have $\alpha \in \mathrm{suc}$ and $\Phi_2(f(\alpha-1, \in s, \in s), \varGamma).$
By Lemma~\ref{lem:monotonicity_and_invariance_of_H}, we can easily show that $\forall x. \ x \in f(\alpha-1, \in s, \in s) \leftrightarrow I^{\alpha-1}(x).$
Therefore, we get $\Phi_2(I^{\alpha-1}, \varGamma)$, as required.
The other cases are similarly proved by using Lemma~\ref{lem:monotonicity_and_invariance_of_H}.

As for the converse direction, we, for example, assume $\alpha \in \mathrm{Suc} \land \Phi_2(I^{\alpha-1}, \varGamma)$. By item~1, we take sets $s$ and $f$ such that $\max(\alpha, \beta, \gamma) \in s$ and $\mathscr{H}(s,f).$ Then, again by  Lemma~\ref{lem:monotonicity_and_invariance_of_H}, we have $\forall x. \ I^{\alpha-1}(x) \leftrightarrow x \in f(\alpha-1, \in s, \in s)$. Therefore, we obtain $\Phi_2(f(\alpha-1, \in s, \in s), \varGamma)$, and hence it follows that $\varGamma \in f(\alpha,\beta,\gamma).$ Thus, we obtain $I^{\alpha,\beta,\gamma}(\varGamma)$.
The other cases are similarly proved. \qed
\end{enumerate}
\end{pf}

Let $I^{\infty}(x) :\leftrightarrow \exists \alpha. \ I^{\alpha}(x).$
We now define an interpretation $^{+}$ of $\mathsf{VFU}$ into $\mathsf{KPi}$.
The vocabularies of $\mathcal{L}$ are interpreted in exactly the same way as in Lemma~\ref{lem:interpretation_of_KFUPI}.
Then, we let $\mathrm{T}^{+}(x) :\equiv I^{\infty}(x);$ let $\mathrm{U}^{+}(x) :\leftrightarrow \exists a . \ x=\mathsf{l}a \land \forall b. \ I^{\infty}(ab) \lor I^{\infty}(\dot{\neg}(ab)).$

\begin{lem}\label{lem:interpretation_of_VFU_into_KPi}
For each $\mathcal{L}_{FS}$-formula $A$, if $\mathsf{VFU} \vdash A,$ then $\mathsf{KPi} \vdash A^{+}$.
\end{lem}

\begin{pf}
The proof is by directly  running the proof of Theorem~\ref{thm:model_of_VFU} within $\mathsf{KPi}$. 
\begin{description}
\item[$(\mathsf{U} \mhyph \mathsf{Class})$] We want to show $\mathsf{KPi} \vdash (\mathrm{U}(u) \to \mathrm{C}(u))^+$. Thus, taking any term $f$ such that $u = \mathsf{l}f$ and $\forall a. \ I^{\infty}(fa) \lor I^{\infty}(\dot{\neg}(fa))$, we show $\forall a. \ I^{\infty}(\mathsf{l}fa) \lor I^{\infty}(\dot{\neg}(\mathsf{l}fa))$. Then, since $I^{\infty}$ is a $\Sigma$-predicate, we can, by $\Sigma$-reflection, take the least successor ordinal $\alpha$ such that $\forall a. \ I^{<\alpha}(fa) \lor I^{<\alpha}(\dot{\neg}(fa))$.
Therefore, Lemma~\ref{lem:I_is_a_required_derivability} implies that $\forall a. \ I^{\alpha}(\mathsf{l}fa) \lor I^{\alpha}(\dot{\neg}(\mathsf{l}fa))$,
and thus we obtain $\forall a. \ I^{\infty}(\mathsf{l}fa) \lor I^{\infty}(\dot{\neg}(\mathsf{l}fa))$.

\item[$(\mathsf{U} \mhyph \mathsf{True})$] Similar to the above.

\item[$(\mathsf{U} \mhyph \mathsf{Out})$] 
Since $\mathsf{KPi} \vdash (\mathsf{U} \mhyph \mathsf{True})^+$, it suffices to verify $(\mathsf{T} \mhyph \mathsf{Out})^+$, that is, $\mathsf{KPi} \vdash I^{\infty}(\dot{A}) \to A^+$ for each $\mathcal{L}_{FS}$-formula $A$.
For that purpose, we formalise Lemma~\ref{lem_3:model_of_VFU} within $\mathsf{KPi}$, similarly to \cite[Lemma~5.8.2]{cantini1990theory}.
In particular, we can show the following for each natural number $k$:
\[
\mathsf{KPi} \vdash \forall \varGamma \in \mathrm{Seq}^k. \  I^{\infty}(\varGamma) \to \exists x \in \varGamma. \ \mathrm{T}_k(x),
\]
where $\varGamma \in \mathrm{Seq}^k$ means every sentence in $\varGamma$ has the logical complexity $\leq k$; the predicate $\mathrm{T}_k(x)$ is a partial truth predicate such that $\mathsf{KPi} \vdash \mathrm{T}_k(\dot{A}) \leftrightarrow A^+$ for each $\mathcal{L}_{FS}$-formula $A$ with the logical complexity $\leq k$. Then, we have $\mathsf{KPi} \vdash I^{\infty}(\dot{A}) \to A^+$ for each $\mathcal{L}_{FS}$-formula $A$, as required.
\end{description}
The other cases are similarly proved by using Lemma~\ref{lem:I_is_a_required_derivability}. \qed
\end{pf}

Combining Theorem~\ref{thm:lower-bound_of_VFU} with Lemma~\ref{lem:interpretation_of_VFU_into_KPi}, we obtain the proof-theoretic strength of $\mathsf{VFU}$: 
\begin{thm}\label{thm:strength_of_VFU}
$\mathsf{VFU}$ and $\mathsf{T}_0$ are proof-theoretically equivalent.
\end{thm}

%%%%%%%%%%%%%%%%%%%%%%%%%%%%%%%%%%%%%%%%%%%%%%%%%%%%%%%%%%%%%%%%%%%%%%%%%%%%%%%%%%%%%%

\section{Conclusion}
The results of this paper are summarised as follows.
\begin{concl}
All the following theories are proof-theoretically equivalent to $\mathsf{T}_0$:
\begin{itemize}
\item $\mathsf{KFUPI}$ and $\mathsf{PTUPI}$,
\item $\mathsf{KFLU}$ and $\mathsf{PTLU}$,
\item $\mathsf{VFU}$.
\end{itemize}
\end{concl}
The author suggests two directions for future studies. First, given that most of the truth theories have been studied over Peano arithmetic, it would be desirable to find systems over Peano arithmetic that corrrespond to our theories. 
Second, we can consider extending our systems further by stronger universe-generating axioms. 
One example by Cantini~\cite{cantini1996logical} is the Mahlo principle, which is an analogue of the recursively Mahlo axiom in Kripke--Platek set theory.
Therefore, the question is how strong the systems of Frege structure will be by the addition of such a principle.
Furthermore, J{\"a}ger and Strahm \cite{jager2005reflections} formulated an even stronger principle in explicit mathematics; thus, it may be possible to give   its counterpart in the framework of Frege structure.
Of course, philosophical discussions would also be required on how well these are motivated as truth-theoretic principles.


  \bibliographystyle{plain}
  \bibliography{inyou}


%\begin{thebibliography}


%\bibitem{}

%\end{thebibliography}
\end{document}